\documentclass[11pt]{amsart}
\usepackage{amsmath,amsfonts,amsthm,amscd,amssymb,mathrsfs,amssymb,bm}
\usepackage{mathrsfs}
\usepackage{graphicx}
\usepackage{wrapfig}
\usepackage[all]{xy}
\newtheorem{theorem}{Theorem}

\newtheorem{proposition}[theorem]{Proposition}
\newtheorem{lemma}[theorem]{Lemma}
\newtheorem{definition}[theorem]{Definition}
\newtheorem{conjecture}[theorem]{Conjecture}

\newtheorem{remark}[theorem]{Remark}
\newcommand{\aaa}{\alpha}
\newcommand{\bbb}{\beta}

\newcommand{\lmd}{\lambda}
\newcommand{\Lmd}{\Lambda}

\newcommand{\CP}{\mathbb{CP}}
\newcommand{\CC}{\mathbb{C}}

\newcommand{\RR}{\mathbb{R}}
\newcommand{\ZZ}{\mathbb{Z}}
\newcommand{\FF}{\mathbb{F}}

\newcommand{\U}{{\rm{U}}}

\newcommand{\ol}{\overline}
\newcommand{\lra}{\longrightarrow}

\newcommand{\set}{\,|\,}
\newcommand{\proofend}{\hfill$\square$}

\newcommand{\Aut}{{\rm{Aut}}}
\newcommand{\Bs}{{\rm{Bs}}}

\newcommand{\ms}{\mathscr}
\newcommand{\vsp}{\vspace{3mm}}
\newcommand{\hsp}{\hspace{1.7mm}}
\numberwithin{equation}{section}
\numberwithin{theorem}{section}

\begin{document}
\bibliographystyle{alpha} 
\title{Moishezon twistor spaces on $4\CP^2$}
\author{Nobuhiro Honda}
\address{Mathematical Institute, Tohoku University,
Sendai, Miyagi, Japan}
\email{honda@math.tohoku.ac.jp}
\thanks{The author was partially supported by the Grant-in-Aid for Young Scientists  (B), The Ministry of Education, Culture, Sports, Science and Technology, Japan. }
\begin{abstract}
In this paper we classify all Moishezon twistor spaces
on $4\CP^2$.
The classification is given in terms of the structure
of the anticanonical system of the twistor spaces.
We show that the anticanonical map satisfies one of the 
following three properties: (a) birational over the image,
(b) two to one over the image, or
(c) the image is two-dimensional.
We determine structure of the images for each case
in explicit forms.
Then we intensively investigate structure of 
the twistor spaces in the case (b), and 
determine the defining equation of the branch divisor 
of the anticanonical map.
\end{abstract}
\maketitle
\setcounter{tocdepth}{1}
\vspace{-5mm}
\tableofcontents


\section{Introduction}
According to Taubes' theorem \cite{Ta} there exist a huge number of compact oriented 4-manifolds
which admit a self-dual conformal structure.
Associated to any self-dual structure is a 3-dimensional complex manifold of a special kind,
which is so called the {\em twistor space},
but
by a theorem of Hitchin \cite{Hi81} a  compact twistor space does not admit  a K\"ahler metric except
 two 
standard examples.
Also, by a theorem of Campana \cite{C91},  a compact twistor space can be Moishezon only when the 4-manifold
is $S^4$ or $n\CP^2$, the connected sum of $n$ copies of complex projective planes.
By Kuiper \cite{Kui}, the standard metric on $S^4$ is the unique self-dual  structure on $S^4$.
Similarly, by Poon \cite{P86}, the Fubini-Study metric on $\CP^2$ is the unique self-dual structure on $\CP^2$ whose twistor space is Moishezon (or, 
whose self-dual structure is of positive scalar 
curvature).
On $2\CP^2$, in the same paper, Poon constructed a family of Moishezon twistor spaces parametrized by an open interval in $\mathbb  R$, and showed that these are
all Moishezon twistor spaces on $2\CP^2$.
Thus the classification of Moishezon twistor spaces is completely over up to $2\CP^2$.

For the case of $3\CP^2$, by the works of Poon \cite{P92} and Kreussler-Kurke \cite{KK92} a classification was given by means of the complete linear system $|K^{-1/2}|$, where
$K^{-1/2}$ is the natural square root of
the anticanonical line bundle, which is available on any twistor space.
Namely Moishezon twistor spaces on $3\CP^2$ can be classified into two types according to  whether 
$|K^{-1/2}|$ is base point free or not;
if free, the associated map is a degree two morphism onto $\CP^3$ and the branch divisor is a quartic surface,
whereas if not free, the image of the rational map is a non-singular quadric in $\CP^3$,
and the twistor space has to be a LeBrun twistor space constructed in  \cite{LB91}
whose structure is also well understood.

For the case of  $n\CP^2$, $n$ being arbitrary, while a lot of Moishezon twistor spaces are already known,
their classification seems still difficult.
In the first half of this paper
(namely in Sections \ref{s:clsf} and \ref{s:summary}) we provide a classification 
in the case of $4\CP^2$, by means of the anticanonical system
of the twistor spaces.
Because the most important Chern number $c_1^3$
of a twistor space on $n\CP^2$ is given by
$16(4-n)$, the case of $4\CP^2$ is of particular interest. 
We also note that on $4\CP^2$, the system $|K^{-1/2}|$ is not enough for analyzing structure of twistor spaces, since in most cases the system  is  just a pencil.
A simplified version of our classification may be presented as follows:

\begin{theorem}\label{thm:main1}
Let $Z$ be a Moishezon twistor space on $4\CP^2$
and $\Phi$ the anticanonical map of $Z$.
Then exactly one of the following  three situations occurs:
\begin{enumerate}
\item[(a)] $\Phi$ is birational over the image,
\item[(b)] $\Phi$ is (rationally) two to one over a scroll of 2-planes over a conic
(and hence the scroll is in $\mathbb{CP}^4$),
and the branch divisor is an intersection of the scroll with
a quartic hypersurface in $\CP^4$,
\item[
(c)] the image $\Phi(Z)$ is 2-dimensional, and general fibers of \,$\Phi$ \!are non-singular rational curves.
\end{enumerate}
\end{theorem}

Thus the anticanonical map nicely describes the structure of any Moishezon 
twistor spaces on $4\CP^2$.
According to these three structures, let us call  a Moishezon twistor space $Z$ on $4\CP^2$ to be 
{\em birational type, double solid type} or {\em conic bundle type}  if $Z$ belongs to the case (a), (b) or (c) 
respectively. 
Then what we actually do in Section \ref{s:clsf} is to investigate each of these three cases much more in detail and to clarify the structure of the anticanonical map  as follows.

For a twistor space on $4\CP^2$ which is of birational type, we show that the dimension of the anticanonical system is either 8 or 6 as a linear system.
We also show that if the dimension is 8 the anticanonical  image
is a (non-complete) intersection of 10  quadrics in $\CP^8$,
and that if the dimension is 6, the anticanonical image is a complete intersection of three 
quadrics in $\CP^6$ (Theorem \ref{thm:birational}).
The latter spaces are much more general than the former spaces.
Any of these twistor spaces actually exist.

For twistor spaces of  conic bundle type, 
we prove that the dimension of the anticanonical system is 8, 5, or 4.
We also show that if the system is 8-dimensional the anticanonical image is an embedded image of a non-singular quadric  in $\mathbb{CP}^3$ by the Veronese embedding $
\CP^3\subset\CP^9$ induced from the system $|\mathscr O(2)|$,
and that if 5-dimensional the anticanonical image is the Veronese surface $\CP^2$ in $\CP^5$.
We also prove that if the anticanonical system is 4-dimensional the anticanonical image is an intersection of two quadrics in $\CP^4$, and determine their defining quadratic polynomials in explicit forms
(Theorem \ref{thm:conicbundle}).
Thus the structure of the anticanonical image is  determined in explicit forms.
Again we remark that any of these three kinds of twistor spaces actually exist.

For Moishezon twistor spaces on $4\CP^2$ of double solid type,
we classify them into four sub-types I, II, III and IV,
according to the number of the irreducible components
of the base curve of the system $|F|$,
or equivalently, according to the difference of the singularities
of the branch divisor of the anticanonical map.
(These differences are illustrated in Figure \ref{fig:quartics}.)
We mention that type I is most general and  type IV is most special.

After these classifications,
we compute in Section \ref{ss:moduli} the dimension of the moduli spaces
of Moishezon twistor spaces on $4\CP^2$.
In particular, we see that the moduli space takes the highest dimension
for twistor spaces investigated by Campana-Kreussler \cite{CK98} and
double solid twistor spaces of type I.
It is  likely that any two Moishezon twistor spaces are deformation equivalent,
even if we require Moishezon property
to be preserved in the deformation family.

Comparing the classification with the case of $3\CP^2$,
it seems natural to view the twistor spaces of double solid type on $4\CP^2$
 as a right generalization of those on $3\CP^2$ 
having the structure of double covering of $\CP^3$.
(Indeed the name of double solid type is taken from this.)
Since the double covering is mostly determined from the branch divisor,
the main issue is to determine the form of the defining polynomial of the quartic hypersurface in
$\CP^4$ which cuts out the branch divisor as the intersection with 
the scroll, and this is the main topic in the latter half of this paper.
However, there seems to be no easy way for determining this quartic polynomial.
In Section \ref{s:acs}, for this goal we analyze the structure of the anticanonical map for the
twistor spaces of double solid type in detail.
In particular we give a sequence of explicit blowups which eliminates
the base locus of the anticanonical system of the twistor space.
All the operations are presented in figures.

In Section \ref{s:de} we determine the quartic polynomial whose zeros cut out the
branch divisor of the anticanonical map.
In Section \ref{ss:red} we find reducible members of the anticanonical system,
which consist of two irreducible components.
Then  in Section \ref{ss:double} by using these reducible members
as well as the results in Section \ref{s:acs}
we find some particular curves on the branch divisor,
which bring a  strong constraint for the form of the
defining equation of the quartic hypersurface.
Next in Section \ref{ss:Q} we show that 
there exists of a quadric
which contains all these curves on 
the branch divisor.
The defining equation of this quadric  will be
involved in the equation of the quartic hypersurface.
Finally in Section \ref{ss:de} by assembling all these results we
obtain the defining polynomial of the 
quartic hypersurface
(Theorem \ref{thm:main}).

This paper is a combination of the author's three articles
arXiv:1009.3153, 1108.1443, and 1109.5427,
minus some detailed construction of the twistor spaces.
I would like to thank Masaharu Ishikawa for a helpful discussion on
topology of varieties with singularities.
Also I am grateful to the referee for giving me an opportunity to
combine the three articles, which should certainly be much more suitable for presentation.

\vsp
\noindent
{\bf Notations.}
The natural square root  $K^{-1/2}$ of the anticanonical line bundle, 
which is available on any twistor space, is 
denoted by the letter $F$.
This is called the fundamental line bundle.
The degree of a divisor on a compact twistor space
means the intersection number with twistor lines.
For a line bundle $\mathscr L$ on a compact complex manifold, we write $h^i(\mathscr L)= \dim H^i(\mathscr L)$.
The dimension of a linear system $|\mathscr L|$ always refers $h^0(\mathscr L) -1 $.
${\rm{Bs}} |\mathscr L|$ denotes the base locus of the  linear system $|\mathscr L|$.
For a non-zero element $s\in H^0(\mathscr L)$, $(s)$ means the zero divisor of $s$.
For a non-negative integer $k$, $\FF_k$ denotes
the ruled surface $\mathbb P(\ms O(k)\oplus \ms O)$ over $\CP^1$.
In particular $\FF_0$ means $\CP^1\times\CP^1$.
 An $(a,b)$-curve means a curve on $\FF_0$ whose bidegree is $(a,b)$.
A divisor and the associated line bundle
are often written by the same letter.

\section{Classification of Moishezon twistor spaces on $4\CP^2$}
\label{s:clsf}
\subsection{The base locus of the fundamental system when $h^0(F)=2$.}\label{ss:bc}
Let $Z$ be a (not necessarily Moishezon) twistor space on $4\mathbb{CP}^2$.
We suppose that the corresponding self-dual structure on $4\CP^2$ is of {\em positive type}\, in the sense that
the conformal structure is
represented by a Riemannian metric 
whose scalar curvature is positive.
 (Note that by a theorem of Poon \cite{P88} this positivity always holds if $Z$ is Moishezon.)
Then by the Hitchin's vanishing theorem \cite{Hi80} and the Riemann-Roch formula for the
line bundle $F$, we have
\begin{align}\label{RR1}
h^0(F) - h^1(F) = 2.
\end{align}
In particular, we always have $h^0(F)\ge 2$.
The structure of Moishezon twistor spaces satisfying the inequality $h^0(F)>2$ is well understood by the following result of Kreussler \cite{Kr98}:
\begin{proposition}\label{prop:Kr}
Let $Z$ be a  twistor space on $4\CP^2$ of positive type as above.
\begin{enumerate}
\item[ (i)]
If $h^0(F)\ge 4$, then $h^0(F)= 4$ and $Z$ is a LeBrun twistor space constructed in \cite{LB91}.
\item[(ii)] If $h^0(F)=3$, then $Z$ is Moishezon if and only if $\Bs\,|F|\neq\emptyset$.
Moreover, in this situation, $Z$ has to be a twistor space studied by Campana and Kreussler in \cite{CK98}.
\end{enumerate}
\end{proposition}

Recall that LeBrun twistor spaces are
Moishezon  satisfying $h^0(F)=4$,
and that the twistor spaces 
studied by Campana-Kreussler are 
also Moishezon and satisfy $h^0(F)=3$.

Because generic twistor spaces on $4\mathbb{CP}^2$
satisfy $h^0(F)=2$, the twistor spaces in Proposition \ref{prop:Kr}
 should be  quite special ones among all Moishezon twistor spaces
on $4\CP^2$.
In order to classify all Moishezon twistor spaces on $4\CP^2$, 
it remains to investigate those which satisfy
$h^0(F)=2$.
Let $Z$ be such a twistor space.
Then because $Z$ is not a LeBrun twistor space,  by Poon \cite[Theorem 3.1]{P92}, $Z$ does not have a pencil of degree-one divisors.
Hence  general members of the pencil $|F|$ are irreducible.
Take any real irreducible member $S$ of this pencil.
Then by a theorem of Pedersen-Poon \cite{PP94}, $S$ is a non-singular rational surface satisfying 
$K_S^2=0$.
For this surface we have the following.
(See also \cite{Kr98}.)

\begin{proposition}
\label{prop:cycle}
Let $Z$ be a Moishezon twistor space on $4\mathbb{CP}^2$ satisfying $h^0(F)=2$,
and  $S$  any  real irreducible member of the pencil $|F|$ as above.
Then $h^0(S,K_S^{-1}) = 1$.
Let $C$ be the unique anticanonical curve on $S$.
Then $\Bs\,|F| = C$ and $C$ is a cycle of non-singular rational curves. 
Moreover, the number $m$ of the irreducible components of $C$ is even with $4\le m\le 12$.
\end{proposition}

Here, by a {\em cycle} of non-singular rational curves, we mean a connected reduced 
divisor on a non-singular surface, which is of the form
$ \sum_{i=1}^m C_i$ with $ m\ge 2$
where $C_1,\dots,C_m$ are mutually distinct smooth rational curves which satisfy
\begin{itemize}
\item
when $m\ge 3$, $(C_i, C_{i+1})_S=1$ holds for any $1\le i\le m$
where $C_{m+1}$ means $C_1$,
\item when $m=2$, $C_1$ and $C_2$
intersect transversally at two points.

\end{itemize}

\vsp
\noindent {\em Proof of Proposition \ref{prop:cycle}.}
By \eqref{RR1} and the assumption
$h^0(F)=2$ we have $H^1(F)=0$.
So the exact sequence 
$ 0 \to \mathscr O_Z\stackrel{\otimes s}{\to} F \to F|_S\simeq K_S^{-1} \to 0$,
where $s$ being a section of $F$ satisfying $(s)=S$, gives an exact sequence
\begin{align}
 0 \lra \mathbb C \stackrel{s\otimes}{\lra} H^0(F) \lra H^0(K_S^{-1}) \lra 0.
\end{align}
This means $h^0(K_S^{-1}) = 1$ and  $\Bs\,|F| = C$.
Let $\epsilon:S\to \FF_0$ be 
a birational morphism preserving the real structure which maps twistor lines in $S$ to   $(1,0)$-curves (see \cite{PP94}). 
Then since the image $\epsilon_*(C)$ (= the image which takes the multiplicity into account) is also an anticanonical curve, 
it is a $(2,2)$-curve on $\FF_0$, which is real.
From the absence of real points on $S$ and $\FF_0$,
we can readily deduce that
 if $\epsilon_*(C)$ has a non-reduced component,
the component must be the twice of a real $(1,0)$-curve.
But if this is really the case,
$\epsilon$ always blows up points on the real $(1,0)$-curve,
which readily means  $h^0(K_S^{-1})=3$. 
This contradicts $h^0(K_S^{-1})=1$. 
Therefore $\epsilon_*(C)$ does not have a non-reduced component.
Hence possible forms of the $(2,2)$-curve $\epsilon_*(C)$ are the following three situations:
\begin{eqnarray}
C_1 + \ol C_1, 
&
C_1\in|\mathscr O(1,1)|, \label{001}
\\
C_1+C_2,  
&
 C_1\in|\mathscr O(1,2)|,
\, C_2 \in|\mathscr O(1,0)|, \label{002}
\\
C_1 + C_2 + \ol C_1 + \ol C_2, 
&
C_1\in |\mathscr O(1,0)|, C_2\in |\mathscr O(0,1)|,
\label{003}
\end{eqnarray}
where in each possibility two components of the same bidegree
are distinct.
But if \eqref{002} would be the case,
since $\epsilon$ cannot blowup points on
the $(1,0)$-curve $C_2$
which is necessarily transformed into
a twistor line on $S$,
we would have $h^0(K_S^{-1}) \ge h^0(\mathscr O(1,0))=2$.
Hence the  situation \eqref{002} cannot happen under our
hypothesis.
For the  situation \eqref{003}  it is obvious that the image $\epsilon_*(C)=\epsilon(C)$ is a cycle of smooth rational curves.
The same conclusion holds for the first possibility \eqref{001} since the two curves $C_1$ and $\ol C_1$ cannot touch at a point
as  $\FF_0$ does not have a real point. 
Hence in both cases
the image $\epsilon_*(C)=\epsilon(C)$ is a cycle of smooth rational curves.
Moreover since $C\in |K_S^{-1}|$, each step
in the birational morphism $\epsilon:S\to\FF_0$  blows up a point on (the inverse image of) $\epsilon(C)$.
This means that $C$ is also a cycle of smooth rational curves and  that the number $m$ of irreducible components
satisfies $2\le m\le 4 + 8 = 12$.
Obviously $m$ is even by the presence of real structure.
It remains to see $m\neq 2$.
If $m=2$, every steps in $\epsilon$ blow up a smooth point of the anticanonical cycle \eqref{001}.
If we use the same letters to mean
the strict transforms of $C_1$ and $\ol C_1$
into $S$,
this means  that $(C_1)_S^2=(\ol C_1)_S^2 = -2$, which implies that
$(C_1+ \ol C_1, C_1)_S = (C_1+\ol C_1,\ol C_1)_S = 0$.
Hence the restriction of the anticanonical bundle
 $\mathscr O_S(C_1+\ol C_1)$ on $S$ to the cycle $C_1 + \ol C_1$ is topologically trivial. 
From this it readily follows that
as a function of an integer $l$, the dimension $h^0(lK_S^{-1})$ increases at most linearly as $l\to \infty$.
This implies that $h^0(Z, lF)$ 
can increase at most quadratically as $l\to \infty$.
This contradicts our assumption that $Z$ is Moishezon.
Hence $m\neq 2$, as claimed.
\proofend

\vsp
By Proposition \ref{prop:cycle} Moishezon twistor spaces on $4\CP^2$ satisfying $h^0(F)=2$ can be classified 
in terms of the number $m$ of  irreducible components of the base curve  of the pencil $|F|$. 
In the sequel we  write the number $m$ by $2k$, where $2\le k\le 6$.

Since the cycle $C$ consists of $2k$ irreducible components,
there are exactly $k$ ways  to divide $C$ into connected halves.
We write $L_i$ ($1\le i\le k)$ for the twistor line going through 
the point $C_i\cap C_{i+1}$, where we read $C_{k+1} = \ol C_1$.
Then for each $i$, there exists a 
pair of divisors $\{S_i^+, S_i^-\}$ which satisfies the following 
properties that $$S_i^-=\ol S_i^+, \,
S_i^++S_i^-\in |F| {\text{\, and\, }} 
S_i^+\cap S_i^-=L_i$$
(see \cite{Kr98}).
The obvious decomposition 
$C = (S_i^+\cap C)\cup (S_i^-\cap C)$ realizes
the divisions into halves. 
These are all degree-one divisors on $Z$,
and play an important role for deriving 
the defining equation of the image under
the anticanonical map of the twistor spaces.
We make distinction for the two components $S_i^+$ and $S_i^-$ by 
promising that $S_i^-$ contains the component $C_1$ of $C$.

Let $\epsilon:S\to\FF_0$ be the birational morphism as in the above proof.
Then with an aid of the real structure we can factorize $\epsilon$ as 
\begin{align}\label{dcp1}
\epsilon =\epsilon_4\circ\epsilon_3\circ\epsilon_2\circ\epsilon_1,
\end{align}
where each $\epsilon_i$ blows up a real pair of points on the anticanonical cycle.
Moreover, as we have denied the possibility $m=2$, even when 
$\epsilon(C)$ is as in  \eqref{001},
some  $\epsilon_i$ has to blowup a pair of singular points of the curve
$C_1+\ol C_1$.
From this, by choosing a different blowdown to $\FF_0$,
 it follows that the situation \eqref{001} is absorbed in the situation \eqref{003}.
 Thus for any Moishezon twistor space on $4\CP^2$ with $h^0(F)=2$ 
we can  suppose that the image $\epsilon(C)\subset \FF_0$ is a real $(2,2)$-curve consisting of four irreducible components.
Hence in the following we choose the birational morphism
$\epsilon$ in a way that the image $\epsilon(C)$ is 
a cycle of {\em four}\, rational curves
as in \eqref{003}.

We will often make use of the following easy property for the anticanonical system
$|2F|$ on $Z$, which is valid only on $4\CP^2$:
\begin{proposition}\label{prop:acs}
Let $Z$ be a Moishezon twistor space on $4\CP^2$ which satisfies $h^0(F)=2$.
Then if $S\in |F|$ is a smooth member,  we have 
\begin{align}\label{acs1}
h^0(2F) =  h^0(2K_S^{-1}) + 2.
\end{align}
\end{proposition}

\proof
As $H^1(F)=0$, from the exact sequence $0 \lra F \stackrel{s\otimes}{\lra} 2F \lra 2K_S^{-1}\lra 0$, where $s\in H^0(F)$  satisfies $(s)=S$ as before,
by using $h^0(F)=2$ we get an exact sequence
\begin{align}\label{acs2}
0 \lra H^0(F) \stackrel{s\otimes}{\lra} H^0(2F) \lra H^0(2K_S^{-1}) \lra 0.
\end{align}
This implies \eqref{acs1}.
\proofend

\vsp
Let $S^2H^0(F)$ be the symmetric product of
the 2-dimensional space $H^0(F)$.
As $h^0(F) = 2$, this is a $3$-dimensional subspace of 
$H^0(2F)$.
Then we have the following basic commutative
diagram of rational maps:
\begin{equation}\label{CD1}
 \CD
 Z@>\Phi_{|2F|}>> \mathbb P H^0(2F)^*\\
 @V\Phi_{|F|} VV @VV\pi V\\
 \mathbb PH^0(F)^*@>g >> \mathbb P\, S^2H^0(F)^*.
 \endCD
 \end{equation}
Here, the asterisk means the dual 
vector space,
$\pi$ is the linear projection 
which is the dual of the inclusion $S^2H^0(F)\subset
H^0(2F)$, and $g$ is an embedding 
of $\mathbb PH^0(F)^*=\CP^1$ as a conic
in $\mathbb P\, S^2H^0(F)^*=\CP^2$.
Evidently whole of the diagram \eqref{CD1} is
preserved by the natural real structure 
induced from that on the twistor space $Z$.

\begin{definition}\label{def:scroll}
{\em
Throughout this article we use the following notations
and terminologies.
\begin{enumerate}
\item[(i)] We always use the letter $\Phi$ to mean the anticanonical map from a twistor space.
\item[(ii)]
We denote the image conic of $g$ by the letter $\Lmd$.
\item[(iii)]
When $h^0(2F) = 5$,  we always denote
the pre-image $\pi^{-1}(\Lmd)\subset\CP^4$ by $Y$, and 
call it as a {\em scroll of planes over the conic} $\Lmd$.
\item[(iv)]
Under the situation of (iii), we denote 
the indeterminacy locus of the projection $\pi$ by $l$
and call it as the {\em  ridge} of the scroll $Y$.
As $\pi$ is a projection from
$\CP^4$ to $\CP^2$, this is a line in $\CP^4$.
\end{enumerate}}
\end{definition}

From the definition, the conic $\Lmd$ 
can be regarded as the 
parameter space of the pencil $|F|$.
Also, when $h^0(2F)=5$,
the target space of the projection 
$\pi:\CP^4\to\CP^2$ can be 
regarded as the space of 
2-planes containing  the line $l$.
Moreover, from the commutativity of the diagram \eqref{CD1} we always have $\Phi_{|2F|}(Z)\subset
\pi^{-1}(\Lmd)$ for the anticanonical image.
Furthermore, by the exact sequence \eqref{acs2},
for any smooth member $S$ of the pencil
$|F|$, 
the restriction of the anticanonical map 
$\Phi_{|2F|}$ to $S$ can be identified with
the bi-anticanonical map of $S$.
In the rest of this section 
using these as basic tools we investigate structure of the anticanonical map of  $Z$ for each of the five values of $k$.

\subsection{Classification in the case $k=6$}\label{ss:k6}
From the proof of Proposition \ref{prop:cycle},
$k$ attains the maximal
value $6$ iff any $\epsilon_i$ in the factorization \eqref{dcp1}  blows up a (real) pair of singular points of the cycle;
namely exactly when $S$ is a toric surface.

It is elementary to see that without loss of generality
we may suppose that for the image $\epsilon_4(C)$, which is 
necessarily a cycle of 10 smooth rational curves, 
the sequence obtained by putting the self-intersection numbers 
of the components is given by $$-3,-1,-2,-2,-1,-3,-1,-2,-2,-1$$ up to permutations and reversing the order.
Therefore all the freedom for giving
the toric $S$ lie on the choice of a pair of points which are blown-up by  $\epsilon_4$.
It is also elementary to see that there are exactly three
 choices of $\epsilon_4$ which give
mutually non-isomorphic toric surfaces, and 
that the sequence  of the self-intersection numbers of  components
of $C$ for these three choices can be listed as:
\begin{align}
\label{string1}-4,   -1,   -2,   -2,   -2,   -1,   -4,   -1,   -2,   -2,   -2,   -1,\\
\label{string2}-3,   -2,   -1,   -3,   -2,   -1,  -3,   -2,   -1,  -3,   -2,   -1,\\
\label{string3}-3,    -1, -3,    -1,   -3,    -1,   -3,    -1,   -3,    -1,   -3,    -1.   
\end{align}
It is immediate to see that for the case 
 \eqref{string1}, we have $h^0(K_S^{-1})=3$ which 
contradicts our assumption.
Hence it suffices to consider the other two cases \eqref{string2} and \eqref{string3}.
The following properties on these 
two toric surfaces are included in 
\cite[Lemmas 2.4 and 4.1]{Hon_JAG}:

\begin{lemma}
\label{lemma:toric1}
If the sequence for the toric surface $S$ is as in \eqref{string2},  we have
the following.
\begin{enumerate}
\item[(i)] $h^0(2K_S^{-1})= 5$ and\, $\Bs\,|2K_S^{-1}| = C-(C_3+C_6+\ol C_3 + \ol C_6)$.
\item[(ii)] After removing this  curve,
$|2K_S^{-1}|$ is base point free.
\item[(iii)]  
The associated morphism $S\to\CP^4$ is birational over the image.
\end{enumerate} 
\end{lemma}

\begin{lemma}
\label{lemma:toric2}
If the sequence for the toric surface $S$ is as in \eqref{string3},  we have
the following.
\begin{enumerate}
\item[(i)] $h^0(2K_S^{-1})= 7$ and $\Bs\,|2K_S^{-1}| = C_1+C_3+C_5+\ol C_1+\ol C_3 + \ol C_5$.
\item[(ii)] After removing this  curve,
$|2K_S^{-1}|$ is base point free.
\item[(iii)]  
The associated morphism $S\to\CP^6$ is birational over the image.
\end{enumerate} 
\end{lemma}

By using these lemmas, we show the classification in the case $k=6$:

\begin{proposition}\label{prop:k6}
Let $Z$ be a Moishezon twistor space on $4\CP^2$ satisfying $h^0(F)=2$.
Suppose that the anticanonical cycle $C$ of  
a real irreducible member $S\in |F|$
consists of  $12$ irreducible components.
Then one of the following two situations  happens:
\begin{enumerate}
\item[(i)] 
If the self-intersection numbers for the cycle $C$ in the toric surface $S$ are as in \eqref{string2}, then $h^0(2F)= 7$ and  the anticanonical map $\Phi:Z\to \CP^6$ is birational over the image.
Further, the image is a complete intersection of three quadrics, whose defining equations are of the form
\begin{align}\label{J1}
z_0^2 = z_1 z_2,\,\,
z_3z_5=q_1(z_0,z_1,z_2),\,\,
z_4z_6=q_2(z_0,z_1,z_2),
\end{align}
where $z_0,z_1,\cdots,z_6$ are homogeneous coordinates
on $\CP^6$ and $q_1$ and $q_2$ are  
homogeneous quadratic polynomials
in $z_0,z_1,z_2$ with real coefficients.
\item[(ii)]
If the self-intersection numbers for the cycle $C$ are
as in \eqref{string3}, then
 $h^0(2F)=9$, and the anticanonical map $\Phi:Z\to \CP^8$ is birational over the image.
Further, the image satisfies the following 
10 quadratic equations:
\begin{align}\label{J2}
z_0^2 = z_1 z_2,\,\,
z_3z_6=q_1(z_0,z_1,z_2),\,\,
z_4z_7=q_2(z_0,z_1,z_2),\,\,
z_5z_8=q_3(z_0,z_1,z_2),\\
\label{J2'}
z_3z_4 = f_1(z_0,z_1,z_2)z_8,\,\,
z_4z_5 = f_2(z_0,z_1,z_2)z_6,\,\,
z_5z_3 = f_3(z_0,z_1,z_2)z_7,\\
\label{J2''}
z_6z_7 = f_1(z_0,z_1,z_2) z_5,\,\,
z_7z_8 = f_2(z_0,z_1,z_2)z_3,\,\,
z_8z_6 = f_3(z_0,z_1,z_2)z_4.
\end{align}
where all the $q_i$-s and \,$f_i$-s are quadratic and linear polynomials respectively
in $z_0,z_1,z_2$ with real coefficients.
\end{enumerate}
In particular, when $k=6$, the anticanonical map of $Z$ is
always birational over the image.
\end{proposition}

\proof
The assertions for the value of 
$h^0(2F)$ immediately follow from
Proposition \ref{prop:acs} and
Lemmas \ref{lemma:toric1}, \ref{lemma:toric2}.
The birationality of the anticanonical 
map $\Phi$ also follows from 
the two lemmas and the commutative diagram \eqref{CD1}.
So it remains to derive the defining equations of the anticanonical image
$\Phi(Z)$.
For both situations (i) and (ii),
let $u_1,u_2$ be any real basis of $H^0(F)$ and 
put $z_0=u_1u_2,z_1=u_1^2, z_2= u_2^2$.
These are real basis of 
the 3-dimensional subspace $S^2H^0(F)$, and
subject to the relation $z_0^2=z_1z_2$.
The choice of the remaining basis of $H^2(2F)$ depends
on the situation.

We first consider the situation (i).
By the combinatorial method introduced in \cite[Section 2, Procedure (A)]{Hon-III}, 
it is not difficult to see  that the following 
two divisors are members of the anti-canonical system $|2F|$:
\begin{align}\label{redmembers4}
X_1:= S_1^+ + S_2^+ + S_3^- + S_6^-,
\quad
X_2:= S_3^+ + S_4^+ + S_5^+ + S_6^-.
\end{align}
(Although presence of a $\CC^*\times\CC^*$-
or a $\CC^*$-action 
on the twistor spaces is supposed in \cite{Hon-III}, for the purpose of finding members of the multiple system $|mF|$,
we just need an action of these groups on the divisor $S$ and we do {\em not}
need that on the twistor spaces.
We also note that the divisors \eqref{redmembers4} are
identical with those in \cite[Lemma 2.7]{Hon_JAG},
where they were found rather incidentally.)
Next let $z_3\in H^0(2F)$
(resp.\,$z_4\in H^0(2F)$) be an element 
which satisfies $(z_3) = X_1$ 
(resp.\,$(z_4)=X_2$).
We further put $z_5:=\ol{\sigma^*z_3}$ and 
$z_6:=\ol{\sigma^*z_4}$,
where $\sigma$ is the real structure on 
the twistor space $Z$.
Then obviously $(z_5)= \ol X_1$ and 
$(z_6)= \ol X_2$, and 
the two products  $z_3z_5$  
and $z_4z_6$ are non-zero real elements of the symmetric product $S^4H^0(F)\subset H^0(4F)$.
Hence there exist quadratic polynomials
$q_1,q_2\in \RR[z_0,z_1,z_2] $ which satisfy
\begin{align}\label{J11}
z_3z_5 = q_1(z_0,z_1,z_2),\quad
z_4z_6 = q_2(z_0,z_1,z_2).
\end{align}
(From the relation $z_0^2 = z_1 z_2$,
the polynomials $q_1$ and $q_2$ are determined only up to the ideal
$(z_0^2-z_1z_2)$.)
Thus if $W$ denotes the variety in $\CP^6$
 defined by the three equations \eqref{J1},
the inclusion $\Phi(Z)\subset W$ holds.
We can write $\Lmd$ for the conic $z_0^2=z_1z_2$
in $\CP^2$ (see the diagram \eqref{CD1}).
Then for each $\lmd\in \Lmd$
the intersection $W\cap \pi^{-1}(\lmd)$ is 
an intersection of the two quadrics
\eqref{J11} in $\pi^{-1}(\lmd)=\CP^4$.
Hence $W\cap \pi^{-1}(\lmd)$ is always a
quartic surface (in $\CP^4$).
This means $\dim W=3$.
Moreover since both of the 
polynomials $q_1$ and $q_2$ in \eqref{J11}
are non-zero, $W\cap \pi^{-1}(\lmd)$ is
irreducible for general $\lmd\in \Lmd$.
Hence $W$ is irreducible.
Therefore as $\Phi$ is birational,
we must have $\Phi(Z) = W$, as desired.

Next we consider the situation (ii).
For this case, again by 
\cite[Section 2, Procedure (A)]{Hon-III}, 
we deduce that the following three divisors 
belong to $|2F|$:
\begin{align}\label{redmembers5}
X_1:=  S_2^+ + S_3^+ + S_5^+ + S_6^-,
\quad
X_2:= S_1^- + S_3^- + S_4^+ + S_6^+,
\quad
X_3:= S_1^+ + S_2^- + S_4^- + S_5^-.
\end{align}
(Again these are identical with the ones
presented in \cite[Lemma 4.5]{Hon_JAG}.)
Let $z_3,z_4$ and $z_5$ be elements
of $H^0(2F)$ which satisfy $(z_3) = X_1,
(z_4) = X_2$ and $(z_5) = X_3$,
and put 
$z_6:=\ol{\sigma^*z_3},
z_7:=\ol{\sigma^*z_4}$ and 
$z_8:=\ol{\sigma^*z_5}$.
Then by the same reason to the case (i) above,
there exist quadratic polynomials
$q_1,q_2,q_3\in \RR[z_0,z_1,z_2]$ 
which satisfy 
$$
z_3 z_6 = q_1 (z_0,z_1,z_2),\quad
z_4 z_7 = q_2 (z_0,z_1,z_2),\quad
z_5 z_8 = q_3 (z_0,z_1,z_2).
$$ 
Hence the image $\Phi(Z)$ satisfies the four equations
\eqref{J2}.
For deriving the next three equations \eqref{J2'},
we notice the relations
\begin{align}
X_1 + X_2 &= \ol X_3 + (S_3^+ + S_3^-)
+ (S_6^++S_6^-),\label{ntre1}\\
X_2 + X_3 &= \ol X_1 + (S_1^+ + S_1^-)
+ (S_4^++S_4^-),\label{ntre2}\\
X_3 + X_1 &= \ol X_2 + (S_2^+ + S_2^-)
+ (S_5^++S_5^-),\label{ntre3}
\end{align}
which are directly obtained from \eqref{redmembers5}.
From \eqref{ntre1}, 
since $(S_3^+ + S_3^-)
+ (S_6^++S_6^-)\in |S^2H^0(F)|$,
there exists
a linear polynomial $f_1\in \RR[z_0,z_1,z_2]$ 
which satisfies $z_3z_4 = f_1 (z_0,z_1,z_2)z_8.$
Similarly, from \eqref{ntre2} and \eqref{ntre3}
we obtain that 
there exist
linear polynomials $f_2,f_3 \in \RR[z_0,z_1,z_2]$ 
which satisfy $z_4z_5 = f_2(z_0,z_1,z_2)z_6$
and $z_5z_3 = f_3(z_0,z_1,z_2)z_7$.
Hence $\Phi(Z)$ satisfies the three equations
\eqref{J2'}.
By taking  $\ol\sigma^*$ for these, we obtain the final three equations
\eqref{J2''}.
(Note that $q_i$-s in \eqref{J2'}
and \eqref{J2''}
are common since they are of real coefficients.)
Thus we obtain that $\Phi(Z)$ satisfies
all the 10 equations \eqref{J2},
\eqref{J2'} and \eqref{J2''}.
\proofend

\vsp
We remark that those twistor spaces on $4\CP^2$
which have these two kinds of toric surfaces
as members of $|F|$
 were investigated in \cite{Hon_JAG}, under a stronger assumption
that the $\mathbb C^*\times\mathbb C^*$-action on the toric surface
extends to that on the twistor space $Z$.
The above combinatorial proof makes quite simpler
 the argument for obtaining 
defining equations of the anticanonical models
of the twistor spaces of Joyce metrics \cite{J95} given in \cite[Sections 2.2 and  4]{Hon_JAG}.
In view of the similarities
of the defining equations as well as
the fact that 
the anticanonical models actually
have effective $\CC^*\times\CC^*$-action
from the equations \eqref{J1}--\eqref{J2''},
it seems quite likely that
 one can show that  
the twistor spaces in Proposition \ref{prop:k6}
are exactly the twistor spaces of Joyce metrics,
by analyzing the anticanonical map more in detail.
But here we do not try to do so because 
such an identification seems to require 
considerable efforts.

\subsection{Classification in the case $k=5$}\label{ss:k5}
As before let $Z$ be a Moishezon twistor
space on $4\CP^2$ which satisfies 
$h^0(F) = 2$, and $S\in |F|$ 
a real irreducible member.
Suppose that the number $m=2k$ of the 
irreducible components of the 
anticanonical cycle $C$ on $S$ is $10$.

Let $\epsilon:S\to\FF_0$ be the birational morphism as in Section \ref{ss:bc} and again decompose it as $\epsilon =\epsilon_4\circ\epsilon_3\circ\epsilon_2\circ\epsilon_1$
as in \eqref{dcp1},
where each $\epsilon_i$ blows up a real pair of points belonging to the anticanonical cycle.
Further, as already mentioned
just before Proposition \ref{prop:acs}, we can suppose that the image $\epsilon(C)\subset\FF_0$ is a cycle of four rational curves (as in \eqref{003}). 
Then from the assumption $k=5$, 
 there exists a unique index $j$ such that $\epsilon_j$ blows up a pair of {\em smooth} points  on the anticanonical cycle.
Without loss of generality we may
suppose that $j=4$.
Under this situation we put $S_1:=\epsilon_4(S)$.
By the choice $S_1$ has an anticanonical cycle which consists of 10 irreducible components, and it is a toric surface.
We note that regardless of the blowup points of $\epsilon_4$,
the resulting surface $S$ admits a non-trivial $\CC^*$-action.
As mentioned at the beginning of Section \ref{ss:k6},  possible structure of the toric surface $S_1$ is unique,
and the self-intersection numbers
of the 10 components of the anticanonical cycle
$\epsilon_4(C)\subset S_1$ in order are given by,
up to permutations and reversing the order,
\begin{align}\label{string4}
-3,-1,-2,-2,-1,-3,-1,-2,-2,-1.
\end{align}
Let $C_1,\cdots,C_5,\ol C_1,\cdots,\ol C_5$ be the  components
of $\epsilon_4(C)$ arranged in this order (so  $C_1^2=\ol C_1^2=-3$ for instance.)
Note that th toric surface $S_1$ has a holomorphic involution which commutes with the real structure and 
which exchanges $C_1,C_2,C_3,C_4$ and $ C_5$
with $\ol C_1, C_5,C_4,C_3$ and $C_2$ respectively.
Thanks to this symmetry we can suppose that $S$ is obtained from $S_1$ by blowing-up a point on $C_1,C_2$ or $C_3$ and its conjugate point,
which are not on the singularities of the 
cycle $\epsilon_4(C)$. 
But if  $\epsilon_4$ blows up a point on $C_1$, then $S$ satisfies $h^0(K_S^{-1})=3$ and this contradicts  $h^0(K_S^{-1})=1$.
So $\epsilon_4$ blows up a point of $C_2\cup C_3$, which is not on the 
singularities of the anticanonical cycle.

\begin{figure}
\includegraphics{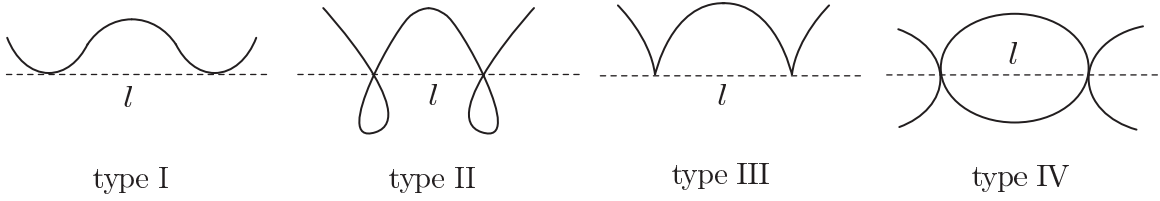}
\caption{The branch quartic curves of the 
bi-anticanonical map $\phi:S\to\CP^2$
which is of degree two, in relation
with the special line $l=\phi(C)$.
See Definition \ref{def:DS} for the types.}
\label{fig:quartics}
\end{figure}

Suppose that $\epsilon_4$ blows up a  point of $C_2$.
Then the structure of $S$ is the same as that of a real irreducible member of $|F|$
in the twistor spaces studied in \cite{HonDSn1} (specialized to the case of $4\CP^2$, of course).
In particular from \eqref{string4} the sequence of the self-intersection numbers of the components 
of $C$ becomes
\begin{align}\label{string5}
-3, -2, -2, -2, -1, -3, -2, -2, -2, -1.
\end{align}
Hence by \cite[Proposition 2.1]{HonDSn1} we have
\begin{lemma}\label{lemma:SIV}
The bi-anticanonical system of this surface $S$
satisfies the following.
\begin{enumerate}
\item[(i)]
$h^0(2K_S^{-1}) = 3$.
\item[(ii)] $\Bs\,|2K_S^{-1}| = \sum_{i=1}^4(C_i + \ol C_i)$,
and after removing this base curve, the bi-anticanonical system is base point free.
\item[(iii)] The  morphism $\phi:S\to\CP^2$ induced by 
the bi-anticanonical system is of degree two, and 
the branch divisor is a reducible  quartic curve.
\item[(iv)] The morphism $\phi$ contracts the connected curve
$\ol C_5\cup C_1\cup C_2\cup C_3$ and the conjugation curve
to mutually distinct points.
\item[(v)] The morphism $\phi$ maps $C_4$ and $\ol C_4$ to
an identical line isomorphically.
\item[(vi)] The branch quartic curve has $A_3$-singularities
at the two points which are the images of the two curves in (iv).
(See Figure \ref{fig:quartics}, type IV.)
\end{enumerate}
\end{lemma}

Next we suppose that $\epsilon_4$ blows up a  point of $C_3$.
In this case, again from \eqref{string4}, the sequence of the self-intersection numbers of the anticanonical cycle 
on $S$ becomes
\begin{align}\label{string6}
-3,-1,-3,-2,-1,-3,-1,-3,-2,-1.
\end{align}
For this surface $S$ we have the following

\begin{lemma}\label{lemma:83}
The bi-anticanonical system of $S$ satisfies the following.
\begin{enumerate}
\item
[(i)] $h^0(2K_S^{-1})=5$.
\item
[(ii)]   $\Bs\,|2K_S^{-1}| = C_1 + C_3 + C_4 
+ \ol C_1 + \ol C_3 + \ol C_4$, and after removing this,
the system becomes base point free.
\item
[(iii)] If $\phi:S\to \CP^4$ is the
morphism induced by the bi-anticanonical system, then $\phi$ is birational over the image  $\phi(S)$, and 
$\phi(S)$
is an intersection of two 
quadrics in $\CP^4$.
\end{enumerate}
\end{lemma}

\proof
In this proof for simplicity we write $K^{-1}$
to mean $K_S^{-1}$.
The curves $C_1, C_3, C_4$ and their conjugations can be readily seen to be 
fixed components of $|2K^{-1}|$ by computing the intersection numbers.
Let $B$ be the sum of these 6 curves (with each coefficient being one).
Then we easily have $(2K^{-1}-B,\,C_i)=(2K^{-1} - B, \,\ol C_i)= 0$
for $i=2,5$.
Hence noticing $(2K^{-1}-B) -( C_2 +C_5 + \ol C_2 + \ol C_5 ) \simeq K^{-1}$,
 we obtain an exact sequence
\begin{align}\label{ses:25}
0 \lra K^{-1} \lra 2K^{-1}
-B \lra 
\mathscr O_{C_2 \,\cup\, C_5 \,\cup\, \ol C_2 \,\cup\, \ol C_5} \lra 0,
\end{align}
where the second arrow is a multiplication of a defining section of the divisor
$C_2 +C_5 + \ol C_2 + \ol C_5$.
Further we readily have $h^0(K^{-1})=1$ and hence $h^1(K^{-1})=0$ by Riemann-Roch.
Therefore from the cohomology exact sequence of \eqref{ses:25} we obtain $h^0(2K^{-1}) = 1+4=5$,
and also that ${\rm{Bs}}\, |2K^{-1} - B|$ is disjoint from $C_2 \cup C_5 \cup \ol C_2 \cup \ol C_5$. 
Hence, as $(2K^{-1} - B,\, C_3) =0$, 
 ${\rm{Bs}}\, |2K^{-1} - B|$ is disjoint from $C_3\cup \ol C_3$ also.
 Therefore  if ${\rm{Bs}}\, |2K^{-1} - B|$ is non-empty, it must be points on 
 $C_1\cup C_4 \cup \ol C_1 \cup \ol C_4$ which are not singularities of the anticanonical cycle $C$.
 But this cannot happen since $S$ has a non-trivial $\mathbb C^*$-action which 
does not fix $C_1\cup C_4\cup \ol C_1\cup
\ol C_4$, and this $\CC^*$-action automatically lifts on the system $ |2K^{-1} - B|$,
 so that its base locus must be $\CC^*$-invariant.
 Hence $\Bs\,|2K^{-1}-B| = \emptyset$,
meaning (ii).
 
 To show (iii), we notice that the two divisors
 $$
 f:= C_4 + 2 C_5 + \ol C_1 + \ol C_2 
 \hsp{\text{and}}\hsp
\ol f:=\ol  C_4 + 2 \ol  C_5 +  C_1 +  C_2
 $$
on $S$ are mutually linearly equivalent, and that $2C-f+\ol f$ and $2C + f - \ol f$ are effective divisors,
 where $C$ is the unique anticanonical curve (cycle) on $S$ as before. 
 Therefore we can consider the linear subsystem of $|2K^{-1}|$ generated by the three curves 
 \begin{align}
 \label{3curves}
 2C-f+\ol f,\,\, 2C + f - \ol f
  \hsp{\text{and}}\hsp
  2C.
 \end{align}
By looking at the multiplicities of the components, we readily see that 
 these three divisors are linearly independent.
Hence these generate a 2-dimensional subsystem
of $|2K^{-1}|$.
Write $|V|$ for this subsystem.
If $\phi:S\to \CP^4$ denotes the bi-anticanonical map and 
$p:\CP^4\to \CP^2$ means the projection
corresponding to the inclusion $V\subset H^0(2K^{-1})$, 
the rational map $S\to \CP^2$ induced by the subsystem
$|V|$ is naturally factorized as 
$p\circ\phi$.

On the other hand, from the 
generators \eqref{3curves}, we readily see
that the fixed components of the subsystem
$|V|$ are the curve $B + C_2 + C_3 + \ol C_2 + \ol C_3$. 
Removing this from \eqref{3curves} 
we get that the movable part of $|V|$ is 
generated by 
$$
2f,\, 2\ol f
\hsp{\text{and}}\hsp
  f + \ol f.
$$
Recalling that $f$ and $\ol f$ are linearly 
equivalent, this means that the 2-dimensional system $|V|$ is composed with
a pencil $|f|=|\ol f|$.
Therefore the rational map associated to $|V|$ is 
factorized as $S\to \CP^1
\stackrel{\iota}{\to} \CP^2$,
where $S\to \CP^1$ is the map associated to
$|f|$, and $\iota:\CP^1\to \CP^2$ is an embedding
as a conic.
Hence we have obtained the commutative diagram
\begin{equation}\label{019}
 \CD
 S@>\phi>> \CP^4\\
 @VVV @VV{p}V\\
 \CP^1@>{\iota}>> \CP^2
 \endCD
 \end{equation}
This diagram implies that the bi-anticanonical image $\phi(S)$ is 
contained in a quadric  $p^{-1}(\iota(\CP^1))$. 
Hence in order to show that $\phi(S)$ is an intersection of two quadrics,
it is enough to show that the image under $\phi$ of a general fiber of the vertical map $S\to \CP^1$
is a conic (in a fiber plane of the projection $p$).

For this 
let $D$ be a general fiber of $S\to\CP^1$.
Then $D\in |f|$.
We can readily see that  $D$ is a smooth rational curve.
Further we can compute $(2K^{-1}-B,\, D)=2$.
Hence we obtain an 
exact sequence
$$
0 \lra 2K^{-1}-B-D \lra 
2K^{-1}-B \lra \mathscr O_D(2) \lra 0.
$$
By computing intersection numbers, we readily deduce that 
the system $|2K^{-1}-B-D|$ has $|f|$ 
as a movable part.
This implies $h^0(2K^{-1}-B-D)=2$.
On the other hand by Riemann-Roch we 
get $\chi(2K^{-1}-B-D)=2$.
Further $H^2(2K^{-1}-B-D)=0$.
These imply $H^1(2K^{-1}-B-D)=0$.
Therefore the restriction map
$H^0(2K^{-1}-B)\to H^0(\ms O_D(2))$ is 
surjective.
This means that the image $\phi(D)$ is a conic.
Hence $\phi(S)$ is an intersection of two quadrics.
The surjectivity also implies
that the restriction $\phi|_D:D\to \phi(D)$ is biholomorphic.
Hence by the diagram \eqref{019},  $\phi$ is birational over the image $\phi(S)$,
finishing a proof of the lemma.
\proofend

\vsp

Now we are able to give a classification
 in the case $k=5$:

\begin{proposition}\label{prop:k5}
Let $Z$ be a Moishezon twistor space on $4\CP^2$ satisfying $h^0(F)=2$.
Suppose that  the anticanonical cycle $C$ in 
a real irreducible member $S\in|F|$ consists
of $10$ irreducible components.
Then one of the following two situations  happens:
\begin{enumerate}
\item[(i)] $h^0(2F)= 7$, and  the anticanonical map $
\Phi:Z\to \CP^6$ is birational to the image.
Further, the image $\Phi(Z)$ is a complete intersection of three quadrics.
\item[(ii)] $h^0(2F)=5$, the image of the anticanonical map $\Phi:Z\to \CP^4$ is  a scroll of planes over a conic, and
the anticanonical map is two to one over the scroll.
Moreover the branch divisor of $\Phi$ 
 is an intersection of the scroll with 
a single quartic hypersurface.
\end{enumerate}
\end{proposition}
\noindent
\proof
First suppose that the sequence of 
the self-intersection numbers for the 
components of $C$ is given by \eqref{string5}.
Then by  Proposition \ref{prop:acs} and Lemma \ref{lemma:SIV} (i) we obtain $h^0(2F) = 5$.
Further, since the restriction of $\Phi$ to any smooth member $S\in |F|$ 
is exactly the bi-anticanonical map by the exact sequence \eqref{acs2},
from the diagram \eqref{CD1} and Lemma \ref{lemma:SIV} (iii),
we conclude that $\Phi$ is two to one over the scroll $Y$ of planes over the conic.
Further,  again by  the lemma, we obtain that 
the branch divisor of the anticanonical map $\Phi:Z\to Y\subset\CP^4$
is an intersection of the scroll by a quartic hypersurface.
Thus $Z$  satisfies all the properties in the case (ii) of the proposition.

Next suppose that 
the sequence of 
the self-intersection numbers for $C$ is given by \eqref{string6}.
Then by Proposition \ref{prop:acs}  and Lemma \ref{lemma:83} (i)  we have $h^0(2F) = 7$.
Further by the commutative diagram \eqref{CD1},
the image $\Phi(Z)$ is
contained in 
the quadric $\pi^{-1}(\Lmd)
\subset\CP^6$.
Moreover 
 the restriction of $\Phi$ to a smooth member $S$ of $|F|$ is
exactly the bi-anticanonical map of $S$
by the same reason to above.
Hence in order to prove that the image $\Phi(Z)$ is a complete intersection of three quadrics
(in $\CP^6$), it suffices to see that a general member of the pencil $|F|$ is mapped onto
a complete intersection of two quadrics in a fiber of $\pi$.
But this is exactly the assertion (iii) of Lemma \ref{lemma:83}.
Hence $\Phi(Z)$ is a complete intersection
of three quadrics in $\CP^6$.
Finally, the birationality of $\Phi$ now follows immediately from the diagram 
\eqref{CD1} and the birationality of 
the bi-anticanonical map $\phi$
proved as Lemma \ref{lemma:83} (iii).
Thus $Z$ satisfies all the properties 
in (i) of Proposition \ref{prop:k5}.
\proofend

\vsp
We note that the twistor spaces studied in \cite{HonDSn1} (specialized to the case of $4\CP^2$)
fall into the case (ii) in
Proposition \ref{prop:k5}.
On the other hand, the twistor spaces in the case (i) seem to have 
not appeared in the literature.
But in light of the presence of a $\CC^*$-action
on the divisor $S$ which is not semifree, it is quite likely that these twistor spaces admit a non-semifree $\CC^*$-action, and can be 
considered as a mild specialization of the twistor spaces
constructed in another paper \cite{Hon-I}.

\subsection{Classification in the case $k=4$}\label{ss:k4}

Let $Z$ be a Moishezon twistor
space on $4\CP^2$ which satisfies 
$h^0(F) = 2$ and $S\in |F|$ 
a real irreducible member.
Suppose that the number $m=2k$ of the 
irreducible components of the 
anticanonical cycle $C$ on $S$ is $8$.

As in the cases of $k=6$ and $k=5$, 
we take the 
birational morphism $\epsilon:S\to \FF_0$ such that 
$\epsilon(C)$ consists of four components as in \eqref{003}.
Again we factorize 
$\epsilon$ as $\epsilon =\epsilon_4\circ\epsilon_3\circ\epsilon_2\circ\epsilon_1$ as in \eqref{dcp1}. 
 Then by the assumption 
$k=4$ we can suppose that $\epsilon_2\circ\epsilon_1$  blows up the four singular points of the cycle $\epsilon(C)$ on $\FF_0$, and that $\epsilon_3$ and $\epsilon_4$ blow up smooth points 
of the anticanonical cycle.
Let $S_1$ be the toric surface obtained from $\FF_0$ by applying
the above blowups $\epsilon_2\circ\epsilon_1$,
and we name the components of the anticanonical cycle on $S_1$ as
$C_1,C_2,C_3,C_4,\ol C_1,\ol C_2,\ol C_3,\ol C_4$ in a way that adjacent curves intersect as before and that $C_1$ is a $(-2)$-curve.
By the real structure one of the two blown-up points of $\epsilon_3$ has to be
on $C_1\cup C_2\cup C_3\cup C_4$, and the same for  $\epsilon_4$.
Hence the surface $S$ is determined by choosing two points on $C_1\cup C_2\cup C_3\cup C_4$, which are not  the singularities
 of the cycle.
In the following for simplicity we do not mention the conjugate 
blown-up points.

First suppose that $\epsilon_4\circ\epsilon_3$ blows up two points on $C_1$;
here we are allowing these points to be the same point, in which case $\epsilon_4$ blows up
the intersection  of the exceptional curve with the strict transform of $C_1$.
Including this case we readily see that the resulting $S$ satisfies $h^0(K_S^{-1})=3$, contradicting $h^0(F)=2$.
Hence by an obvious holomorphic symmetry of the toric surface $S_1$, we also
obtain that  $\epsilon_4\circ\epsilon_3$ cannot
blow up two points on $C_3$.

Next suppose that $\epsilon_4\circ\epsilon_3$ blows up one point on $C_1$ and one point on $C_2$.
We note that in this case  $H^0(\Theta_S)=0$ holds, where $\Theta_S$ denotes the 
tangent sheaf of $S$.
Then using the same letters for the strict transforms into $S$, the sequence of the self-intersection numbers of the 8 curves $C_1,C_2,\cdots,\ol C_3,\ol C_4$ on $S$ is given by
\begin{align}\label{string7}
-3,-2,-2,-1,-3,-2,-2,-1.
\end{align}

\begin{lemma}\label{lemma:SIII}
The bi-anticanonical system of this surface $S$
satisfies the following.
\begin{enumerate}
\item[(i)]
$h^0(2K_S^{-1}) = 3$.
\item[(ii)] $\Bs\,|2K_S^{-1}| = \sum_{i=1}^3(C_i + \ol C_i)$,
and after removing this  curve, the bi-anticanonical system is base point free.
\item[(iii)] The morphism $\phi:S\to\CP^2$ induced by 
the bi-anticanonical system is of degree two, and 
the branch divisor is an irreducible  quartic curve.
\item[(iv)] The morphism $\phi$ contracts the connected curve
$\ol C_4\cup C_1\cup C_2$ and the conjugate curve
to mutually distinct points.
\item[(v)] The morphism $\phi$ maps $C_3$ and $\ol C_3$ to
an identical line isomorphically.
\item[(vi)] The branch quartic curve has $A_2$-singularities
at the two points which are the images of the curves in (iv).
(See Figure \ref{fig:quartics}, type III.)
\end{enumerate}
\end{lemma}

We omit a proof of this lemma as it can be
shown in a standard way.

Next suppose that $\epsilon_4\circ\epsilon_3$ blows up one point on $C_1$ and one point on $C_3$.
In this case we again have $H^0(\Theta_S)=0$.
Then this time the sequence of the self-intersection numbers of the 8 components $C_1,C_2,\cdots,\ol C_3, \ol C_4$ 
on $S$ is given by
\begin{align}\label{string8}
-3,-1,-3,-1,-3,-1,-3,-1.
\end{align}
\begin{lemma}\label{lemma:bas2}
The bi-anticanonical system of this surface $S$ has the following properties.
\begin{enumerate}
\item
[(i)] $h^0(2K_S^{-1})=5$.
\item
[(ii)] The fixed components of \,$|2K_S^{-1}|$ are $C_1+C_3+\ol C_1+\ol C_3$.
\item
[(iii)] After subtracting these components, the system is base point free.
\item
[
(iv)] The morphism $\phi:S\to\CP^4$ induced by the system is birational over the image and contracts 
the four $(-1)$-curves $C_2, C_4, \ol C_2, \ol C_4$.
Further the image $\phi(S)\subset\CP^4$ is a quartic surface.
\item
[(v)] If the birational morphism $\phi$ does not contract any other curves, then the quartic surface
$\phi(S)$ is an intersection of two 
quadrics in $\CP^4$.
\end{enumerate}
\end{lemma}

Note that in (iv) we are not claiming that the birational morphism
$\phi$ does not contract any curve other than the four
$(-1)$-curves.
Actually if there exists a $(1,1)$-curve
on $\FF_0$ passing through the four 
blown-up points of $\epsilon_4\circ\epsilon_3$,
then its strict transform into $S$ is contracted to 
an ordinary double point.
Later (in the proof of Proposition \ref{prop:k4}) we will show that 
this kind of curve does not exist
for a {\em general} smooth member $S$ of the pencil $|F|$.

\vsp
\noindent {\em Proof of Lemma \ref{lemma:bas2}.}
Again in this proof we write $K^{-1}$ for $K_S^{-1}$.
For (i)--(iii), as $(K^{-1}, C_i)=(K^{-1},\,\ol C_i)=-1$ for $i=1,3$, these curves are
base curves of $|2K^{-1}|$.
Put $B:= C_1+C_3+\ol C_1+\ol C_3$.
Then since 
$2K^{-1} - B - ( C_2+C_4+\ol C_2 + \ol C_4) = K^{-1}$,
we have the exact sequence 
\begin{align}\label{ses:4}
0 \,\lra\, K^{-1} \,{\lra}\, 2K^{-1} -B\, \lra \,\mathscr O_{C_2 \,\cup\, C_4\,\cup\,\ol C_2 \,\cup\, \ol C_4}\,\lra\, 0.
\end{align}
From this noting $H^1(K^{-1})=0$, we obtain $h^0(2K^{-1})=5$ and that $|2K^{-1}
-B|$ has no fixed component.
These mean (i) and (ii).
Suppose ${\rm Bs}\,|2K^{-1}-B|\neq\emptyset$.
Then again by the cohomology exact sequence of \eqref{ses:4}, 
it is contained in the support of $B$, and also it is disjoint from $C_2\cup C_4\cup\ol C_2 \cup \ol C_4$.
Since the restriction of the line bundle 
$2K^{-1} -B$ to the $(-3)$-curves $C_i$ and $\ol C_i$ ($i=1,3$)
is of degree 1, base points
of $|2K^{-1}-B|$ appears at one point at most on each of these curves.
By reality, the number of base points is either 2 or 4.
If it is 2, the base points are resolved by blowing-up the 2 points,
and consequently we obtain a morphism $\tilde S\to \CP^4$, where
$\tilde S$ is the 2 points blowup of $S$.
Then since $(2K^{-1}-B)^2= 4$ on $S$, the self-intersection number of the 
base point free linear system on $\tilde S$ is 2.
This means that the image of $\tilde S$ in $\CP^4$ is a quadratic surface.
But a quadratic surface in $\CP^4$ necessarily degenerates (see \cite{R}),
and hence the number of the base points of $|2K^{-1}-B|$ cannot be 2.
If the number of the base points is 4,
 $|2K^{-1}-B|$ has one base point on each of the four $(-3)$-curves,
and by blowing-up these four points we obtain another 
surface $\tilde S$ and a fixed point free linear system on $\tilde S$. 
Since the self-intersection number of this system is zero, the image of the morphism $\tilde S\to \CP^4$ 
induced by the system is a curve.
Let $E$ be any one of the four exceptional curves of the blowup $\tilde S\to S$.
Then the intersection number of $E$ with the free system on $\tilde S$ is 1.
Hence the image of $E$ under the morphism $\tilde S\to \CP^4$ must be a line.
These mean that the image curve in $\CP^4$ is a line, which cannot happen.
Thus the number of the base points of $|2K^{-1}-B|$ cannot be 4.
This implies that  $|2K^{-1}-B|$ is base point free.
Hence we obtain (iii) of the lemma.

Let $\phi:S\to\CP^4$ be the morphism induced by $|2K^{-1}-B|$.
Since $(2K^{-1}-B)^2=4$, the image $\phi(S)$ is a quadratic surface and $\phi$ is of degree 2, or
otherwise $\phi(S)$ is a quartic surface and $\phi$ is birational.
But as above it cannot be a quadratic surface.
Therefore $\phi$ is birational and $\phi(S)$ is a non-singular quartic surface.
Then as $(2K^{-1}-B,\, C_i) = (2K^{-1}-B,\,\ol C_i) = 0$ for $i=2,4$,
these four $(-1)$-curves are contracted to points by $\phi$.
Thus we get (iv) of the lemma.

Finally, for (v), suppose that $\phi$ does not contract any curve except the four $(-1)$-curves $C_2,C_4,\ol C_2$ and $\ol C_4$.
Then the quartic surface $\phi(S)$ is clearly non-singular.
It is classically known (see \cite{R}) that 
 a smooth quartic surface in $\CP^4$ is either
the image of 
a Veronese surface under a projection  $\CP^5\to \CP^4$,
or a (complete) intersection of two quadrics.
Further these two quartic surfaces can be distinguished by
the genus of a general hyperplane section:
for the Veronese surface the genus is zero,
while  for the complete intersection of the quadrics
the genus is one.
Hence to show (v) it is enough to show that the last genus is one for the present
surface $\phi(S)$.
By Bertini's theorem a general hyperplane section of $\phi(S)$
is irreducible and non-singular. 
Moreover such a hyperplane section is
clearly a biholomorphic image of
a member of $|2K^{-1}-B|$.
Hence it suffices to see that the
arithmetic genus of the last system is one.
In fact, we have
\begin{align*}
\frac12 \left(2K^{-1}-B,\,2K^{-1}-B+K\right)+1 & =
\frac12 \left(2K^{-1}-B\right)^2 + \frac12 \left(2K^{-1}-B,K\right)+1\\
&= \frac12 4 + 0 + \frac12 (-4) + 1 = 1.
\end{align*}
Hence we obtain the assertion (v) of the lemma.
\proofend

\vsp
Next for another possibility suppose that $\epsilon_4\circ\epsilon_3$ blows up two points on $C_2$.
Then the resulting surface $S$ has a non-trivial $\CC^*$-action, and the sequence of the self-intersection numbers of the 8 curves $C_1,C_2,\cdots,\ol C_3,\ol C_4$ on $S$ becomes
\begin{align}\label{string9}
-2,-3,-2,-1,-2,-3,-2,-1.
\end{align}
This surface $S$ is identical to the one
constructed in \cite[Section 2]{Hon-I}
(if specialized to the case of $4\CP^2$).
\begin{lemma}\label{lemma:base3}
For the anticanonical system of this surface $S$ we have the following:
\begin{enumerate}
\item[(i)] $h^0(2K_S^{-1}) = 3$.
\item[(ii)] $\Bs\,|2K_S^{-1}| = C_1+2C_2+C_3+\ol C_1+ 2\ol C_2+\ol C_3$.
\item[(iii)] After removing these curves, the system is base point free.
\item[(iv)] This free system is composed with
a pencil $|C_3 + 2 C_4 + \ol C_1| = 
|\ol C_3 + 2 \ol C_4 +  C_1|$
and hence the bi-anticanonical morphism
$\phi:S\to\CP^2$ 
is factorized as $S\to \CP^1\to\CP^2$,
where $S\to \CP^1$ is a morphism associated to
the last pencil, and $\CP^1\to\CP^2$ is an embedding as a conic.
\end{enumerate}
\end{lemma}
We omit a proof of this lemma.

As a final possibility, if $\epsilon_4\circ\epsilon_3$ blows up one point on $C_2$ and one point on $C_4$,
all the components of $C$  become $(-2)$-curves.
In this case the line bundle $K_S^{-1}|_C$ clearly becomes topologically trivial, and it follows that 
$h^0(lK_S^{-1})$ increases at most linearly as a function of $l$.
This means that $Z$ is non-Moishezon.
By symmetries of the toric surface $S_1$
(which was obtained from 
$\FF_0$ by blowing up the four singular points
of the anticanonical cycle; see 
an explanation at the beginning of this subsection), we do not need to consider the remaining cases, because they fall into 
one of the  above cases.
Thus we have seen that if $Z$ is Moishezon and $h^0(F)=2$, then $Z$
has a member $S\in |F|$ whose anticanonical cycle $C$
satisfies \eqref{string7}, \eqref{string8} or \eqref{string9}.
 
By using Lemmas \ref{lemma:SIII}, \ref{lemma:bas2} and 
\ref{lemma:base3}, we show the classification in the case $k=4$:
\begin{proposition}\label{prop:k4}
Let $Z$ be a Moishezon twistor space on $4\CP^2$ satisfying $h^0(F)=2$.
Suppose that the
 anticanonical cycle $C$ in a real irreducible member $S\in |F|$
consists of 8 irreducible components.
Then one of the following three situations  happens:
\begin{enumerate}
\item[(i)] $h^0(2F)=7$ and  the anticanonical map
$\Phi:Z\to\CP^6$ is birational over the image. Moreover
the image $\Phi(Z)$ is a complete intersection of three quadrics.
\item[(ii)] $h^0(2F)=5$,  the image of the anticanonical map $\Phi:Z\to \CP^4$ is  a scroll of planes over a conic, and
$\Phi$ is rationally 2 to 1 over the scroll.
Moreover  the branch divisor of $\Phi$ 
is an intersection of
 the scroll with a single quartic hypersurface.
\item[(iii)] $h^0(2F)= 5$, and the image of the anticanonical map $Z\to \CP^4$ is a complete
intersection of two quadrics, whose defining equations
are of the form
\begin{align}\label{ci}
z_0^2 = z_1z_2,\,\,z_3z_4 = q(z_0,z_1,z_2),
\end{align} 
where $z_0,z_1,\cdots,z_4$ are homogeneous coordinates on
$\CP^4$ and $q$ is a quadratic polynomial in $z_0,z_1,z_2$ with real coefficients.
Moreover, a general fiber of $\Phi$ is a
non-singular rational curve.
\end{enumerate}
\end{proposition}

\begin{remark}
{\em As will be turn out at the end of the next subsection,
the case (iii) is the unique situation where
the  image of the anticanonical map is not 3-dimensional,
provided that $h^0(F)=2$.
In other words, except this case, the anticanonical image is 3-dimensional when $h^0(F)=2$.
}
\end{remark}
\noindent
{\em Proof of Proposition \ref{prop:k4}.} 
If the 
sequence of the self-intersection numbers
of the 8 components of the cycle $C$ is 
as in \eqref{string7},
the twistor space $Z$ satisfies all the properties in (ii) of the proposition 
by  Lemma \ref{lemma:SIII} (iii), Proposition \ref{prop:acs} and the diagram \eqref{CD1}.

Next we consider the case where 
the sequence for the cycle $C$  is as in 
\eqref{string9}.
In this case,
by Proposition \ref{prop:acs} and Lemma 
\ref{lemma:base3} (i) we have 
$h^0(2F) = 5$.  
Further by Lemma \ref{lemma:base3} (iv) 
and the commutative diagram \eqref{CD1}, the anticanonical image $\Phi(Z)\subset\CP^4$ is 2-dimensional.
For the defining equations
of this surface, 
by the same reason to the proof
of Proposition \ref{prop:k6},
we have the first equation $z_0^2 = z_1z_2$ for appropriate real basis
$z_0,z_1,z_2$ of $S^2H^0(F)$.
For the second equation,
by using the $\CC^*$-action on $S$, again as in the proof of 
Proposition \ref{prop:k6},
we can show that the reducible divisor
$$
X:= S_1^+ + S_2^+ + S_3^+ + S_4^-
$$
belongs to $|2F|$.
Let $z_3\in H^0(2F)$ be an element
which satisfies $(z_3) = X$,
and put $z_4=\ol{\sigma^*z_3}$.
Then obviously $z_3z_4$ is 
a real element of $S^4H^0(F)$, so 
there exists a quadratic polynomial
$q(z_0,z_1,z_2)\in\RR[z_0,z_1,z_2]$ 
which satisfies $z_3z_4 =q(z_0,z_1,z_2)$.
This is the second equation in \eqref{ci}.
For the assertion about fibers of $\Phi$, it suffices to show that general fibers 
of the bi-anticanonical map  $\phi$ of 
$S$ are smooth rational curves.
But by Lemma \ref{lemma:base3} (iv) the movable part of $|2K_S^{-1}|$ is
generated by the pencil 
$|C_3 + 2C_4 + \ol C_1|$,
and it is easy to see that general members of
this pencil are  non-singular rational curves.
Hence general fibers of $\phi$ are non-singular rational curves.
Thus we have obtained that $Z$ satisfies all the properties in (iii) of the proposition.

We are left to show that 
if the twistor space $Z$ has $S\in |F|$ whose sequence is as in \eqref{string8},
then $Z$ satisfies the properties  (i)
of the proposition.
The assertion $h^0(2F) = 7$ follows
immediately from Proposition \ref{prop:acs} 
and Lemma \ref{lemma:bas2} (i).
In order to show that $\Phi(Z)$ is an 
intersection of three quadrics in $\CP^6$,
we need to show that 
for a  general   $S\in |F|$,
$S$ satisfies the assumption in (v) of Lemma \ref{lemma:bas2}. 
For this,
we note that
from the fact that
 $C_i$ and $\ol C_i$ are $(-3)$-curves on $S$
for $i\in \{1,3\}$, 
it readily follows that 
 these are base curves of $|2F|$, and 
the normal bundle
of these curves in $Z$ is given by 
\begin{align}\label{nbl}
N_{C_i/Z}\simeq N_{\ol C_i/Z}\simeq
\mathscr O(-2)^{\oplus 2} {\text{ for }} i\in \{1,3\}.
\end{align}
Hence 
if we blow up these curves, each exceptional divisor is 
biholomorphic to $\FF_0=\CP^1\times\CP^1$.
Then we have the following

\begin{lemma}\label{lemma:nc}
Suppose that the twistor space $Z$ has 
a real irreducible member $S\in|F|$
whose sequence is as in \eqref{string8}.
Let $\mu:\tilde Z\to Z$ be the blow-up at the base curve $B:=C_1\cup C_3\cup \ol C_1\cup\ol C_3$, $E_i$ and $\ol E_i$ $(i\in\{1,3\})$  the exceptional divisors over $C_i$ and $\ol C_i$ respectively, and put $E=E_1+E_3+\ol E_1+\ol E_3$.
Then the system $|\mu^*(2F)-E|$ is base point free and 
the resulting morphism 
$\tilde \Phi:\tilde Z\to \CP^6$
contracts $E_i\simeq\CP^1\times\CP^1$
 and $\ol E_i$ $(i=1,3)$
to $\CP^1$ along the projection which is different from the original projections
$E_i\to C_i$ and $\ol E_i\to \ol C_i$.
Moreover, other than these four exceptional divisors, the morphism $\tilde \Phi$ does not contract any divisor (in $\tilde Z$) to a lower-dimensional subvariety.  
\end{lemma}

\proof
By the exact sequence \eqref{acs2}  
we have $\Bs\, |2F| = \Bs\, |2K_S^{-1}|$, and 
the latter is the curve $B$ by Lemma \ref{lemma:bas2} (ii).
For any of the four exceptional divisors
of the blowup $\mu$, let $\mathscr O(0,1)$ be
the fiber class of the restriction of $\mu$.
Then from \eqref{nbl} we have $N_{E_i/\tilde Z}\simeq\mathscr O(-1,-2)$.
Hence we obtain
\begin{align}\label{nb2}
\left(\mu^*(2F)-E\right)|_{E_i}
\simeq
\mu^*\left(2K_S^{-1}|_{C_i}\right)\otimes N_{E_i/\tilde Z}^{-1}
\simeq
\mathscr O(0,-2) \otimes \mathscr O(1,2)
\simeq \mathscr O(1,0). 
\end{align}
Therefore by restricting the line bundle
$\mu^*F-E$ to $E_1\sqcup\ol E_1$ we have the  exact sequence
\begin{align}\label{rest60}
0 \,\lra \,\mu^*(2F) - E - E_1 - \ol E_1
\,\lra \,\mu^*(2F) - E \,\lra
\ms O_{E_1}(1,0)\oplus\ms O_{\ol E_1}(1,0)
\,\lra 0.
\end{align}
We clearly have $h^0(\tilde Z,\mu^*(2F) - E) = h^0 (Z,2F)=7$.
On the other hand, if $\tilde S$ means the 
strict transform of $S$ into $\tilde Z$
which is clearly biholomorphic to $S$, we have  
$$\mu^*(2F) - E - E_1 - \ol E_1|_{\tilde S}
\simeq 2K_S^{-1}-C-C_1-\ol C_1,$$ and we
easily have $ h^0(2K_S^{-1}-C-C_1-\ol C_1)=1$.
From this noting $\tilde S\in |\mu^*F - E|$ we readily obtain 
$h^0(\mu^*(2F) - E - E_1 - \ol E_1)\le 3$.
Then the cohomology exact sequence
of \eqref{rest60} means that 
the restriction map 
$H^0(\mu^*(2F) - E)\to H^0
(\ms O_{E_1}(1,0)\oplus\ms O_{\ol E_1}(1,0))$ is 
surjective
(and also that 
$h^0(\mu^*(2F) - E - E_1 - \ol E_1)= 3$).
This means that $\Bs\,|\mu^*(2F)-E|\cap (E_1\cup
\ol E_1)=\emptyset.$
If we replace the role of $E_1+\ol E_1$ by
$E_3 + \ol E_3$, we  obtain
$\Bs\,|\mu^*(2F)-E|\cap (E_3\cup
\ol E_3)=\emptyset.$
Therefore we conclude $\Bs\,|\mu^*(2F)-E|=\emptyset$.
Further \eqref{nb2} implies that the morphism
$\tilde{\Phi}$ contracts $E_i\simeq\CP^1\times\CP^1$ to $\CP^1$ along the projection
which is different from the original
projection $E_i\to C_i$.

It remains to see that the morphism $\tilde{\Phi}$ does not
contract any divisor other than the four components of $E$.
Let $D$ be such an irreducible divisor.
If $D$ is real,
then $D\in |\mu^*(kF)-l_1E_1-l_1\ol E_1
-l_2E_2- l_2\ol E_2|$ for some $k\ge 1$ and $l_1,l_2\ge 0$,
and $(\mu^*2F-E)^2\cdot D =0$ 
has to hold by the contractedness property.
For computing this intersection number, 
we notice that, as $\mu^*F|_{E_i} \simeq \mu^* (K_S^{-1}|_{C_i}) \simeq \mathscr O(0,-1)$ for $i\in\{1,3\}$,  we have
$(\mu^*F)^2 \cdot E_i = (\mu^*F|_{E_i})^2 = 0$ on $\tilde Z$.
Hence we have, for $i\in \{1,3\}$,
\begin{align*}
(\mu^*2F-E)^2 \cdot E_i & =  4 \mu^* F^2 \cdot E_i - 4 \mu^* F \cdot E \cdot E_i
+ E^2\cdot E_i\\
&= 4 \cdot 0 - 4 \mu^*F \cdot E_i^2 + E_i^3 \\
&= -4 \mathscr O(0,-1) \cdot \mathscr O (-1,-2) + \mathscr O(-1,-2)^2 \\
&= -4+4=0,
\end{align*}
and the same for $\ol E_i$.
From these, recalling $F^3=0$ (as we are over $4\mathbb{CP}^2$), we obtain
\begin{align*}
(\mu^*2F-E)^2\cdot D &= (\mu^*2F-E)^2 \cdot (\mu^*(kF)-l_1E_1-l_1\ol E_1
- l_2E_2-l_2\ol E_2) \\
&= (\mu^*2F-E)^2 \cdot \mu^*(kF) \\
&= (E_1^2 + \ol E_1^2 + E_2^2 + \ol E_2^2)
\cdot \mu^*(kF) 
- 4k \,\mu^*F^2\cdot (E_1+\ol E_1
+ E_2 + \ol E_2) \\
&= 4k \left(
\mu^*F|_{E_1},E_1|_{E_1}
\right)_{E_1} \\
&= 4k \mathscr O(0,-1)\cdot \mathscr O(-1,-2) = 4k.
\end{align*}
Therefore we have $(\mu^*2F-E)^2\cdot D >0$.
Hence $D$ is not contracted to a curve or a point by $\tilde\Phi$.
When $D$ is not real, by applying the above computations
to $D+\ol D$ instead of $D$,
we obtain $(\mu^*2F-E)^2\cdot (D+ \ol D) >0$.
So since $(\mu^*2F-E)^2\cdot D
= (\mu^*2F-E)^2\cdot \ol D$, 
we again conclude $(\mu^*2F-E)^2\cdot D >0$.
Therefore in the non-real case too, $D$ cannot be contracted to a point or a curve by 
the morphism $\tilde\Phi$.
Thus we have obtained Lemma \ref{lemma:nc}.
\proofend

\vsp\noindent
{\em Continuation of the proof of 
Proposition \ref{prop:k4}.} 
Arguing the case where the sequence of 
the self-intersection numbers is as in \eqref{string8}
as before,
we claim that for a general real member $S\in|F|$, the bi-anticanonical
map $\phi:S\to \CP^4$ does not contract any curve other than
the four $(-1)$-curves $C_2,C_4,\ol C_2$ and $\ol C_4$.
If this is not the case, then from the diagram \eqref{CD1} and  since the restriction of the 
anticanonical map $\Phi:Z\to \CP^6$ to $S$ is exactly
the bi-anticanonical map $\phi$, the morphism
$\tilde{\Phi}$ must contract a {\em divisor}, say $D$,  to a lower-dimensional subvariety.
But by Lemma \ref{lemma:nc} such a divisor $D$ 
must be one of $E_i$ and $\ol E_i$ $(i\in\{1,3\})$.
Thus we conclude that for a general member $S\in |F|$
the bi-anticanonical map $\phi$ is just
a (simultaneous) blowdown of the $(-1)$-curves
$C_2,C_4,\ol C_2$ and $\ol C_4$.
This means that for these surfaces the image $\phi(S)\subset\CP^4$ 
is a {\em non-singular} surface.
The latter surface is of degree 4 by Lemma \ref{lemma:bas2} (iv).
Hence by Lemma \ref{lemma:bas2} (v) we can conclude that 
$\phi(S)$ is an intersection of two quadrics.
Thus we obtain that for a generic
member $S$ of the pencil $|F|$, 
 $\phi(S)$ is a complete intersection of two quadrics.
Hence so is $\Phi(S)$.
Then by the diagram \eqref{CD1}, we obtain that 
$\Phi(Z)$ is a complete intersection of three
quadrics. Thus we have seen
that $Z$ satisfies all the  properties 
 in (i) of  Proposition \ref{prop:k4}.
\proofend

\vsp
We note that from the normal bundle of $E_i$ in $\tilde Z$ and also from
\eqref{nb2}
the image $\tilde{\Phi}(E_i)$ and $\tilde{\Phi}(\ol E_i)$ ($i\in\{1,3\}$) are lines in $\CP^6$ and the anticanonical image $\Phi(Z)=\tilde{\Phi}(Z)$ has ordinary double points along these 4 lines.
We also note that $\tilde{\Phi}$ contracts not only the four exceptional divisors of $\mu$ but also the strict transforms of the 
curves $C_i$ and $\ol C_i$ for $i\in\{2,4\}$.
With these information, it might be interesting to find explicit form of quadratic polynomials
which define the anticanonical model.

\vsp
We note that 
the twistor spaces in (iii) of Proposition \ref{prop:k4} include those studied in 
\cite{Hon-I},
while
the the twistor spaces in (i) and (ii) 
of the proposition  seem to have
not appeared in the literature.

\subsection{Classification in the cases  $k=3$ and $k=2$}\label{ss:k32}
The conclusion in these two cases is  simple:
\begin{proposition}\label{prop:k32}
Let $Z$ be a Moishezon twistor space on $4\CP^2$ satisfying $h^0(F)=2$.
Suppose that the
 anticanonical cycle $C$ in $S\in |F|$ consists
of 6 or 4 irreducible components.
Then we always have $h^0(2F)=5$.
Further the image of the anticanonical map $\Phi:Z\to \CP^4$ is  a scroll of planes over a conic, and
the anticanonical map is 2 to 1 over the scroll.
Moreover the branch divisor of $\Phi$ is an intersection 
of the scroll with a single quartic hypersurface.
\end{proposition}

\proof
Let  $\epsilon =\epsilon_4\circ\epsilon_3\circ\epsilon_2\circ\epsilon_1
:S\to\CP^1\times\CP^1$
be the decomposition as in \eqref{dcp1}.
If $k=3$, we can suppose that $\epsilon_1$ blows up a pair of singularities of the cycle $\epsilon(C)\subset\FF_0$ 
(which consists of four irreducible components
as in \eqref{003}) and that 
the remaining $\epsilon_i$-s blow up  pairs of smooth points of the resulting anti-canonical cycle.
Like before let $S_1$ be the toric surface obtained
by applying the above $\epsilon_1$, and
$C_1,\cdots,\ol C_3$  the cycle of six $(-1)$-curves on $S_1$.
If $\epsilon_4\circ\epsilon_3\circ\epsilon_2$ blows up three points on $C_1$,
we get $h^0(K_S^{-1})=3$, which contradicts our assumption.
If $\epsilon_4\circ\epsilon_3\circ\epsilon_2$ blows up two points on $C_1$,
by a holomorphic symmetry of $S_1$ we can suppose that it blows up one point on $C_2$.
Then the self-intersection numbers for the 
anticanonical cycle  $C$ on $S$ become
\begin{align}\label{string10}
-3,-2,-1,-3,-2,-1.
\end{align}
For this surface we have the following
\begin{lemma}\label{lemma:SII}
The bi-anticanonical system of this surface $S$
satisfies the following.
\begin{enumerate}
\item[(i)]
$h^0(2K_S^{-1}) = 3$.
\item[(ii)] $\Bs\,|2K_S^{-1}| = C_1 + C_2 + \ol C_1 + \ol C_2$,
and after removing this  curve, the bi-anticanonical system is base point free.
\item[(iii)] The  morphism $\phi:S\to\CP^2$ induced by 
the bi-anticanonical system is of degree two, and 
the branch divisor is an irreducible  quartic curve.
\item[(iv)] The morphism $\phi$ contracts the two connected curves
$\ol C_3\cup C_1$ and $C_3\cup \ol C_1$
to mutually distinct points.
\item[(v)] The morphism $\phi$ maps $C_2$ and $\ol C_2$ to
an identical line isomorphically.
\item[(vi)] The branch quartic curve
of $\phi$ has $A_1$-singularities
at the two points which are the images of the curves in (iv).
\end{enumerate}
\end{lemma}
Omitting a proof of this lemma
and continuing the proof of 
Proposition \ref{prop:k32},
again by the diagram \eqref{CD1} and the 
above lemma we obtain that
the image of $Z$ under the anticanonical map $\Phi$ is a 
scroll of 2-planes over a conic,
 the map is rationally 2 to 1 over the scroll,
 and that the branch divisor is a cut of the scroll by
 a quartic hypersurface.
Thus $Z$ satisfies all the properties in the proposition.

Alternatively if $\epsilon_4\circ\epsilon_3\circ\epsilon_2$ blows up one point on each of $C_1, C_2$ 
and  $C_3$, then the line bundle
$K_S^{-1}|_C$ becomes topologically trivial, which again means 
that $Z$ is non-Moishezon.
By symmetries of the above toric surface $S_1$ we do not need to consider other possibilities, and we have proved the assertions 
of Proposition \ref{prop:k32} for the case $k=3$.

If $k=2$, any of the four $\epsilon_i$-s blows up smooth points of the cycle $\epsilon(C)\subset\FF_0$.
As before write this cycle as $C_1+C_2+\ol C_1+\ol C_2$.
If $\epsilon$ blows up 4 points on $C_1$, then we get $h^0(K_S^{-1})=3$ which is a contradiction.
If $\epsilon$ blows up exactly 3 points on $C_1$, then it also blows up one point on $C_2$,
and 
the sequence of the self-intersection numbers for the 
anticanonical cycle $C$ on $S$ is given by
\begin{align}\label{string11}
-3,-1,-3,-1.
\end{align}
For this surface we have the following lemma, for which
we again omit a proof.
\begin{lemma}\label{lemma:SI}
The bi-anticanonical system of this surface $S$
satisfies the following.
\begin{enumerate}
\item[(i)]
$h^0(2K_S^{-1}) = 3$.
\item[(ii)] $\Bs\,|2K_S^{-1}| = C_1 + \ol C_1$,
and after removing this  curve, the bi-anticanonical system is base point free.
\item[(iii)] The  morphism $\phi:S\to\CP^2$ induced by 
the bi-anticanonical system is of degree two, and 
the branch divisor is an irreducible  quartic curve.
\item[(iv)] The morphism $\phi$ maps the connected curves
$C_2$ and $\ol C_2$
to mutually distinct points.
\item[(v)] The morphism $\phi$ maps $C_1$ and $\ol C_1$ to
an identical line isomorphically.
\item[(vi)] The branch quartic curve
of $\phi$ is tangent to the last
line 
at the two points which are the images of the curves in (iv).
\end{enumerate}
\end{lemma}
From this lemma, again by using the diagram \eqref{CD1}, we obtain the asserted
properties on the anticanonical system of the twistor spaces.

Finally if $\epsilon$ blows up exactly 2 points on $C_1$, then it also blows up exactly 2 points on $C_2$,
which implies that the line bundle $K_S^{-1}|_C$ is topologically trivial.
So again this cannot happen under the Moishezon assumption.
Clearly these cover all possible situations
for the case $k=2$,
and we finish a proof of Proposition 
\ref{prop:k32}.
\proofend

\vsp
We remark that when $k=2$, the twistor spaces in 
Proposition \ref{prop:k32} are identical 
to those treated in \cite[Theorem 7.11]{P92},
where it was proved that the twistor spaces
are Moishezon by showing that the anticanonical
map has 3-dimensional image.
On the other hand, when $k=3$,
the twistor spaces in 
Proposition \ref{prop:k32} seem to have
not appeared in the literature.

\section{Summary 
on the classification, and the dimension of the moduli spaces}
\label{s:summary}
\subsection{Summary and conclusions
on the classification}
\label{ss:smmr}

In Sections \ref{ss:k6}--\ref{ss:k32}  we have classified Moishezon twistor spaces on $4\CP^2$
which satisfy $h^0(F)=2$.
In order to place those satisfying $h^0(F)>2$ in the same perspective,
we first investigate the anticanonical system of those twistor spaces.

As in Proposition \ref{prop:Kr}, if a twistor space on $4\CP^2$ satisfies 
$h^0(F)>2$, then $Z$ is either a LeBrun twistor space or a Campana-Kreussler twistor space.
For LeBrun twistor spaces,  by putting $n=4$ in \cite[Proposition 2.2 and Corollary 2.3]{HonLB2}, we have
$h^0(2F) =h^0(\FF_0,\mathscr O(2,2)) = 9$ and the system $|2F|$ is generated by
$|F|$.
Since $h^0(F)=4$ and $|F|$ induces a
rational map onto a smooth quadric in $\CP^3$,
we obtain the commutative diagram of left
\begin{equation}\label{017}
 \CD
 Z_{\rm{LB}}@>\Phi_{|2F|}>> \CP^8\\
 @V \Phi_{|F|} VV @VV{\iota}V\\
 \CP^3 @>>|\mathscr O(2)|> \CP^9,
 \endCD
\hspace{15mm}
  \CD
 Z_{\rm{CK}}@>\Phi_{|2F|}>> \CP^5\\
 @V \Phi_{|F|} VV @|\\
 \CP^2 @>>|\mathscr O(2)|> \CP^5,
 \endCD
 \end{equation}
where $\iota$ is an embedding as a hyperplane.
In particular, the image $\Phi_{|2F|}(Z)$ is an embedded image  of a smooth quadric in $\CP^3$ under the Veronese image $\CP^3\subset \CP^9$, which is necessarily contained in a  hyperplane in $\CP^9$.

We next see that the situation is  similar  for the Campana-Kreussler twistor spaces.
Let $Z$ be such a twistor space on $4\CP^2$.
Then we have $h^0(F)=3$, and by  \cite[Proposition 1.3 (ii)]{CK98} 
the rational map $Z\to\CP^2$ associated to $|F|$ is surjective.
Also by  \cite[Lemma 1.8]{CK98} we have $h^0(2F)=6$.
These imply that $H^0(2F)$ is generated by $H^0(F)$.
Hence analogously to the left diagram in  \eqref{017} we get the right commutative diagram in \eqref{017}.
 Therefore the image by the anticanonical map 
of $Z$ is the 
Veronese surface $\CP^2$ in $\CP^5$.

Thus in contrast sharply with the case $h^0(F)=2$, for twistor spaces satisfying $h^0(F)>2$,
the anticanonical system $|2F|$ is always generated by $|F|$.

With these remarks, we can summarize our classification 
in the following simple form:
\begin{theorem}\label{thm:smr}
Let $Z$ be a Moishezon twistor space on $4\CP^2$.
Then the anticanonical map $\Phi$ of $Z$ satisfies 
 one of the following three properties:
\begin{enumerate}
\item
[(a)] $\Phi$ is birational over the image,
\item
[(b)] $\Phi(Z)$ is the scroll\, $Y\subset\CP^4$ over the conic (see
Definition\,\ref{def:scroll}), and $\Phi$ is  (rationally) two to one over the scroll,
\item [(c)]  $\Phi(Z)$ is a rational surface.
\end{enumerate}
\end{theorem}

\proof
The situation (a) happens in Proposition \ref{prop:k6} (i) and (ii) (the case $k=6$),
Proposition \ref{prop:k5} (i) (among the case $k=5$) and
Proposition \ref{prop:k4} (i) (among the case $k=4$),
and there exists no other Moishezon twistor space whose  $\Phi$ is birational,
because the invariant $k$  catches all Moishezon twistor spaces satisfying
$h^0(F)=2$ by Proposition \ref{prop:cycle}, and for Moishezon twistor spaces with $h^0(F)>2$,
$\Phi$ is not birational by the above considerations for the twistor spaces
of LeBrun's and Campana-Kreussler's.
The situation (b) occurs in Proposition \ref{prop:k5} (ii) 
(among the case $k=5$) and
Proposition \ref{prop:k4} (ii) (among the case $k=4$) and 
Proposition \ref{prop:k32} (the cases $k=3$ and $k=2$),
and these are all such examples.
The remaining Moishezon twistor spaces are those in Proposition \ref{prop:k4} (iii)
(among the case $k=4$) and the twistor spaces of LeBrun's and
Campana-Kreussler's.
We have already verified that all these twistor spaces satisfy
$\dim \Phi(Z)=2$.
Further, for LeBrun's and Campana-Kreussler's, $\Phi(Z)$ is
rational since it is $\FF_0$ and $\CP^2$ respectively.
For the twistor spaces in Proposition \ref{prop:k4} (iii),
the surface $\Phi(Z)$ is also rational since the defining equations
\eqref{ci} means that $\Phi(Z)$ is birational to a conic bundle
over the conic $\Lmd\simeq\CP^1$.
Thus these three kinds of twistor spaces satisfy the property (c).
This finishes a proof of the theorem.
\proofend 

\begin{definition}\label{def:type}{\em
We call Moishezon twistor space $Z$
on $4\CP^2$ to be
{\em birational type, double solid type},
or {\em conic bundle type} if $Z$ satisfies
(a), (b) or (c) in Theorem \ref{thm:smr} respectively.
}
\end{definition}

In terms of these types, we can summarize the results in Sections \ref{ss:k6}--\ref{ss:k32}
and an explanation at the beginning of this section  as follows:

\begin{theorem}\label{thm:birational}
If a Moishezon twistor space $Z$
on $4\CP^2$ is of birational type, then we have $h^0(2F)=9$ or
$h^0(2F) = 7$.
If $h^0(2F)=7$, the birational image $\Phi(Z)\subset\CP^6$
is a complete intersection of three quadrics.
If $h^0(2F) = 9$, the birational image $\Phi(Z)\subset\CP^8$
is a non-complete intersection of 10 quadrics.
\end{theorem}

\begin{theorem}\label{thm:conicbundle}
If a Moishezon twistor space $Z$
on $4\CP^2$ is of conic bundle type, then
$h^0(2F) = 9,6$ or $5$,
and general fibers of $\Phi$ are non-singular rational curves.
Further,
the image rational surface $\Phi(Z)$ can be explicitly described as follows:
\begin{enumerate}
\item[(i)] If $h^0(2F) = 9$, 
it is a smooth hyperplane section of the Veronese embedding of $\CP^3$ 
into $\CP^9$ by $|\mathscr O(2)|$.
\item[(ii)] If $h^0(2F) = 6$, it is a Veronese surface in $\CP^5$.
\item[(iii)] If $h^0(2F) = 5$, it is an intersection of two quadrics whose defining equations
can be taken as $z_0^2 = z_1z_2$ and $z_3z_4= q(z_0,z_1,z_2)$ in homogeneous coordinates,
where $q$ is a quadratic polynomial in $z_0,z_1,z_2$ with real coefficients.
\end{enumerate}
\end{theorem}

\begin{theorem}\label{thm:DS}
If a Moishezon twistor space $Z$
on $4\CP^2$ is of double solid type,
the branch divisor of $\Phi$ is an intersection of the
scroll with a single quartic hypersurface.
\end{theorem}

In view of the last theorem, 
since the double covering is mostly determined
from the branch divisor, it is a basic question
 to determine the defining equation of 
the quartic hypersurface in $\CP^4$
which cuts out the branch divisor of the
anticanonical map.
This is exactly the question pursued by 
 Poon and Kreussler-Kurke in \cite{P92, KK92} for the twistor spaces of $3\CP^2$
which have a structure of double covering
of $\CP^3$.
We shall give an answer to this question in Section \ref{s:de},
after eliminating the base locus of the anticanonical system
in Section \ref{s:acs}.
To this end, 
for twistor spaces of double solid type,
 we introduce the following definition that will be
used in the rest of this paper:
\begin{definition}\label{def:DS}
Let $Z$ be a twistor space on $4\CP^2$ which is of double 
solid type. Then according to the number $4,6,8$ or $10$ 
of the irreducible components of the cycle $C$
(namely according to $k=2,3,4$ or $5$ respectively),
we call $Z$ to be of {\em type I, II, III or IV} respectively.
\end{definition}

As we have already mentioned, it is quite certain
that most twistor spaces of birational type and conic bundle type are identical to 
known classes of twistor spaces as follows:
\begin{conjecture}\label{conj}
Let $Z$ be a Moishezon twistor space on $4\CP^2$.
Then the following would be true: 
\begin{enumerate}
\item If  $Z$ is birational type with $k=6$, then
$Z$ is a twistor space of a Joyce metric
\cite{J95}.
In particular, they should admit an effective
$\CC^*\times\CC^*$-action.
\item If  $Z$ is  conic bundle type with $k=4$,
then $Z$ is isomorphic to one of the  twistor spaces
constructed in \cite{Hon-I}.
In particular, they should admit a
(non-semifree) $\CC^*$-action.
\item If  $Z$ is of conic bundle type
with $k=5$, it is a degenerate form of 
the twistor spaces in \cite{Hon-I}.
In particular they also admit 
a (non-semifree) $\CC^*$-action, and
their minitwistor spaces
with respect to the $\CC^*$-action are isomorphic to
the ones in \cite{Hon-I}.
\end{enumerate}
\end{conjecture}
The point is to show that the $\CC^*\times \CC^*$- or $\CC^*$-action on a real irreducible member
$S\in |F|$ extends to that on the twistor space.
In view of the Pontecorvo's theorem
\cite{Pont92},
it seems quite natural to expect that 
there always exists a natural linear injection 
$H^0(\Theta_S)\to H^0(\Theta_Z)$,
which implies the above extension property. 

\subsection{Dimension of the 
moduli spaces}
\label{ss:moduli}
In this subsection 
based on the results in Section \ref{s:clsf}
we compute dimension of the moduli spaces of 
Moishezon twistor spaces on $4\CP^2$.
We begin with 
 the following proposition 
which implies that the inverse of 
the above conjectural injection   
exists for Moishezon twistor spaces which satisfy
$h^0(F)=2$:  

\begin{proposition}\label{prop:aut1}
Let $Z$ be a Moishezon twistor space on $4\CP^2$ which satisfies
$h^0(F)=2$. 
Then for any real irreducible member $S$ of the pencil $|F|$,
there is a natural linear injection $H^0(\Theta_Z)\to H^0(\Theta_S)$.
\end{proposition}

\proof
By Proposition \ref{prop:cycle}, under the assumption
$h^0(F) = 2$, the base locus $\Bs\,|F|$ is a cycle of smooth
rational curves, and it consists of $2k$ irreducible components,
where $2\le k\le 6$.
Further the pencil $|F|$ has precisely $k$ reducible members.
Then the existence of the natural linear map 
$H^0(\Theta_Z)\to H^0(\Theta_S)$ is obvious in the case $k\ge 3$
since $\Aut_0Z$,
the identity component of the holomorphic 
automorphism group of $Z$,
 naturally acts linearly on $H^0(F)$
 and preserves each
reducible members.
When $k=2$, let $S_1^++S_1^-$ and $S_2^++ S_2^-$ be
the reducible members of $|F|$ as before.
Then these are also invariant under the action of $\Aut_0Z$, and moreover,
the intersection twistor line $L_i$
is invariant under this action.
On the other hand, it can be readily seen that 
in this case any  $S_i^+$ (and $S_i^-$) is obtained from
$\CP^2$ by blowing-up 4 points, and that 
precisely 3 of the 4 points are on a line.
If $\Aut_0(S_i^+,L_i)$ means
the identity component of the group
of holomorphic automorphisms of $S_i^+$ 
which preserve $L_i$,
this means $\Aut_0(S_i^+,L_i)=\{e\}$. 
Since the homomorphism $\Aut_0Z\to \Aut_0(S_i^{+},L_i)$ is clearly injective
(by going down to $4\CP^2$), this means  $\Aut_0Z=\{e\}$.
Therefore we obtain $H^0(\Theta_Z)=0$.
(Namely if $k=2$, the asserted injection is
a trivial map from $0$ to $0$.)
\proofend

\vsp
 
 Proposition \ref{prop:aut1} implies that $H^0(\Theta_S)=0$ means
$H^0(\Theta_Z) = 0$.
As we have mentioned in Sections \ref{ss:k6}--\ref{ss:k32}
where we presented an explicit construction of $S$
as a blowup of $\FF_0$,
this is the case  exactly when 
\begin{itemize}
\item $k\in\{2,3\}$, or 
\item $k=4$ and $Z$ is either a birational type or a double solid type.
\end{itemize}
(For all other remaining cases
we have $H^0(\Theta_S)\neq 0$.)

We compute the dimension of the moduli
space for those twistor spaces.
Let $Z$ be any one of these Moishezon twistor spaces
and $S\in |F|$ a real irreducible member.
Let $\Theta_{Z,S}$ be the subsheaf
of $\Theta_Z$
formed by germs of vector fields
which are tangent to $S$.
We also write $\Theta_Z(-S)$ for
$\Theta_Z\otimes\ms O_Z(-S)$. 
Then by \cite[Theorem 5.1]{Fu04}, 
we have
\begin{align}\label{van2}
H^2(\Theta_Z) = H^2(\Theta_{Z,S}) = 
H^2 (\Theta_Z(-S))=0.
\end{align}
Further by computations in 
\cite[p.146]{Hi81}, we have 
$\chi(\Theta_Z) = 15- 7\cdot 4=-13,$
where $\chi$ denotes the Euler characteristic.
Then since we always have $H^3(\Theta_Z)=0$ and are assuming 
$H^0(\Theta_Z)=0$,
we obtain
\begin{align}\label{tRR0}
h^1(\Theta_Z) = 13. 
\end{align}
We also have, for the rational surface $S$, 
by the Riemann-Roch formula
\begin{align}\label{acRR}
 h^1(K_S^{-1}) =  h^0(K_S^{-1})-1
\end{align} 
and
\begin{align}\label{tRR}
 h^1(\Theta_S) =  10,
\quad H^2(\Theta_S) = 0.
\end{align} 
Then by various standard exact sequences of sheaves including the above ones,  noting $N_{S/Z}\simeq K_S^{-1}$
and using \eqref{van2},
we obtain the following commutative diagram of cohomology groups
on $Z$ and $S$:
$$
\begin{CD}
@. @. 0\\
@. @. @VVV\\
@. 0 @>>> H^0(K_S^{-1}) @= H^0(K_S^{-1}) @>>> 0\\
@. @VVV @VV{\delta}V @VVV @.\\
0@>>> H^1(\Theta_Z(-S)) @>>> H^1(\Theta_{Z,S})  @>{\alpha}>> H^1(\Theta_S)
@>>>0\\
@. @| @VV{\beta}V @VVV @.\\
0 @>>> H^1(\Theta_Z(-S)) @>>> H^1(\Theta_Z) 
@>>> H^1(\Theta_Z|_S) @>>>  0\\
@. @VVV @VVV @VVV @.\\
 @. 0 @>>> H^1(K_S^{-1}) @=
H^1(K_S^{-1})
@>>>0\\
@. @. @VVV @VVV\\
@. @. 0 @.0
\end{CD}
$$
From the middle column of this diagram,
by using \eqref{tRR0} and \eqref{acRR}, we obtain 
\begin{align}
h^1(\Theta_{Z,S}) &= h^1(\Theta_{Z})
+ h^0(K_S^{-1}) - h^1 (K_S^{-1})
\notag\\
&= 13 + 1 = 14.\label{modui1}
\end{align}
Hence from the second row 
(including the map $\aaa$) and \eqref{tRR}, we obtain
\begin{align}
h^1(\Theta_Z(-S)) &= h^1 (\Theta_{Z,S})
- h^1(\Theta_S)\notag\\
&= 14 - 10 = 4.\label{ker01}
\end{align}

In order to compute the dimension of the moduli spaces,
we  recall from Section \ref{s:clsf} that
the geometric structure of the twistor spaces is 
determined from the sequence of 
the self-intersection numbers of the cycle $C$ in $S\in |F|$,
and  the complex structure
of $S$ is 
determined from the  blowups 
$\epsilon:S\to \mathbb F_0$.
For the actual computation, we first
obtain for each fixed sequence
the dimension of the $\RR$-linear subspace $V$ 
of $H^1(\Theta_S)$ which is 
generated  by deformations of $S$
through the Kodaira-Spencer map,
where we consider all deformations of $S$
which keep the sequence invariant.
Then from obvious meanings of the maps $\aaa$ and $\bbb$
 (see the above diagram),  
the subspace $\beta(\aaa^{-1}(V)^{\sigma})$
in $H^1(\Theta_Z)$,
where $\aaa^{-1}(V)^{\sigma}$ denotes the 
real part of $\aaa^{-1}(V)$ in $H^1(\Theta_{Z,S})$,
 can be 
considered as the tangent space of the moduli space
of the (Moishezon) twistor spaces which have
the sequence as the self-intersection
numbers of the cycle $C$ in $S$.

For example, 
we compute the dimension of the 
$\RR$-subspace $V\subset H^1(\Theta_S)$
when the sequence for the cycle $C$ is as in \eqref{string11};
namely when $Z$ is a double solid
twistor space of type I (see
Definition \ref{def:DS}).
In this case all freedom for constructing
$S$ from the blowups $\epsilon:S\to \FF_0$ is
to move three points on $C_1$ and one point
on $C_2$ on the initial surface $\mathbb F_0$, since the remaining four points
are determined from these four points
by the real structure
(see Section \ref{ss:k32}).
But not all these moving contribute
for deforming the complex structure of $S$;
namely,  the $\U(1)\times\U(1)$-action 
on $\FF_0$ which preserves the cycle $C_1+C_2+
\ol C_1 + \ol C_2$  
clearly erases the effect of the moving.
Consequently, the $\RR$-subspace $V\subset H^1(\Theta_S)$ generated by
all real deformations of $S$ 
which keep the sequence \eqref{string11} fixed is ($2\cdot 4-2=)\,6$-dimensional.
Therefore by using \eqref{ker01}, we have 
$\dim_{\RR}\aaa^{-1}(V)^{\sigma} = 6+4=10$.
Next we see that this $\RR$-subspace
$\aaa^{-1}(V)^{\sigma}$  contains  
the real part of the image of
the vertical map $\delta$ in the commutative diagram.
The last image clearly corresponds to
deformations of the pair $(Z,S)$ which are
obtained by just moving $S$ in the real part of the pencil $|F|$,
while the complex structure of the twistor space $Z$ is fixed.
On the other hand, it is immediate to see that 
any real irreducible member $S_{\lmd}\in |F|$ also contains the 
curves $C_1, C_2,\ol C_1$ and $\ol C_2$ as
$(-3), (-1), (-3)$ and $(-1)$-curves. 
This means that the $\RR$-subspace in 
$H^1(\Theta_S)$ generated by 
the $\RR$-pencil $\{S_{\lmd}\set
\lmd\in \Lmd$ and $\lmd$ is real$\}$
is contained in the subspace $V$. 
This means
$\aaa\circ\delta(H^0(K_S^{-1})^{\sigma})\subset V$ and hence
 the desired inclusion 
$\delta(H^0(K_S^{-1})^{\sigma})\subset \aaa^{-1}(V)^{\sigma}$.
Therefore we obtain that the $\RR$-subspace
$\bbb(\aaa^{-1}(V)^{\sigma})\subset H^1(\Theta_Z)$ is isomorphic to
the quotient space
$\aaa^{-1}(V)^{\sigma}/\delta(H^0(K_S^{-1})
^{\sigma})$.
Since $h^0(K_S^{-1}) = 1$, we obtain
that the moduli space of the twistor spaces
of double solid type I is $(10-1=)\,9$-dimensional.

In a similar way, concerning other twistor spaces of double solid type,
for the dimension of the $\RR$-subspace $V\subset H^1(\Theta_S)$,
we have $\dim V= 4$ for the case of type II,
 $\dim V= 2$ for the case of type III.
(We also have $\dim V= 1$ for the case of type IV.
But we are excluding this case because  $H^0(\Theta_S)\neq 0$.)
From these, we obtain that the moduli space is
$7$-dimensional for the case of type II,
and $5$-dimensional for the case of type III.
For the twistor spaces of birational type which satisfy
$k=4$, the $\RR$-subspace $V\subset H^1(\Theta_S)$ is
readily seen to be 2-dimensional.
Hence again by using \eqref{ker01} the moduli space becomes $(2+4-1=)\,5$-dimensional.
Thus we have computed the 
dimension of the moduli space
for those satisfying $H^0(\Theta_S)=0$
and $h^0(F) = 2$.

A similar argument works for
computing the dimension of the moduli 
space of Campana-Kreussler twistor spaces
on $4\CP^2$ which satisfy $H^0(\Theta_Z)=0$.
First we readily see that 
for the $\RR$-subspace $V\subset
H^1(\Theta_S)$ we have
$\dim V=7$.
Further we have $h^0(K_S^{-1}) = 2$
and the image $\delta(H^0(K_S^{-1}))$ is 
again contained in $\aaa^{-1}(V)$
by the same reason to the above argument.
Hence the dimension of the moduli space is computed as $7+4 - 2 = 9$.
We note that by the same argument the 
moduli space of Campana-Kreussler twistor spaces on $n\CP^2$ 
is $(3n-3)$-dimensional.
We also mention that this can obtained in a simpler way by 
computing $h^1(\Theta_{Z,C_0})$, where $C_0$ is the 
base curve of the net $|F|$, instead of passing to the surface $S$.

Up to Conjecture \ref{conj}, the dimension of the moduli spaces 
of  
Moishezon  twistor spaces on $4\CP^2$  satisfying
$h^0(F)=2$  can be  listed up as in 
the following table:
 \begin{center}
\begin{tabular}{c||c|c|c}
& birational type & double solid type & conic bundle type\\
\hline
$k=6$ & 3-dim. ($\CC^*\times\CC^*$) & - & -\\
\hline
$k=5$ & 4-dim. ($\CC^*$) &  4-dim. ($\CC^*$)  & -
  \\ \hline
$k=4$ & 5-dim. & 5-dim. & 6-dim. ($\CC^*$)
\\ \hline
$k=3$ & - & 7-dim. & -
\\ \hline
$k=2$& - & 9-dim. & -\\
\end{tabular} 
\end{center}
Here, hyphen means non-existence, and 
the group in parentheses is the identity component of the holomorphic
automorphism group of the twistor space.

\section{Anticanonical system of the twistor spaces
of double solid type}\label{s:acs}

Our final goal in the rest of this paper is to obtain 
defining equation of the quartic hypersurface in $\CP^4$
which cuts out the branch divisor of the anticanonical map
of the twistor spaces of double solid type.
However there seems to be no shortcut
for obtaining this equation.
In this section for that goal
we eliminate the base locus of the 
anticanonical system for each type 
I, II, III and IV.

Let $Z$ be a Moishezon twistor space on $4\CP^2$ which is
of double solid type.
Then according to the (sub)type introduced in 
Definition \ref{def:DS}, the anticanonical cycle
$C$ of a real irreducible member $S\in |F|$
has the following sequence as the self-intersection 
numbers in $S$ of the irreducible components
(see \eqref{string11}, \eqref{string10}, \eqref{string7} and
\eqref{string5}):
\begin{align}\label{I}
-3,-1,-3,-1 \hsp&{\text{for type I}},\\
-3,-2,-1,-3,-2,-1 \hsp&{\text{for type II}},\label{II}\\
-3,-2,-2,-1,-3,-2,-2,-1 \hsp&{\text{for type III}},\label{III}\\
-3,-2,-2,-2,-1,-3,-2,-2,-2,-1 \hsp&{\text{for type IV}}.\label{IV}
\end{align}
As before these indicate that, for example in the case of type II,  the self-intersection numbers in $S$
of the components $C_1,C_2,C_3,\ol C_1,\ol C_2,\ol C_3$  are
$-3,-2,-1,-3,-2,-1$ respectively.
The structure of the bi-anticanonical systems of $S$ is 
presented  in Lemmas \ref{lemma:SI}, \ref{lemma:SII}, \ref{lemma:SIII} and 
\ref{lemma:SIV} respectively.

As in Section \ref{ss:bc}, there are precisely
$k$ reducible members
 of the pencil $|F|$,
where as before $k$ is the half of the number
$m$ of the components of the cycle $C$.
We still denote $\{S_i^+ + S_i^-\set 1\le i\le k\}$
for these members,
where we are keeping the rule for the indices and
the superscripts.
Besides the notations introduced in 
Definition \ref{def:scroll}, we use the 
following notations throughout 
Sections \ref{s:acs} and \ref{s:de}.
(Recall that the conic $\Lmd$ is naturally identified 
with the parameter space of the pencil $|F|$.)

\begin{definition}
\label{def:term2}
{\em
Let $Z$ be a twistor space on $4\CP^2$ of double solid type.
\begin{enumerate}
\item[(i)] 
We denote by $B$ for the branch divisor 
of the anticanonical map $\Phi:Z \to Y\subset\CP^4$.
\item[(ii)]
We denote by $\lmd_i\in\Lmd$ for the point 
corresponding to the reducible member
$S_i^+ + S_i^-$.
\item[(iii)]
We denote by $P_i$ for the fiber plane
$\pi^{-1}(\lmd_i)$ in the scroll $Y$ over $\Lmd$.
\end{enumerate}
}
\end{definition}

Since $\lmd_i$ is a real point, $P_i$ is a real plane.
Concerning the image of the reducible members
under the anticanonical map $\Phi$, by the commutative diagram \eqref{CD1}, 
we have  
\begin{align}\label{incl}
\Phi ( S_i^+\cup S_i^- ) \subset P_i.
\end{align}
Since $\Phi$ is surjective to $Y$,
one might first think that the coincidence
$
\Phi ( S_i^+\cup S_i^- ) = P_i
$
necessarily holds. But this is not always true
since the following situation can happen:
an exceptional divisor of an elimination
of the base locus of the anticanonical system
is mapped surjectively to $P_i$ and instead
$S_i^+$ and $S_i^-$ are mapped to a proper subvariety
in $P_i$.
We will see that this is actually the case
for the case of type II, III and IV for one 
particular reducible member.

Before proceeding to the actual analysis of the anticanonical system,
we mention a property which helps to understand the 
structure of the branch divisor $B$ of the anticanonical map, at least geometrically:

\begin{proposition}\label{prop:wsing}
For a general point $\lmd\in\Lmd$,
the intersection $B\cap \pi^{-1}(\lmd)$
is a quartic curve in the plane $\pi^{-1}(\lmd)$ which touches 
the ridge $l$ at two points,
or which has singularities
at two points on $l$.
Further, these two points are independent of the choice of $\lmd\in\Lmd$.
\end{proposition}

\proof
From the commutative diagram \eqref{CD1},
for a general $\lmd\in\Lmd$,
the corresponding member $S\in |F|$ is mapped 
to the plane $\pi^{-1}(\lmd)$, and
 from the exact sequence \eqref{acs2} the map is
identified with the bi-anticanonical map
of $S$.
Therefore
each plane $\pi^{-1}(\lmd)$ possesses a line which is special 
in relation with the branch quartic curve of the 
bi-anticanonical map,
in the sense that branch curve touches the line
at two points for the case of type I,
and it has singularities at two points on the line
for the case of type II, III and IV
((vi) of Lemmas \ref{lemma:SI}, \ref{lemma:SII}, \ref{lemma:SIII} and \ref{lemma:SIV}.
See also Figure \ref{fig:quartics}).
From these lemmas, in terms of the bi-anticanonical map
$\phi:S\to\CP^2$, 
these two points are concretely written as
$\phi(C_i)$ and $\phi(\ol C_i)$,
where $C_i$ and $\ol C_i$ are
the components of the cycle $C$ whose
self-intersection numbers are $(-1)$.
Since the anticanonical map $\Phi$ is
well-defined on a general point
of these two components,
we can replace the above  $\phi(C_i)$ and $\phi(\ol C_i)$ by 
$\Phi(C_i)$ and $\Phi(\ol C_i)$ respectively.
This means that the special line is  independent of the 
choice of $\lmd$.
Because a common line for the planes $\pi^{-1}(\lmd)$ is nothing but the ridge $l$,
this means that the special line on the plane $\pi^{-1}(\lmd)$ is
exactly $l$.
Thus we have obtained the proposition.
\proofend

\vsp
If we denote by 
$q$ and $\ol q$ for the two common points on the ridge $l$, the proposition 
implies  
\begin{align}\label{q}
B\cap l = \{q, \ol q\}.
\end{align}
Also from the proposition, the branch divisor $B$
 is obtained from the branch curve of 
the bi-anticanonical map of $S$ (which is illustrated in 
Figure \ref{fig:quartics}) by `rotating around' 
the ridge $l$.
Thus it is quite obvious that the branch divisor of the 
anticanonical map has the  worst singularities at the two points $q$ and $\ol q$.

We begin an analysis of the anticanonical map.
From the exact sequence \eqref{acs2}, 
we have $\Bs\,|F|= \Bs\,|K_S^{-1}|$.
The latter base locus is explicitly given in the above lemmas in Section \ref{s:clsf}.
In particular, the base locus is a string of smooth rational curves
which is a part of the cycle $C$ .
For the case of type I for example, the base locus 
can be eliminated by blowing-up the base curve $C_1\sqcup\ol C_1$ just once.
However, for other types II, III and IV the base locus cannot be eliminated
by blowing-up just once, and 
it seems difficult to keep track of the 
base locus if we initially blowup the base curves themselves.

So instead let $\hat Z\to Z$ be the blowup at the whole of the cycle $C$, and $E_i$ and $\ol E_i$ $(1 \le i\le k)$
the exceptional divisors over the curves $C_i$ and $\ol C_i$ respectively.
Then the indeterminacy locus of 
the rational map $Z\to \CP^1$
induced by the pencil $|F|$
is eliminated and 
the composition $\hat Z\to Z\to \CP^1$
is a morphism.
(See Figures \ref{fig:blowups_I}--\ref{fig:blowups_IV} for configuration of
these divisors on $\hat Z$.)

From the fact that $C$ constitutes a cycle, it follows that 
any $E_i$ and $\ol E_i$ are biholomorphic to $\FF_0$.
If we denote the strict transforms of the components $S_i^+$ and $S_i^-$
by the same letters, then the four divisors
$S_i^+, S_i^-, E_i$ and $E_{i+1}$ in $\hat Z$,
where $E_{k+1}$ means $\ol E_1$, share a unique point, 
for which we denote by $p_i$.
(In Figures \ref{fig:blowups_I}--\ref{fig:blowups_IV}
these are indicated by dotted points.)
Then the new variety $\hat Z$ has ordinary double points at $p_i$ and $\ol p_i$, $1\le i\le k$.
For each of them we take a small resolution such that
\begin{itemize}
\item
for $1\le i\le k-1$ the pair of the divisors $E_i$ and $S_i^+$ are blown-up at the point $p_i$,
\item
for $i=k$, the alternative pair of divisors $\ol E_1$ and $S_k^-$ are blown-up at the point $p_k$.
\end{itemize}
(See Figures \ref{fig:blowups_I}--\ref{fig:blowups_IV}.)
For the points $\ol p_i$, $1\le i\le k$, we take a small resolution which is uniquely 
determined from that for $p_i$ by the real structure,
and let $\Delta_i$, $1\le i\le k$, be the exceptional curve over the point $p_i$.
We denote by $Z_1\to \hat Z$ for the simultaneous small resolutions of all these $2k$ ordinary double points,
and write $\mu_1: Z_1\to Z$ to mean the composition $Z_1\to \hat Z\to Z$.
From the choice, the divisor $E_1$ is blown-up at the two points $p_1$ and $\ol p_k$,
and the divisor $E_i$, $1<i<k$, is blown-up at one point $p_i$, while 
$E_k$ is not affected.

Next let $E$ be the total sum of
all the exceptional divisors
of the birational morphism $\mu_1$,
and 
let $f_1:Z_1\to \CP^1$ be the morphism associated to the free system
$|\mu_1^*F-E|$.
Then we have a basic relation
\begin{align}\label{rel1}
\mu_1^*F \simeq f_1^*\mathscr O(1) + E.
\end{align}
This  is very useful for computing  intersection numbers on the manifold $Z_1$.

Also in the following we often consider  line bundles over the exceptional divisors $E_i 
\subset Z_1$, $1\le i\le k$.
To express them, we write $\mathscr O_{E_i}(a,b)$
or just $\mathscr O(a,b)$ to mean the pull-back of 
the line bundle $\mathscr O(a,b)$ over $\FF_0\simeq E_i\subset \hat Z$ 
by (the restriction of) the small resolution $Z_1\to \hat Z$,
where we promise that the line bundle $\mathscr O(0,1)$ is the fiber class
of the projection $(\hat Z \supset) E_i\to C_i$.
This means that for any member $S\in |F|$, the intersection of
its strict transform into $Z_1$ with the divisor $E_i\subset Z_1$ is a member of 
$|\mathscr O(1,0)|$ in the 
above notation for  line bundles 
over $E_i$.

\subsection{Elimination of the base locus -- the case of type I}\label{ss:elimI}
Because the divisor $E_1+ \ol E_1$ is 
a fixed component of the pullback system $|\mu_1^*(2F)|$, 
we put $$\mathscr L_1:= \mu_1^*(2F) - E_1 - \ol E_1.$$
From the relation \eqref{rel1} this can be rewritten as
\begin{align}\label{L1I}
\mathscr L_1 = f_1^*\mathscr O(2)+ E_1 + 2E_2 + \ol E_1 + 2\ol E_2.
\end{align}
Then the
base locus of the anticanonical system
is already removed at this stage: 
\begin{proposition}\label{prop:baseI}
The linear system $|\mathscr L_1|$ on $Z_1$ is base point free.
\end{proposition}

\begin{figure}
\includegraphics{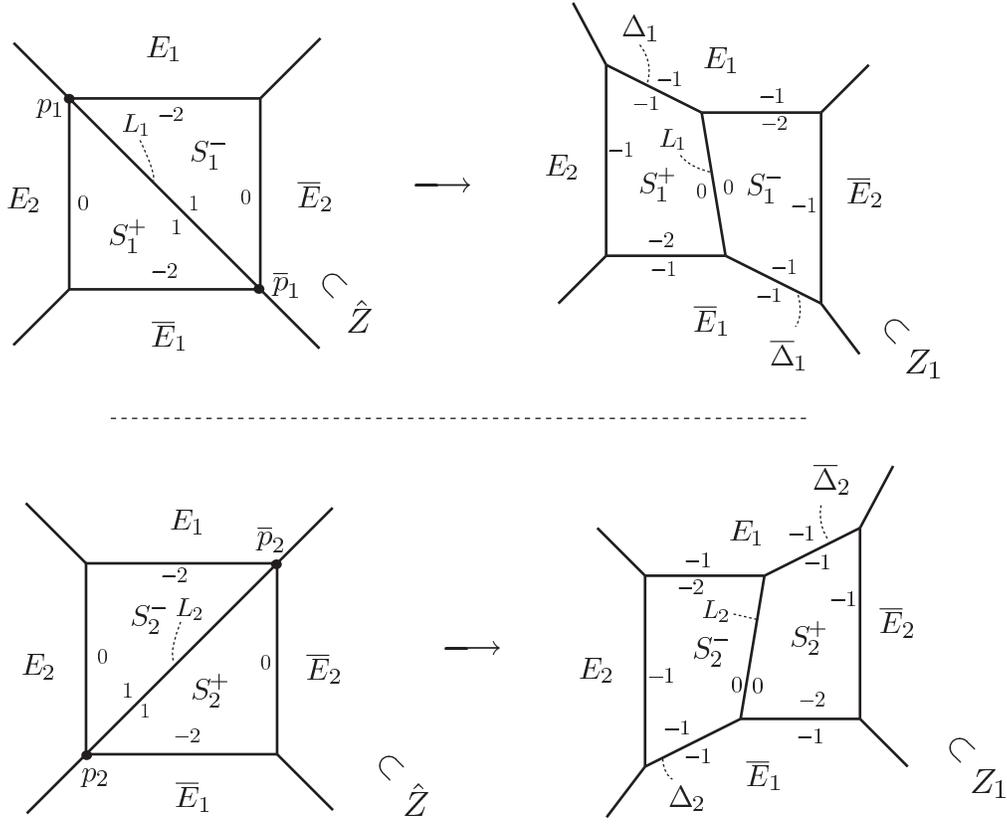}
\caption{The blowups for the case of type I.}
\label{fig:blowups_I}
\end{figure}

\proof
Since $\Bs\,|2F| = C_1\sqcup \ol C_1$,
we clearly have $\Bs\,|\mathscr L_1| \subset E_1\sqcup\ol E_1$.
We calculate the restrictions of $\ms L_1$ to the divisors $E_1$ and $E_2$
in concrete forms, by computing the intersection numbers
of $\ms L_1$ with curves on these divisors which generate $H^2(E_i,\ZZ)$.
First we compute for the exceptional curve $\Delta_1\subset E_1$.
For this, by identifying curves on the original divisor $S_1^+$ in $Z$ 
and its strict transform into $Z_1$, we temporary write 
$C_2 = S_1^+\cap E_2$ and $\ol C_1 = S_1^+\cap \ol E_1$ in this proof.
Then noting that $\Delta_1 = S_1^+\cap E_1$ and that
$\Delta_1$ does not intersect a general fiber of $f_1$ (see Figure \ref{fig:blowups_I}) so that
 $(f_1^*\mathscr O(2), \Delta_1)_{Z_1}=0$,
by using \eqref{L1I}, we have
 \begin{align}
\left(\mathscr L_1 , \Delta_1 \right)_{Z_1}  
&= \left(E_1 + 2 E_2  + \ol E_1 + 2 \ol E_2,\, \Delta_1\right)_{Z_1} 
\notag\\
&= \left(\Delta_1 + 2 C_2 + \ol C_1,\, \Delta_1\right)_{S_1^+} 
\notag\\
&= -1 + 2 + 0 = 1.\label{int0001}
\end{align}
For the exceptional curve $\ol\Delta_2$, since $\ol E_1\cap\ol\Delta_2
= E_2\cap\ol\Delta_2=\emptyset$, we have
\begin{align}
\left(\mathscr L_1 ,\, \ol{\Delta}_2\right)_{Z_1}  
&= \left(f_1^*\mathscr O(2) + E_1 + 2 E_2  + \ol E_1 + 2 \ol E_2,\, \ol{\Delta}_2\right)_{Z_1} 
\notag\\
&= \left( E_1, \, \ol{\Delta}_2 \right)_{ Z_1 } 
+ 2 \left(\ol E_2, \, \ol{\Delta}_2\right) _{Z_1} \notag\\
&=
(E_1, \,\ol{\Delta}_2)_{S_2^+}+2
= (-1) +2 = 1.\label{int00025}
\end{align}
Next for the curve $E_1\cap E_2$, 
noting that this is a section of 
$f_1:Z_1\to\CP^1$ and the self-intersection
numbers in $E_1$ and $E_2$ are $-1$ and $0$ respectively, 
we compute
\begin{align}
\left(\mathscr L_1 ,\, E_1\cap E_2\right)_{Z_1}  
&= \left(f_1^*\mathscr O(2) + E_1 + 2 E_2  + \ol E_1 + 2 \ol E_2,\, E_1\cap E_2\right)_{Z_1} 
\notag\\
&= \left(f_1^*\mathscr O(2), \, E_1\cap E_2 \right)_{ Z_1 } 
+ \left(E_1, E_1\right)_{E_2} + \left(2E_2,E_2\right)_{E_1}
\notag\\
&= 2 + 0 + 2 (-1) = 0.\label{int0002}
\end{align}
Similarly for the section $E_1\cap \ol E_2$ of $f_1$,
because it is a $(-1)$- and $0$-curve
on $E_1$ and $\ol E_2$ respectively,
we have
\begin{align}
(\mathscr L_1 ,\, E_1\cap \ol E_2)_{Z_1}  
&= \left(f_1^*\mathscr O(2), \, E_1\cap \ol E_2 \right)_{ Z_1 } 
+ \left(E_1, E_1\right)_{\ol E_2} + \left(2\ol E_2,\ol E_2\right)_{E_1}
\notag\\
&= 2 + 0 + 2 (-1) = 0.\label{int00026}
\end{align}
From \eqref{int0001}--\eqref{int00026}
we readily obtain 
\begin{align}\label{rest001}
\mathscr L_1|_{E_1} \simeq \mathscr O(1,1)-\Delta_1 - \ol{\Delta}_2.
\end{align}
In a similar way it is possible to compute the intersection numbers
of $\ms L_1$ with curves on $E_2$ 
which generate $H^2(E_2,\ZZ)$,
which turn out to be zero.
Hence we obtain
\begin{align}\label{rest001.5}
\mathscr L_1|_{E_2} \simeq  \mathscr O.
\end{align}

If $E_1$ were a fixed component of the system $|\ms L_1|$,
as $E_1\cap E_2\neq\emptyset$,
 \eqref{rest001.5} would mean that $E_2$ is also a fixed component.
This contradicts $\Bs\,|\mathscr L_1| \subset E_1\sqcup\ol E_1$,
so $E_1$ is not a fixed component of $|\ms L_1|$.
Then since the line bundle $\mathscr O(1,1)-\Delta_1 - \ol{\Delta}_2$ on $E_1$ is
a fiber class of a surjective morphism $p:E_1\to\CP^1$,
 \eqref{rest001} implies that if $|\ms L_1|$ has a base point,
then the base locus has to consist of fibers of the morphism $p$.
Further, any fiber of $p$ intersects  a member of $|\ms O_{E_1}(1,0)|$ from  intersection number.
On the other hand, by \eqref{acs2}, for any smooth member $S\in |F|$,
the restriction map 
$H^0(Z,2F)\to H^0(S,2F)$ is surjective.
Moreover through the birational morphism $\mu_1:Z_1\to Z$ this map can be
naturally identified with the restriction map
$H^0(Z_1,\mathscr L_1)\to H^0(S,\mathscr L_1)$,
where we are using the same letter to mean the strict transform
of the original $S\subset Z$ into $Z_1$.
Therefore the latter map is also surjective.
Since the system $|\ms L_1|_S|$ is exactly the movable part of the bi-anticanonical system of $S$
from the definition of $\ms L_1$, and since the movable part is base point free
by Lemma \ref{lemma:SI} (ii), 
this means that for any smooth member $S\in |F|$,
the system $|\ms L_1|$ does not have a base point on the 
strict transform $S\subset Z_1$.
Therefore since fibers of the above morphism $p$ intersects
$S$ as above, the system $|\ms L_1|$ 
cannot have a base point on $E_1$.
Hence $\Bs\,|\ms L_1|=\emptyset$, as required.
\proofend

\vsp
Let $\Phi_1$ be the morphism associated to 
the free system $|\ms L_1|$ on $Z_1$.
From the construction $\Phi_1$ factors as $\Phi\circ\mu_1:
Z_1\to Z\to \CP^4$.
As seen in the above proof, 
on a general fiber $S$ of $f_1$, the morphism $\Phi_1$ 
is exactly the bi-anticanonical map of $S$.
In particular it is  degree-two over a fiber plane of 
the projection $\pi: Y\to \Lmd$.
Behavior of $\Phi_1$ over the exceptional divisors of $\mu_1$ 
is described as follows:  
\begin{proposition}\label{prop:elimI2}
We have $\Phi_1(E_1) = \Phi_1 (\ol E_1) = l$,
where $l$ is the ridge of the scroll $Y$ as before.
Further $\Phi_1(E_2\sqcup\ol E_2)=\{q,\ol q\}$,
where $q$ and $\ol q$ are the intersection 
$B\cap l$ (see Proposition \ref{prop:wsing}
and a comment after it).
\end{proposition}

\proof
From the proof of Proposition \ref{prop:baseI},
the restriction map 
$H^0 (Z_1, \ms L_1) \to H^0 (E_1, \ms L_1)$
is surjective.
By \eqref{rest001}, the latter space is 2-dimensional.
Hence the image $\Phi_1(E_1)$ is  a line.
By reality, the image $\Phi_1(\ol E_1)$ is also a line.
If we write $\phi:S\to\CP^2$ for the bi-anticanonical map
of a smooth member $S\in |F|$ as before, these lines can  be naturally identified
with the 
images $\phi(C_1)$ and $\phi(\ol C_1)$ respectively,
and, these lines are identical by Lemma \ref{lemma:SI}.
Hence $\Phi_1(E_1)$ and $\Phi_1(\ol E_1)$ must be
an identical line,
and it has to be the ridge $l$ from the proof of Proposition
\ref{prop:wsing}.
On the other hand the image $\Phi_1(E_2)$ is a point by \eqref{rest001.5}. 
Hence $\Phi_1(\ol E_2)$ is also a point.
These two points have to be 
the points $q$ and $\ol q$,
because they are naturally identified
with the images $\phi(C_2)$ and $\phi (\ol C_2)$,
which are identical to $\Phi_1(E_2)$
and $\Phi_1(\ol E_2)$ respectively.
\proofend

\vsp
By using the freeness of the system $|\ms L_1|$,
we derive information about the restriction of $\Phi_1$ to (the strict transforms of) the divisors $S_i^+$ and
$S_i^-$.
Recall from \eqref{incl} that for $i\in\{1,2\}$ we have
the inclusion $\Phi_1(S_i^+\cup S_i^-)\subset P_i$.

\begin{proposition}\label{prop:elimI3}
On the divisor $S_i^+$ ($i\in\{1,2\}$), the morphism 
$\Phi_1$ is birational over the plane $P_i$.
Moreover, the image $\Phi_1(L_i)$ is a conic on the plane,
where as before $L_i$ is the intersection $S_i^+\cap S_i^-$ which is
a twistor line.
\end{proposition}

\proof
Since $\Phi_1$ is already a morphism, we  have the coincidence
$\Phi_1(S_i^+\cup S_i^-)=P_i.$
Hence from reality we have $\Phi_1(S_i^+) = \Phi_1( S_i^-)=P_i.$
Therefore since the degree of $\Phi_1:Z_1\to Y$ is two, 
we obtain the desired birationality over the plane $P_i$.
For the latter assertion, noting that $L_i$ intersects
$E_1$ transversally at one point but does not intersect $E_2$ (see Figure \ref{fig:blowups_I}), we calculate
\begin{align*}
(\ms L_1, L_i)_{Z_1} &= 
\left(E_1 + 2E_2 + \ol E_1 + 2\ol E_2, \, L_i\right) _{Z_1} \\
&= \left( E_1 + \ol E_1, \, L_i\right) _{Z_1} = 2.
\end{align*}
This directly shows that $\Phi_1(L_i)$ is a conic.
\proofend

\subsection{Elimination of the base locus -- the case of type II}
\label{ss:elimII}
Let $Z$ be a double solid twistor space on $4\CP^2$ 
which is of type II.
Then by Lemma \ref{lemma:SII}, the pullback system $|\mu_1^*(2F)|$ on $Z_1$ has
the divisor $E_1+ E_2 + \ol E_1 + \ol E_2$ as fixed components at least.
Hence similarly to the case of type I, we put
$$
\ms L_1 : = \mu_1^* (2F) - (E_1+ E_2 + \ol E_1 + \ol E_2).
$$ 
Again from the relation \eqref{rel1} this can be rewritten as 
\begin{align}
\label{L1II}
\ms L_1= f_1^*\ms O(2) + E_1+E_2+ 2E_3 + \ol E_1+ \ol E_2+ 2 \ol E_3.
\end{align}
Concerning the base locus of the system $|\ms L_1|$, in contrast with the case of type I, we have
the following.

\begin{figure}
\includegraphics{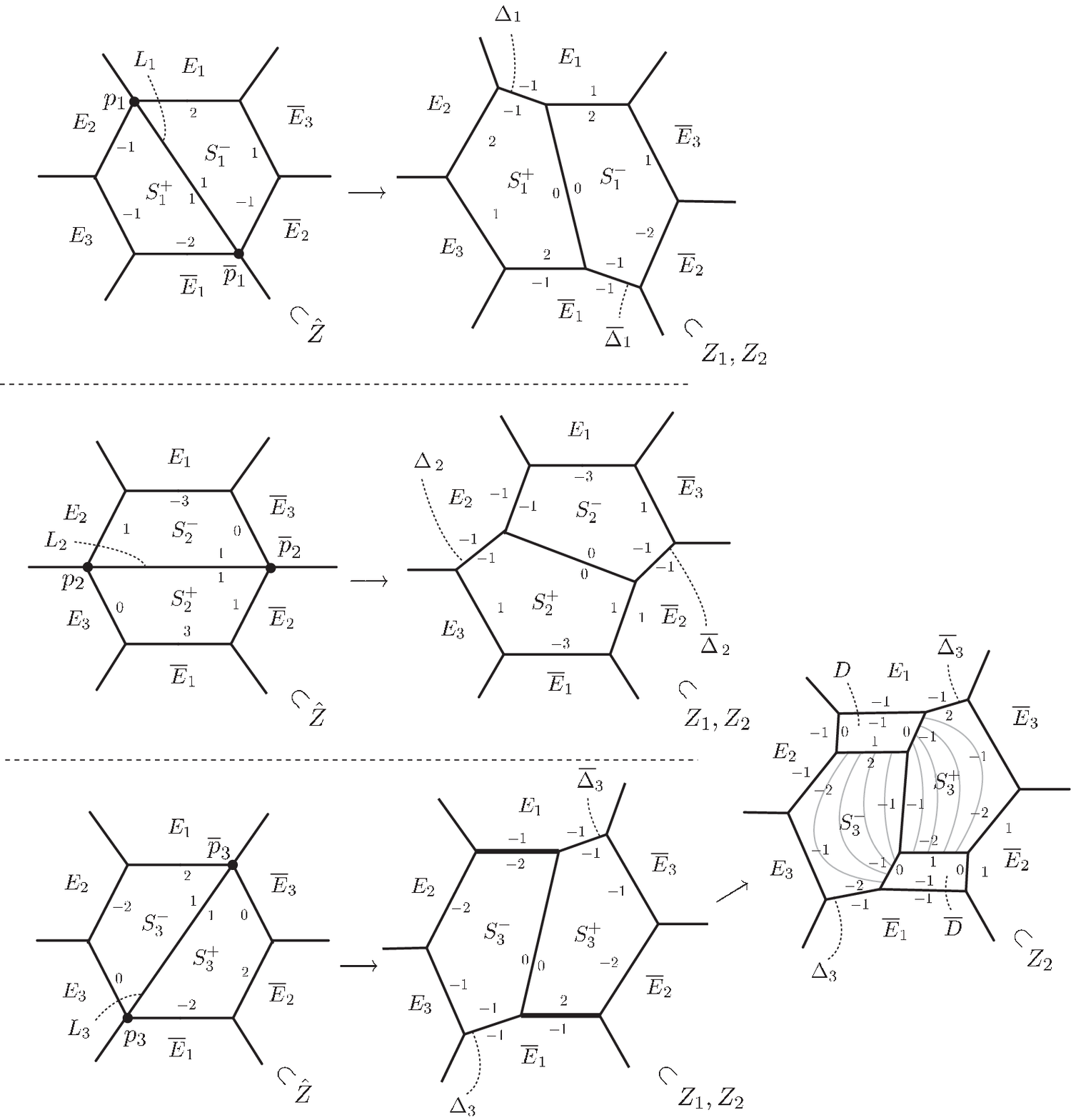}
\caption{The blowups for the case of type II.
}\label{fig:blowups_II}
\end{figure}

\begin{proposition}\label{prop:elimII}
Let $Z_1$ and $\ms L_1$ be as above.
\begin{enumerate}
\item
[(i)] The linear system $|\mathscr L_1|$ has no  fixed component.
\item
[(ii)] The two smooth rational curves 
$(S_3 ^- \cap E_1)$ and  $ (S_3 ^+ \cap  \ol E_1)$ are 
base curves of $|\mathscr L_1|$.
{\em(In Figure \ref{fig:blowups_II}, these
curves are displayed as bold segments.
Also, at this stage we do not prove that these are all base points
of $|\mathscr L_1|$, although this is actually the case.)}
\end{enumerate}
\end{proposition}

\proof
For (i), since $\Bs\,|2F| = C_1\cup C_2 \cup \ol C_1\cup\ol C_2$ on $Z$, 
it is enough to show that the exceptional divisors $E_1$ and $E_2$ are not  fixed components of $|\mathscr L_1|$.
As before we use the same letter $S$ to mean a smooth member of the original
pencil $|F|$ and its strict transform into $Z_1$.
Then again since the restriction map 
$H^0 (Z_1,\mathscr L_1) \to H^0 (S,\mathscr L_1)$ on $Z_1$
is naturally identified with the restriction map
 $H^0(Z, 2F) \to H^0 (S, 2F)$ on $Z$,
 the former map is surjective by the exact sequence \eqref{acs2}.
 Therefore, since the linear system 
$|\mathscr L_1|_S|$
is base point free by Lemma \ref{lemma:SII} (ii),
we obtain that $|\mathscr L_1|$ has no base point on the 
strict transform $S$.
Hence noting that $E_1$ and $E_2$ actually intersect $S$,
these cannot be  a fixed component of $|\mathscr L_1|$,
and we obtain (i).

For (ii), in a similar way to \eqref{int0001}--\eqref{int00026},
we can calculate the intersection numbers of 
the line bundle $\ms L_1 $ with curves on the exceptional divisor $E_i$ which generate
$H^2(E_i,\ZZ)$.
Consequently under the same notation as before, we obtain 
\begin{align}\label{rest003}
\ms L_1 |_ {E_1} \simeq \ms O(1,0) - \ol{\Delta}_3,\hsp
\ms L_1 |_ {E_2} \simeq \ms O(1,1) - {\Delta}_2,\hsp
\ms L_1 |_ {E_3} \simeq \ms O.
\end{align}
In particular, 
if we temporary write $C_1= S_3^-\cap E_1$ in this proof,
from the first one, we obtain $(\mathscr L_1 , {C_1} )_{Z_1}=-1$.
Hence $C_1$ is a base curve of $|\mathscr L_1|$.
By reality $\ol C_1 = S_3^+\cap \ol E_1$ is also a
base curve of $|\mathscr L_1|$.
\proofend

\vsp
By Proposition \ref{prop:elimII}, let $\mu_2:Z_2\to Z_1$ 
be the blowup at the base curves $(S_3 ^- \cap E_1) \sqcup (S_3 ^+ \cap  \ol E_1)$,
and $D$ and $\ol D$ the exceptional divisors
over $S_3 ^- \cap E_1$ and 
$S_3 ^+ \cap \ol E_1$ respectively.
(See the lower right in Figure \ref{fig:blowups_II}.)
From the normal bundle of the curves inside $Z_1$,
we have $D\simeq\ol D\simeq\FF_1$.
We define a line bundle over $Z_2$ by
 $$\ms L_2 := \mu_2^*\mathscr L_1 - (D + \ol D).$$
The morphism $\mu_2$ induces an isomorphism $H^0(\ms L_2)\simeq H^0(\mathscr L_1)$.
Let $f_2:Z_2\stackrel{\mu_2}\to Z_1 \stackrel{f_1}\to\CP^1$ be 
the composition morphism. Then
 the relation \eqref{L1II} is valid for $\ms L_2$ if we replace $f_1$ with $f_2$ since only $E_1\sqcup\ol E_1$
include the center of $\mu_2$, and the coefficients in 
\eqref{L1II}
of these divisors are one.
(But note that on $Z_2$ the surfaces $S_3^-$ and $ E_1$ are separated by $D$
as displayed in the lower right in Figure \ref{fig:blowups_II}.)
The next proposition implies that the blowup 
$\mu_2$ terminates an elimination
of  the anticanonical system of $Z$
in the case of type II:

\begin{proposition}\label{prop:elimII2}
The system $|\ms L_2|$ on $Z_2$ is base point free.
\end{proposition}

\proof
Again we use the same letters to mean the strict transforms into $Z_2$ and compute the restriction
$\ms L_2|_{E_i}$ for $1\le i\le 3$.
The linear system $|\ms O(1,0) - \ol{\Delta}_3|$ on $E_1$
in \eqref{rest003}
consists of the unique curve $S_3^-\cap E_1$, 
which is exactly the center of $\mu_2$.
Therefore by \eqref{rest003}, we obtain that 
$\ms L_2$ is trivial over $E_1$.
Since the center is disjoint from $E_3$,
the third one in \eqref{rest003} implies that 
$\ms L_2$ is still trivial over $E_3$.
On the other hand, noting that the divisor $E_2$ in $Z_1$ is 
blown-up at a point by $\mu_2$ and its exceptional curve
is $D\cap E_2$ 
(see the lower right in Figure \ref{fig:blowups_II}),
we obtain from the middle of \eqref{rest003} that $\ms L_2|_{E_2}$ is isomorphic to 
$\ms O(1,1) - \Delta_2 - (D\cap E_2)$.
Thus we have obtained 
\begin{align}\label{rest004}
\ms L_2 |_ {E_1} \simeq \ms O,\hsp
\ms L_2 |_ {E_2} \simeq \ms O(1,1) - {\Delta}_2
- (D\cap E_2) ,\hsp
\ms L_2 |_ {E_3} \simeq \ms O.
\end{align}
On the other hand, any exceptional divisors $E_i$ or $\ol E_i$
are not fixed components of the system $|\ms L_2|$
by Proposition \ref{prop:elimII} (i).
Therefore
by the first and the third one in \eqref{rest004},
the divisors $E_1$ and $E_3$ must be  disjoint from $\Bs\,|\ms L_2|$.
In particular, as  $\Bs\,|2F| = C_1\cup C_2 \cup \ol C_1\cup\ol C_2$,
we have  
$\Bs\,|\ms L_2| \subsetneq E_2\cup \ol E_2
\cup D\cup \ol D.$
Next we see that $(\Bs\,|\ms L_2| ) \cap E_2 = \emptyset$.
Since $|\mathscr L_1|$ has no fixed point on any smooth fiber of $f_1$
as shown in the proof of Proposition \ref{prop:elimII}, 
the same is true for $|\ms L_2|$ for any smooth fiber
of $f_2$.
On the other hand, 
from the second one in \eqref{rest004},
we obtain that 
the system $|\ms L_2|_{E_2}|$ is a pencil without 
a base point and it induces a surjective morphism 
$E_2\to\CP^1$.
Further
 the intersection $S\cap E_i$ belongs to
$|\ms O(1,0)|$ for $S\in |f_2^*\ms O(1)|$ as before
and hence intersects any fiber of 
the morphism $E_2\to\CP^1$.
These imply that $(\Bs\,|\ms L_2|)\cap E_2=\emptyset$ as claimed.

Finally we show that 
$(\Bs\,|\ms L_2| ) \cap D = \emptyset$,
which is more difficult.
For this we first notice from the lift of \eqref{L1II} to $Z_2$ and 
Figure \ref{fig:blowups_II} that 
the restriction $\ms L_2|_D$ is isomorphic to 
$\ms O_D(E_1 + E_2)$.
The intersection $E_1\cap D$ and $E_2\cap D$ are
a $(-1)$-curve and a $0$-curve in $D$ respectively
(see the lower right in Figure \ref{fig:blowups_II}).
Recalling $D\simeq\FF_1$, these imply that $\ms L_2|_D$ is isomorphic to 
the pullback of the line bundle $\ms O_{\CP^2}(1)$ by
the blowdown $D\to \CP^2$.
Therefore in order to show $(\Bs\,|\ms L_2| ) \cap D = \emptyset$, it suffices to prove that 
the rational map  associated to 
$|\ms L_2|$, for which we denote by $\Phi_2$, satisfies $\dim\Phi_2(D)=2$.
We show this by proving that the image $\Phi_2(D)$ contains two distinct 
lines.

From the construction $\Phi_2$ factors as $Z_2
\stackrel{\mu_1\circ\mu_2}\lra Z\stackrel{\Phi}\to \CP^4$.
Therefore the image of $\Phi_2$ is the scroll $Y=\pi^{-1}({\Lambda})$.
Moreover 
each fiber of the fibration $f_2:Z_2\to\CP^1$ is mapped
(by $\Phi_2$) to a fiber plane of the projection $\pi:Y\to \Lambda$
by the diagram \eqref{CD1}.
Now we show that  the image $\Phi_2(S_3^-)$ is  a line.
 By the lift of \eqref{L1II}, temporary writing $ C_i=S_3^-\cap  E_i$ for $i=2, 3$ on $Z_2$ and recalling that 
we have subtracted the divisor $D$ when defined $\ms L_2$, we 
clearly have (see the lower right in
Figure \ref{fig:blowups_II})
\begin{align}\label{fa89fw}
\ms L_2|_{S_3^-}\simeq \mathscr O_{S_3^-} (  C_2 + 2 C_3 + \Delta_3).
\end{align}
From this it is easy to see that
the system 
$|\ms L_2|_{S_3^-}|$ is a pencil without
a base point, and the associated morphism
$S_3^-\to \CP^1$ has the curves
$D\cap S_3^-$ and $\ol D\cap S_3^-$ as sections.
(In the lower right of Figure \ref{fig:blowups_II},
fibers of the morphism $S_3^-\to \CP^1$ are indicated 
as gray curves.)
In particular, if
the restriction 
$H^0(Z_2,\ms L_2)\to H^0(S_3^-, \ms L_2)\simeq\CC^2$ were not surjective, 
the system $|\ms L_2|$ would have
base points along a fiber of the morphism
$S_3^-\to \CP^1$.
But this cannot happen 
since we already know $\Bs\,|\ms L_2|\subset
D\cup \ol D$.
Therefore the above restriction map is surjective.
This implies that $(\Bs\,|\ms L_2|)\cap
S_3^-=\emptyset$ and $\Phi_2(S_3^-)$ is a line.
Hence $\Phi_2(S_3^-\cap D)$ and $\Phi_2(S_3^-\cap \ol D)$ are also the same line.
In particular $\Phi_2(D)$ contains the line
$\Phi_2(S_3^-)$.

By the real structure, $\Phi_2(S_3^+)$ is also a line,
and this can also be written
alternatively as $\Phi_2(S_3^+\cap  D)$
by taking the conjugation of $\Phi_2(S_3^-)=\Phi_2(S_3^-\cap \ol D)$.
This implies that $\Phi_2(D)$ contains 
the line $\Phi_2(S_3^+)$ too.
If the two lines $\Phi_2(S_3^+)$ and $ \Phi_2(S_3^-)$ coincide,
the image $\Phi(S_3^-\cup S_3^+)$ would be a 1-dimensional linear subspace of the plane $P_3$.
This implies that the kernel of the restriction map
$H^0(Z,2F)\to H^0(S_3^+\cup S_3^-,2F)$ is 3-dimensional,
which contradicts $h^0(F) = 2$.
Thus the two lines  $\Phi_2(S_3^+)$ and $\Phi_2(S_3^-)$
are distinct.
Hence $\Phi_2(D)$ actually contains two lines.
Therefore $\Phi_2(D)=\CP^2$, and we finally obtain 
$(\Bs\,|\ms L_2| ) \cap D = \emptyset$.
Thus we conclude $\Bs\,|\ms L_2|=\emptyset$,
as claimed.  
\proofend

\vsp
By Proposition \ref{prop:elimII2},
the rational map $\Phi_2:Z_2\to Y\subset\CP^4$
associated to $|\ms L_2|$ is a morphism. 
For restrictions of $\Phi_2$ to reducible fibers of the fibration $f_2:Z_2\to\CP^1$, we have the following interesting behavior.

\begin{proposition}\label{prop:elimII3}
Let $\Phi_2:Z_2\to Y\subset \CP^4$ be the 
morphism induced by the free system $|\ms L_2|$ as above.
\begin{enumerate}
\item[(i)] On the two reducible fibers $S_1^+\cup S_1^-$ and $S_2^+\cup S_2^-$,
$\Phi_2$ is birational over the plane
$P_1$ and $P_2$ 
respectively on each irreducible component.
Further, the image of the twistor line $L_i= S_i^+ \cap S_i^-$ ($i\in\{1,2\}$) is a conic in the plane.
\item[(ii)] $\Phi_2(S_3^+)$ and $\Phi_2( S_3^-)$  are mutually distinct  lines in the plane $P_3$.
\end{enumerate}
\end{proposition}

In terms of the original map $\Phi:Z\to\CP^4$, the above (iii) means that 
the anticanonical map contracts $S_3^+$ and $S_3^-$ to lines.
The existence of such divisors is a major difference between the case of type I.

\vsp\noindent
{\em Proof of Proposition \ref{prop:elimII3}.}
The assertion (ii) is already shown in the final part of the  proof of Proposition \ref{prop:elimII2}.
For (i), since $\Phi_2$ is a morphism by 
Proposition \ref{prop:elimII2}, 
from the commutativity $f_2 = \pi\circ\Phi_2:Z_2\to Y\to \Lmd$,
we have the coincidence $\Phi_2(S_i^+\cup S_i^-) = P_i$
for $i\in \{1,2\}$.
Then since $\Phi_2:Z_2\to Y$ is of degree two,
it follows that $\Phi_2$ is birational over
the plane $P_i$ on each $S_i^+$ and $S_i^-$
for $i\in\{1,2\}$.
\proofend

\vsp
The morphism $\Phi_2$ maps
the exceptional divisors of the birational morphism
$Z_2\to Z$ as follows.
\begin{proposition}\label{prop:elimII4}
Let $\Phi_2:Z_2\to Y\subset\CP^4$ be the morphism associated to $|\ms L_2|$ as above.
\begin{enumerate}
\item[
(i)] $\Phi_2(E_1)=\Phi_2(\ol E_3), \Phi_2(\ol E_1)=\Phi_2( E_3)$  and 
$\{\Phi_2(E_1),\Phi_2(\ol E_1)\}=\{q,\ol q\}$,
where $q$ and $\ol q$ are the two points $B\cap l$ as before.
\item[(ii)] $\Phi_2(E_2)=\Phi_2(\ol E_2)=l$.
\item[(iii)] $\Phi_2(D)=\Phi_2(\ol D)=P_3$, and $\Phi_2$ is birational on $D$
and $\ol D$.
\end{enumerate}
\end{proposition}

\proof
We first show (ii).
As  remarked in the proof of Proposition 
\ref{prop:elimII2}, the restriction $|\ms L_2 |_{E_2}|$ is a  pencil without a base point (see \eqref{rest004})
and the restriction map 
$H^0(Z_2,\ms L_2)\to H^0(E_2,\ms L_2)$ is surjective.
These mean that $\Phi_2(E_2)$ is a line.
We show that this line is exactly the ridge $l$. 
For any non-singular fiber $S$ of $f_2$, $\Phi_2|_S$ is naturally identified with the 
bi-anticanonical map of $S$ as before, and $S\cap E_2$ is exactly the component $C_2$ of the cycle $C$, 
which is mapped to a line by Lemma \ref{lemma:SII} (v).
Since the last line is independent of the choice of 
$S$,
the line must be the intersection of fiber planes
of the projection $\pi:Y \to \Lmd$.
Thus  $\Phi_2(E_2)=l$.
Hence since $l$ is real, $\ol E_2$ 
is also mapped to the same line, and we get (ii).

Since $\ms L_2$ is trivial on $E_1$ and $\ol E_3$
as in \eqref{rest004},
each of these are mapped to a point by $\Phi_2$.
As $E_1\cap \ol E_3\neq\emptyset$, these points
must coincide.
Because $\Phi_2(\ol E_3) = \Phi(\ol C_3)$ and 
the latter is exactly the point $q$ or $\ol q$ 
on $l$ as in the proof of
Proposition \ref{prop:wsing}, by Lemma \ref{lemma:SII} (iv) and (vi)
we obtain $\Phi_2(\ol E_3) \in \{q,\ol q\}$.
Hence by the real structure, we obtain the assertion (i)
of the proposition.
Finally (iii) is already shown in the final part of the proof
of Proposition \ref{prop:elimII2}.
\proofend

\subsection{Elimination of the base locus -- the case of type III}
\label{ss:elimIII}
Let $Z$ be a double solid twistor space on $4\CP^2$ 
which is of type III.
The elimination of the base locus of the anticanonical system
$|2F|$ in this case can be obtained basically in the same way to the case of type II.
So proofs are partially sketchy.

By Lemma \ref{lemma:SIII}, the pullback system $|\mu_1^*(2F)|$ has
the divisor $E_1+ E_2 + E_3 + \ol E_1 + \ol E_2 +\ol E_3$ as fixed components at least.
Hence we put
\begin{align}
\ms L_1 &: = \mu_1^* (2F) - (E_1+ E_2 +E_3 + \ol E_1 + \ol E_2
+ \ol E_3)\notag\\
&\simeq 
f_1 ^* \mathscr O(2) + E_1 + E_2 + E_3 + 2 E_4 + \ol E_1 + \ol E_2 + \ol E_3 + 2 \ol E_4,
\label{L1III}
\end{align}
where $f_1:Z_1\to\CP^1$ is the morphism
associated to $|\mu_1^*F-E|$.
We have a natural isomorphism
$H^0(2F)\simeq H^0(\mathscr L_1)$.
Then analogously to Proposition \ref{prop:elimII} we have the following

\begin{figure}
\includegraphics{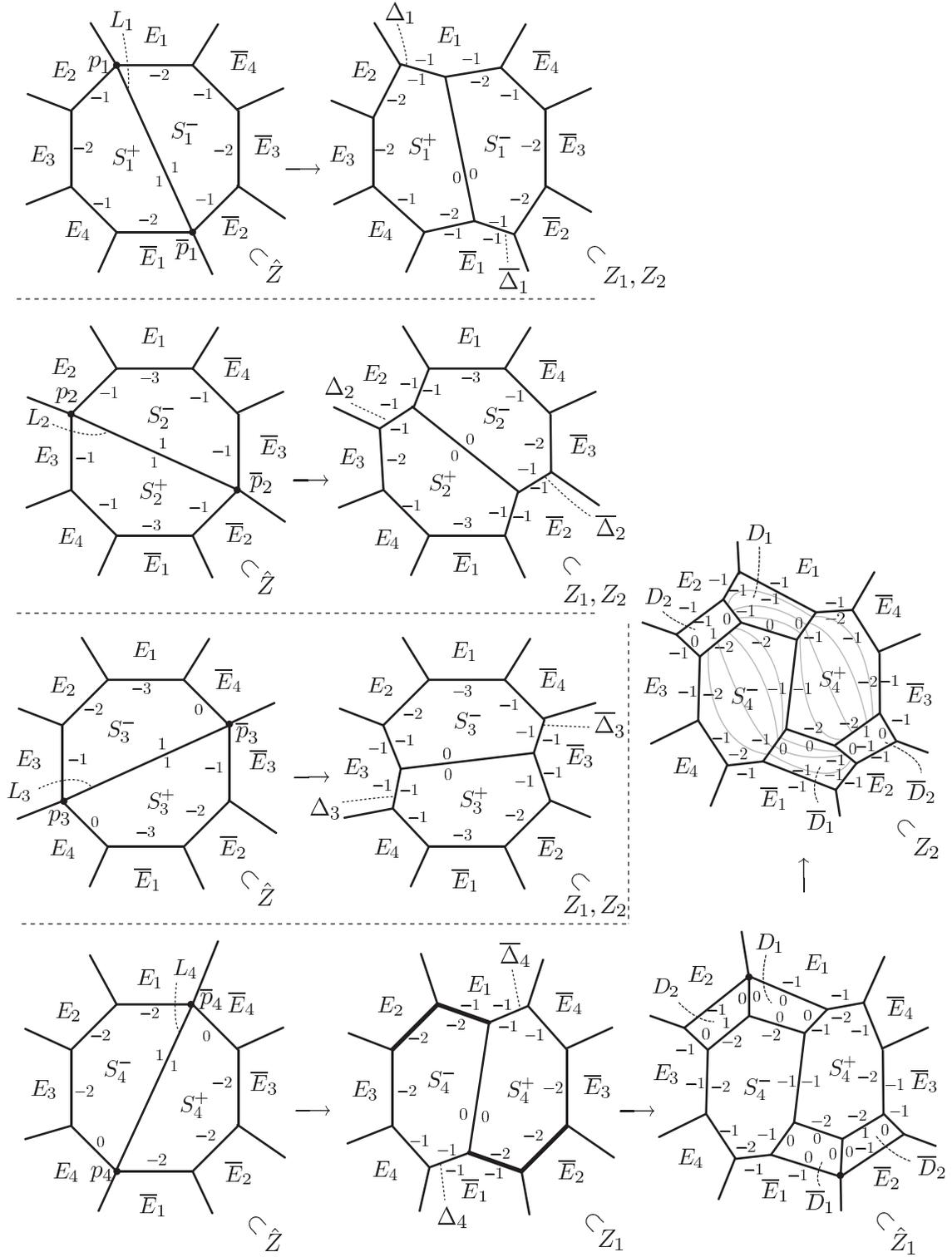}
\caption{The blowups for the case of type III.}
\label{fig:blowups_III}
\end{figure}

\begin{proposition}\label{prop:elimIII}
Let $\ms L_1$ be the line bundle over $Z_1$ as above.
\begin{enumerate}
\item[(i)] The system $|\mathscr L_1|$ does not have a fixed component.
\item[(ii)] The following four smooth rational curves 
\begin{align}\label{4curves}
S_4 ^+ \cap \ol E_1,  \hsp S_4 ^+ \cap \ol E_2,  \hsp
S_4 ^- \cap  E_1,  \hsp S_4 ^- \cap  E_2
\end{align} are 
base curves of $|\mathscr L_1|$.
{\em (Note that the first two curves intersect and the same for the last two curves; see Figure \ref{fig:blowups_III}.)}
\end{enumerate}
\end{proposition}

\proof
The assertion (i) can be seen in completely analogous way to the proof of  Proposition \ref{prop:elimII} (i).
(We just use Lemma \ref{lemma:SIII} (ii) instead of Lemma \ref{lemma:SII} (ii).)
For (ii), 
in a similar way to 
\eqref{int0001}--\eqref{int00026},
we can compute intersection numbers
of the line bundle $\ms L_1$ with
curves on the exceptional divisor $E_i$, by using \eqref{L1III}.
From them we can deduce that the restrictions of the line bundle $\ms L_1$
on $E_i$ are concretely given by 
\begin{align}\label{rest005}
\mathscr L_1|_{E_1} \simeq \mathscr O(1,0) - \ol{\Delta}_4,
\hsp
\mathscr L_1|_{E_2} \simeq  \mathscr O(1,0),
\hsp
\mathscr L_1|_{E_3} \simeq  \mathscr O(1,1) - \Delta_3,
\hsp
\mathscr L_1|_{E_4} \simeq  \mathscr O.
\end{align}
From the first and the second ones,
we obtain that  the  intersection numbers of $\mathscr L_1$ with
the four curves \eqref{4curves} are   $-1,0,-1,0$ respectively.
These imply the assertion (ii).
\proofend

\vsp
Let $\hat Z_1\to Z_1$ 
be the blowup at the base curves \eqref{4curves},
and $D_1$ and $D_2$ the exceptional divisors
over the curves $S_4 ^- \cap  E_1$ and $ S_4 ^- \cap  E_2 $ respectively.
From the normal bundles, we obtain
\begin{align}\label{Ds}
D_1\simeq\FF_0,
\hsp
D_2\simeq\FF_1.
\end{align}
The variety $\hat Z_1$ has ordinary nodes over the singular points of the
center of the last blowup,
which consist of two points in total.
Let $Z_2\to \hat Z_1$ be the small resolution of these two nodes
as displayed in Figure \ref{fig:blowups_III}
(the up-arrow in the lower right).
Let $\mu_2:Z_2\to \hat Z_1\to Z_1$ and 
$f_2:Z_2\to Z_1\to \Lmd$ be the compositions.
Then by the last small resolution $D_1$ is blown up at a point, 
while $D_2$ remains unchanged.
We again use the same letters to denote the strict transforms of these divisors into $Z_2$,
and define a line bundle $\ms L_2$ over $Z_2$ by
$$\mathscr L_2:=\mu_2^*\mathscr L_1 - (D_1+ D_2 + \ol D_1 + \ol D_2).
$$
Then similarly to Proposition \ref{prop:elimII2}, we have
\begin{proposition}\label{prop:elimIII2}
The system $|\mathscr L_2|$ on $Z_2$ is base point free.
\end{proposition}

\proof
By Proposition \ref{prop:elimIII} (i) the divisors $E_i$ and $\ol E_i$ are not  fixed components of $|\mathscr L_2|$ for any $1\le i\le 4$.
From \eqref{rest005}, we have $S_4^-\cap E_1\in |\ms L_1|_{E_1}|$
and  $S_4^-\cap E_2\in |\ms L_1|_{E_2}|$.
From these it follows that 
\begin{align}\label{rest006}
\mathscr L_2|_{E_1} \simeq \mathscr O,
\hsp
\mathscr L_2|_{E_2} \simeq  \mathscr O,
\hsp
\mathscr L_2|_{E_3} \simeq  \mathscr O(1,1) - \Delta_3 - (E_3\cap D_2),
\hsp
\mathscr L_2|_{E_4} \simeq  \mathscr O.
\end{align}
These imply that the system $|\mathscr L_2|$ does not have 
a base point on $E_1, E_2$ nor $E_4$.
Hence, as $\Bs\,|2F| = C_1\cup C_2 \cup C_3\cup \ol C_1\cup\ol C_2 \cup\ol C_3$ on $Z$,
we have
$
\Bs\,|\ms L_2| \subsetneq E_3\cup D_1\cup D_2 \cup 
\ol E_3\cup \ol D_1 \cup \ol D_2.
$
Moreover, from the third one in \eqref{rest006} the system $|\ms L_2|_{E_3}|$ is a base point free pencil, and
the intersection number with $\ms O(1,0)$ 
(on $E_3$) is 1.
Therefore, since $|\ms L_2|$ does not have 
a base point on a smooth $S\in |f_2^*\ms O(1)|$ which restricts to $\ms O(1,0)$, 
it follows that 
$|\mathscr L_2|$ does not have a base point on $E_3$ too.

Let $\Phi_2:Z_2\to Y\subset\CP^4$ be the rational map
associated to $|\mathscr L_2|$.
It remains to see that there is no base point on $D_1$ nor $D_2$.
By \eqref{Ds}, $D_2$ is biholomorphic to one point blownup
of $\CP^2$.
Further, from the lift of \eqref{L1III}
to $Z_2$,
we can deduce that $\ms L_2|_{D_2}$ is 
again isomorphic to the pullback
of $\mathscr O(1)$ by the blowdown $D_2\to \CP^2$.
Also from the lift,
 it is immediate to see from Figure \ref{fig:blowups_III}  that 
 the restrictions of $\ms L_2$
to the divisors $S_4^-$ and $D_1$  are
respectively isomorphic to the pullback of 
$\mathscr O_{\CP^1}(1)$ by  surjective morphisms
$S_4^-\to \CP^1$ and $D_1\to \CP^1$,
and that the two curves $S_4^-\cap D_2$ and $S_4^-\cap\ol D_1$ are
sections of the former morphism $S_4^-\to \CP^1$.
Furthermore the curves $S_4^+\cap D_1$ 
and $D_1\cap D_2$ are sections of
the latter morphism $D_1\to\CP^1$.
(In Figure \ref{fig:blowups_III} fibers of these morphisms
from $S_4^-$ and $D_1$ are
indicated by gray curves.)
From these we obtain that the restrictions of
$\Phi_2$ to $S_4^-$ and ${D_1}$ are exactly 
the above morphisms to $\CP^1$, and that
the images $\Phi_2(S_4^-)$ and $\Phi_2(D_1)$ are lines
in the plane $P_4$.
In particular we have $(\Bs\,|\ms L_2|) \cap D_1 
=\emptyset$.
We also have $\Phi_2(S_4^-) = \Phi_2(S_4^-\cap \ol D_1)$ and
$\Phi_2(D_1) 
= \Phi_2(D_1\cap D_2) = \Phi_2( S_4^+ \cap D_1)$.
From the former, by taking the image under
the real structure, we obtain 
$\Phi_2(S_4^+) = \Phi_2(S_4^+\cap D_1)$.
But by the latter this can be written as
$\Phi_2(D_1)$.
Therefore the above two lines in $P_4$
can be rewritten as $\Phi_2(S_4^-)$ and $\Phi_2(S_4^+)$ respectively,
which means that they are distinct.
Thus  
the image $\Phi_2(D_2)$ contains 
two distinct lines. 
Hence 
$\Phi_2|_{D_2}:D_2\to \Phi_2(D_2)$ is just the blowdown
of the $(-1)$-curve $D_2\cap E_2$.
In particular, $(\Bs\,|\ms L_2|)\cap D_2 = \emptyset$.
Hence we obtain $\Bs\,|\ms L_2| = \emptyset$.
\proofend

\vsp
By the proposition, the map $\Phi_2:Z_2 \to Y\subset\CP^4$ associated to $|\ms L_2|$
is a morphism.

\begin{proposition}\label{prop:elimIII3}
For the images of divisors on $Z_2$ under the morphism $\Phi_2$, we have
the following.
\begin{enumerate}
\item[(i)] If $i\in\{1,2,3\}$, the restrictions of $\Phi_2$ to the  divisors
$S_i^+$ and $S_i^-$ are birational over the plane $P_i$.
Moreover, the image $\Phi_2(L_i)$, where $L_i$ is the twistor line
$S_i^+\cap S_i^-$ as before, is a conic in the plane. 
\item[(ii)] $\Phi_2(S_4^+)$ and
$\Phi_2(S_4^-)$ are lines in $P_4$, which are mutually distinct.
\item[(iii)] If $i\in\{1,2,4\}$, it holds
$\{\Phi_2(E_i),\Phi_2(\ol E_i)\}=\{q,\ol q\}$, where $q$ and $\ol q$ are the intersection $B\cap l$ as before.
\item[(iv)] $\Phi_2(E_3)=\Phi_2(\ol E_3) = l$.
\item[(v)] $\Phi_2(D_1) = \Phi_2(S_4^+)$ and 
$\Phi_2(\ol D_1) = \Phi_2(S_4^-)$
(which are lines by {\rm{(ii)}}),
\item[(vi)]  $\Phi_2$ maps  
the exceptional divisors $D_2$ and $\ol D_2$ birationally to the plane $P_4$.
\end{enumerate}
\end{proposition}

\proof
The claims (ii), (iii), (v) and (vi) are already proved in 
the proof of Proposition \ref{prop:elimIII2}.
The assertions (i) and (iv) can be shown in a similar way to 
Proposition \ref{prop:elimI3}  and Proposition \ref{prop:elimII4} (ii) respectively.
\proofend

\subsection{Elimination of the base locus -- the case of type IV}
\label{ss:elimIV}
An elimination of the base locus of the 
anticanonical system on $Z$
in the case of type IV can be given in a
similar way to the case of type III.
So we just mention  key facts required in the proofs.

By Lemma \ref{lemma:SIV}, the pullback system $|\mu_1^*(2F)|$ has
the divisor $\sum_{i=1}^4(E_i + \ol E_i)$ as fixed components at least.
Let $\ms L_1$ be the line bundle obtained from the pullback $\mu_1^*(2F)$ by subtracting  this divisor.
By the relation \eqref{rel1}, we have
\begin{align}\label{basicIV}
\mathscr L_1 \simeq f_1 ^* \mathscr O(2) + 
\sum_{i=1}^4(E_i+\ol E_i) + 2 ( E_5 + \ol E_5).
\end{align}

\begin{figure}
\includegraphics{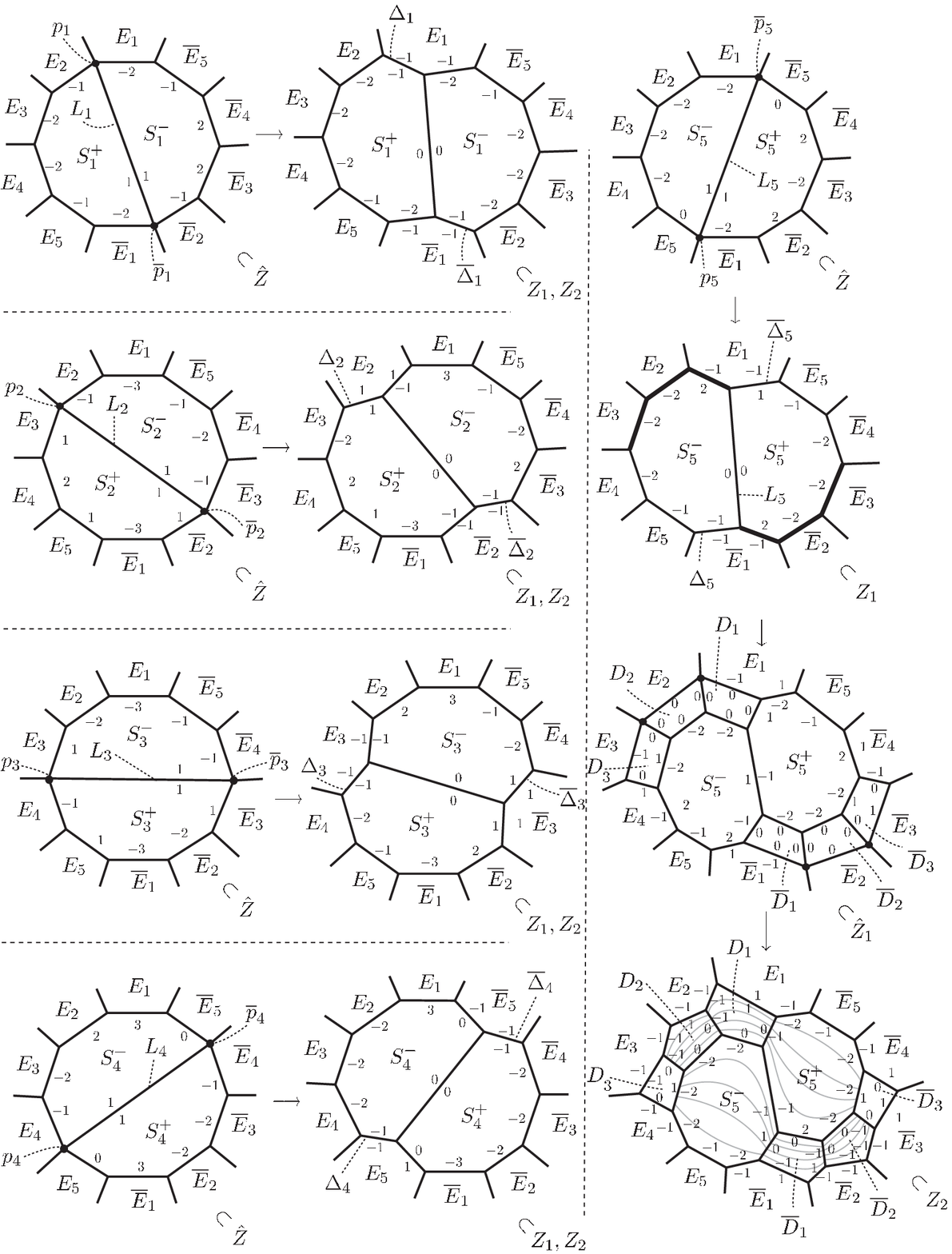}
\caption{The blowups for the case of type IV.
}\label{fig:blowups_IV}
\end{figure}

\begin{proposition}\label{prop:elimIV}
Let $\ms L_1$ be the line bundle over $Z_1$ as above.
\begin{enumerate}
\item[(i)] The system $|\ms L_1|$ has no fixed component.
\item[(ii)] The connected curves
\begin{align}\label{bcIV}
S_5^-\cap (E_1\cup E_2\cup E_3)
\hsp{\text{and}}\hsp
S_5^+\cap (\ol E_1\cup \ol E_2\cup \ol E_3)
\end{align}
are  base curves of $|\ms L_1|$.
{\em(See Figure \ref{fig:blowups_IV}.)}
\end{enumerate}
\end{proposition}

\noindent
{\em Sketch of a proof.}
The assertion (i) is easy.
For (ii), by similar computations  to \eqref{int0001}--\eqref{int00026} we can deduce 
\begin{multline}\label{888}
\ms L_1 |_{E_1} \sim \ms O(1,0) - \ol\Delta_5,
\hsp
\ms L_1 |_{E_2} \sim \ms O(1,0),
\hsp
\ms L_1 |_{E_3} \sim \ms O(1,0),\\
\ms L_1 |_{E_4} \sim \ms O(1,1) - \Delta_4,
\hsp
\ms L_1 |_{E_5} \sim \ms O.
\end{multline}
From these we obtain
\begin{align}
\left(\ms L_1,\, S_5^-\cap E_1\right)_{Z_1} = -1,
\hsp
\left(\ms L_1,\, S_5^-\cap E_2\right)_{Z_1} = 
\left(\ms L_1,\, S_5^-\cap E_3\right)_{Z_1} = 0.
\end{align}
These imply the assertion (ii).
\proofend

\vsp
Let $\hat Z_1\to Z_1$ be the blowup at
the curves \eqref{bcIV}.
Let $D_1,D_2$ and $D_3$ be the exceptional divisors
over the curves $S_5^-\cap  E_1$,
$S_5^- \cap  E_2$ and $S_5^-\cap  E_3$
respectively (see Figure \ref{fig:blowups_IV}).
We have 
\begin{align}
D_1 \simeq D_2 \simeq \FF_0,
\hsp
D_3\simeq \mathbb \FF_1.
\end{align}
The variety $\hat Z_1$ has ordinary double 
points over the singular points of the curves
\eqref{bcIV}.
These consist of four points in total
and are indicated as dotted points in Figure \ref{fig:blowups_IV}.
We take small resolutions of these 
ODP-s in a unique way that each of $E_2$ and $E_3$
(and hence $\ol E_2$ and $\ol E_3$ also) is blown-up at one point
(Figure \ref{fig:blowups_IV}).
Let $Z_2\to\hat Z_1$ be the resulting birational morphism.
Let $\mu_2:Z_2\to \hat Z_1\to Z_1$
be the composition, and define a line bundle  
$\ms L_2$ over $Z_2$ to be the line bundle
  obtained from the pullback $\mu_2^*\ms L_1$
by subtracting the divisor $\sum_{i=1}^3(D_i+\ol D_i)$,
where we are using the same letters to denote
the divisors in $\hat Z_1$ and their strict
transforms into $Z_2$.
Then we have:
\begin{proposition}
The linear system $|\ms L_2|$ over $Z_2$ is base point free.
\end{proposition}
\noindent
{\em Outline of a proof.}
Since the restrictions of $\ms L_1$ 
to the divisors $E_i$, $1\le i\le 3$,
are exactly equivalent to the curves \eqref{bcIV} by \eqref{888},
the transformation $\ms L_2$ is 
trivial over the strict transforms
of these divisors into $Z_2$.
Therefore $\ms L_2|_{E_i}\simeq\ms O$\,
for $i\in \{1,2,3, 5\}$.
By Proposition \ref{prop:elimIV} (i), this means that 
$(\Bs\,|\ms L_2| ) \cap (E_1\cup E_2\cup E_3
\cup E_5) = \emptyset$.
Hence $\Bs\,|\ms L_2|\subsetneq E_4\cup D_1\cup
D_2\cup D_3\cup \ol E_4\cup \ol D_1\cup
\ol D_2\cup \ol D_3$.
Further we have
\begin{align}
\ms L_2|_{E_4} \simeq \ms O(1,1) - \Delta_4 
-(E_4\cap D_3).
\end{align}
Since $(\Bs\,|\ms L_2|)\cap S=\emptyset$
for a smooth $S\in |f_2^*\ms O(1)|$
and $S\cap E_4\in |\ms O(1,0)|$, 
this means that 
$(\Bs\,|\mathscr L_2|)\cap E_4 = \emptyset$.
Also, if we temporary write curves on $S_5^-\,(\subset Z_2)$ by
$\ol C_1 =  S_5^-\cap \ol E_1$ and
$C_i= S_5^-\cap E_i$ for $i\in\{4,5\}$,
by using the lift of the relation \eqref{basicIV} to $Z_2$,
we have
\begin{align}
\ms L_2|_{S_5^-} \simeq \ms O_{S_5^-} 
( C_4 + 2 C_5 + \ol C_1).
\end{align}
From this, the system $|\ms L_2|_{S_5^-}|$ is a base 
point free pencil, whose members are not contained in 
$D_1\cup D_2\cup D_3$.
(Members of this pencil on $S_5^-$ are illustrated
as gray curves in the lower right of Figure \ref{fig:blowups_IV}.)
Therefore since 
$\Bs\,|\ms L_2| \subset
D_1\cup D_2 \cup D_3 \cup 
\ol D_1\cup \ol D_2 \cup \ol D_3$,
we obtain $(\Bs\,|\mathscr L_2|)\cap S_5^- = \emptyset$.
At the same time, if $\Phi_2:Z_2\to Y\subset\CP^4$
denotes the rational map induced by $|\ms L_2|$,
  we obtain that
$\Phi_2(S_5^-)$ is a line.
In a similar way to the proofs of Propositions \ref{prop:elimII2}
and \ref{prop:elimIII2},
we can show that 
the two lines $\Phi_2(S_5^-)$ and $\Phi_2(S_5^+)$ are distinct
by using $h^0(F) = 2$.
Moreover from th above lifted expression for $\ms L_2$,
we see that the system $|\ms L_2|_{D_i}|$,
$i\in\{1,2\}$, is a base point free pencil,
and hence induces a surjective morphism $D_i\to\CP^1$.
(Fibers of these morphisms are indicated by
gray curves in the lower right of Figure \ref{fig:blowups_IV}.)
Furthermore  the curves
$D_1\cap S_5^+$ and $D_2\cap D_1$ 
are sections of the two morphisms.
This implies that if the restriction map
$H^0(Z_2,\ms L_2) \to H^0(D_i,\ms L_2)$ is not surjective
for $i=1$ or $i=2$, then $(\Bs\,|\ms L_2|)\cap S_5^+\neq
\emptyset$.
This contradicts   $(\Bs\,|\mathscr L_2|)\cap S_5^- = \emptyset$ which was already obtained.
Hence the last restriction map is surjective for $i\in\{1,2\}$.
Therefore $\Bs\,|\mathscr L_2|\subset D_3\cup \ol D_3$.
At the same time we also obtain that $\Phi_2(D_1)$ and $\Phi_2(D_2)$ are
lines, both of which are identical to $\Phi_2(S_5^+)$.
Finally the lifted expression for $\ms L_2$ means that 
$\ms L_2|_{D_3}$ is of the form $\eta^*\ms O(1)$,
where $\eta:D_3\simeq\FF_1\to \CP^2$ is a blowdown.
From the above argument the image $\Phi_2(D_3)$
contains two distinct lines $\Phi_2(S_5^-)$ and 
$\Phi_2(S_5^+)$.
This means that $\Phi_2(D_3)$ is 2-dimensional.
Hence $(\Bs\,|\mathscr L_2|)\cap D_3=\emptyset$ and
$\Phi_2|_{D_3}$ can be identified with $\eta$.
Hence $(\Bs\,|\ms L_2|)\cap D_3 = \emptyset$,
and we have obtained 
$\Bs\,|\ms L_2|  = \emptyset$.
\proofend

\vsp
All the following properties concerning behavior of
the morphism $\Phi_2$ on divisors on $Z_2$ are
in effect shown in the above proof:
\begin{proposition}\label{prop:elimIV2}
The morphism $\Phi_2:Z_2\to Y\subset\CP^4$ 
associated to the free system $|\ms L_2|$ satisfies
the following properties.
\begin{enumerate}
\item[(i)]
If $1\le i\le 4$, the restrictions $\Phi_2|_{S_i^+}$ 
and $\Phi_2|_{S_i^-}$ are birational over
the plane $P_i$.
Moreover the image $\Phi_2(L_i)$ is a conic in the plane.
\item[(ii)] The images
$\Phi_2(S_5^+)$ and $\Phi_2(S_5^-)$
are mutually distinct lines in the plane
$P_5$. 
\end{enumerate}
\end{proposition}

\begin{proposition}
The morphism $\Phi_2$ fulfils
the following properties.
\begin{enumerate}
\item[(i)]
$\Phi_2$ maps the two connected divisors
$\ol E_5\cup E_1\cup E_2\cup E_3$
and 
$E_5\cup \ol E_1\cup\ol E_2\cup \ol E_3$
to the two points $q$ and $\ol q$ on the ridge $l$.
\item[(ii)] 
$\Phi_2(E_4)=\Phi_2(\ol E_4) = l$.
\item[(iii)]
$\Phi_2(D_1) = \Phi_2(D_2) = \Phi_2(S_5^+)$
and 
$\Phi_2(\ol D_1) = \Phi_2(\ol D_2) = \Phi_2(
S_5^-)$.
(These are distinct lines by Proposition \ref{prop:elimIV2} (ii).)
\item[(iv)] The restrictions $\Phi_2|_{D_3}$
and $\Phi_2|_{\ol D_3}$ are birational 
over the plane $P_5$,
and these can be identified with 
the blowdown $\FF_1
\to \CP^2$.
\end{enumerate}
\end{proposition}

\section{Defining equation of the quartic hypersurface}
\label{s:de}
The goal in this section is to obtain 
the defining equation of the quartic hypersurface
in $\CP^4$ which cuts out the branch 
divisor of the anticanonical map 
$\Phi:Z\to Y\subset\CP^4$ of the
double solid 
twistor spaces on $4\CP^2$ for arbitrary type. 
\subsection{Reducible members of the anticanonical system}
\label{ss:red}
In this subsection we show that
on twistor spaces of double solid type
which are of type I, II or III, there exist particular 
reducible members of the anticanonical system.
(For the case of type IV there does not
exist such a member.)
These anticanonical divisors will bring a strong constraint for the form of the 
defining equation of the quartic hypersurface.
We show the existence of these divisors by explicitly giving 
the chern class in which their irreducible 
component belongs.
For this, in order to obtain concrete basis of
the cohomology group $H^2(Z,\ZZ)$,
we make use of the smooth member $S$ in the 
pencil $|F|$. 
Recall that 
as explained in Section \ref{s:clsf} there is
a birational morphism
 $\epsilon:S\to \FF_0$.

First we discuss the case of type I.
For this case as in Section \ref{ss:k32} we
can suppose that  the morphism
$\epsilon$ blows up only smooth points of 
the  $(2,2)$-curve 
$C_1 + C_2 + \ol C_1 + \ol C_2$
on $\FF_0$.
Further, 
$\epsilon$ blows up three points on $C_1
\in|\ms O(1,0)|$ 
and one point on $C_2\in |\ms O(0,1)|$.
Let $e_i$,  $i\in\{1,2,3\}$, be the exceptional curve intersecting the curve
$C_1\subset S$, 
and $e_4$ the exceptional curve 
intersecting the curve $C_2\subset S$.
We have
$(e_i, C_1)_S=1$ for $i\in\{1,2,3\}$ 
and $(e_4, C_2)_S=1$.
Next let $\{\alpha_i\set 1\le i\le 4\}$ be  orthonormal basis of $H^2(4\mathbb{CP}^2,\mathbb Z)$ determined from $\{e_1,e_2,e_3,e_4\}$
by the condition that,
letting $t:Z\to 4\mathbb{CP}^2$ be the twistor fibration, $t^*\alpha_i|_S=e_i-\ol e_i$
holds in $H^2(S,\mathbb Z)$.
Then the collection 
$\{ F, t^*\aaa_1,t^*\aaa_2,t^*\aaa_3,t^*\aaa_4\}$
spans the real cohomology group 
$H^2(Z,\RR)\simeq\RR^5$.
Under these settings we have the following proposition.

\begin{proposition}\label{prop:redivisorI}
For each $i\in \{1,2,3\}$ (not including $4$), 
the linear systems $|F+t^*\alpha_i|$ and
$|F-t^*\alpha_i|$ consist of a single element
respectively.
Moreover, all these six divisors are irreducible.
\end{proposition}

\proof
In this proof for simplicity we write $\alpha_i$ for $t^*\alpha_i$.
By the real structure, it is enough to prove the claim for the latter system $|F-\alpha_i|$.
Fix any $i\in\{1,2,3\}$.
We then have 
\begin{align}
(F -   \alpha_i ) | _S &= K_S^{-1} - ( e_i - \ol e_i ) 
\notag\\
&=
\epsilon^* \mathscr O(2,2) - \sum _ {j=1} ^4 ( e_j + \ol e_j ) -  ( e_i - \ol e_i )
\notag\\ 
&=\Big( \epsilon^* \mathscr O(1,0) - \sum_{1\le j\le3} \ e_j \Big) +
 \Big( \epsilon^* \mathscr O(1,2)   -  e_i - e_4 - \ol e_4 - \sum _ {1\le j\le 3,\, j\neq i }  \ol e_j \Big).
 \label{reddd}
\end{align}
From this, the intersection number of 
 $(F - \alpha_i)|_S$ with the curve
$C_1\sim \epsilon^* \mathscr O(1,0) - \sum _ {j=1}^3  e_j$ 
 is easily 
 computed to be $(-2)$.
Hence $ C_1$ is a fixed component of $|(F -   \alpha_i ) | _S|$.
Further counting dimension, the remaining system $| \epsilon^* \mathscr O(1,2)   -  e_i - e_4 - \ol e_4 - \sum _ {j=1,\, j\neq i } ^ 3 \ol e_j | $ consists of a single member, which is the strict transform of a $(1,2)$-curve passing through the 5 points 
$\epsilon( e_i), 
\epsilon(e_4), \epsilon(\ol e_4)$ and 
$\epsilon(\ol e_j)$ with $1\le j\le 3$ and $j\neq i$.
Thus we have $h^0 ( ( F - \alpha_i ) | _S )  = 1$, and also that if $X\in |F-\aaa_i|$ denotes the unique member of the system, then 
the intersection $X|_S$ has a component
which is not contained the cycle $C$.

Next let $s\in H^0(F)$ be an element such that $(s)=S$,
and for $i\in\{1,2,3\}$ we consider the obvious exact sequence 
\begin{align}\label{ses:40}
 0 \,\lra\,
 \mathscr O_Z(-\alpha_i) \,\stackrel{\otimes s}{\lra}\, F \otimes \mathscr O_Z(- \alpha_i) \,\lra\, ( F- \alpha_i) |_S \,\lra\, 0.
\end{align}
By Riemann-Roch formula we have $\chi(\mathscr O_Z(-\alpha_i))=0$ 
and by Hitchin vanishing 
\cite{Hi80} we have $H^2(\mathscr O_Z(-\alpha_i))=0$.
Also $H^0(\mathscr O_Z(-\alpha_i))=H^3(\mathscr O_Z(-\alpha_i))=0$ by a trivial reason.
Hence we obtain $H^1 (\mathscr O_Z(-\alpha_i))=0$.
Thus by the cohomology exact sequence of \eqref{ses:40} we obtain $h^0(F-\alpha_i) = 
h^0( (F-\alpha_i)|_S)$.
Therefore we obtain $h^0(F-\alpha_i) = 1$ for $i\in\{1,2,3\}$.

To see the irreducibility of the unique member of $|F-\aaa_i|$,
let $X$ be the member.
Then as $X$ is of degree-two, if $X$ is reducible,
it consists of two degree-one divisors.
As $S_1^{\pm}$ and $S_2^{\pm}$ are all degree-one divisors
on $Z$, the intersection $X|_S$ must be contained in the
cycle $C$.
However, as we have seen above from \eqref{reddd}
the intersection $X|_S$ has a component
which is not contained in the cycle $C$.
Therefore $X$ cannot be 
the sum of degree-one divisors and hence irreducible.
\proofend

\vsp
Next we consider the case of type II.
In this case as in the proof of 
Proposition \ref{prop:k32} we decompose
the birational morphism 
 $\epsilon: S\to \FF_0$ 
in a way that only
$\epsilon_1$ blows up a singular point
of the curve
$C_1 + C_2 + \ol C_1 + \ol C_2$ on $\FF_0$.
We can suppose that the pair of blowup points
of $\epsilon_1$ is $\{\ol C_2\cap C_1,
C_2\cap \ol C_1\}$.
This means that the $(-1)$-curves
$\ol C_3$ and $C_3$ are
the exceptional curves of $\epsilon_1$.
Then under the self-intersection numbers \eqref{II},
$\epsilon_4\circ\epsilon_3\circ\epsilon_2$ blows up two points on $C_1$
and one point on $C_2$, all of which
are smooth points of the cycle.
Let $e_1$ and $e_2$ be the exceptional curves
of the blown-up points on $C_1$,
and $e_4$ the exceptional curve
of the blownup point on $C_2$.
Let $e_3$ be the exceptional curve of
the point $\ol C_2\cap C_1$.
Then we have
\begin{align}\label{int1}
e_3 = \ol C_3, \,\, (e_1,C_1) = (e_2,C_1) = 1, \,\, (e_4,C_2)=1.
\end{align}
Next let $\{\alpha_1,\alpha_2,\alpha_3,\alpha_4\}$
be orthonormal basis of $H^2(4\CP^2,\ZZ) $ which are uniquely determined from 
$\{e_1,e_2,e_3,e_4\}$ via
the twistor fibration $t:Z\to 4\CP^2$ as in the case of type I. 

\begin{proposition}\label{prop:redivisorII}
For each $i\in \{1,2\}$,
the linear systems $|F+t^*\alpha_i|$ and
$|F-t^*\alpha_i|$ consist of a single element
respectively.
Moreover, all these four divisors are irreducible.
\end{proposition}

\proof
We again write $\aaa_i$ to mean $t^*\aaa_i$.
We first show that $h^0((F-\alpha_1)|_S) =1$.
If we write $(a,b):=\epsilon^*\mathscr O(a,b) \in H^2(S,\ZZ)$,
from \eqref{003} and \eqref{int1}   we have the following relations in $H^2(S,\ZZ)$:
\begin{align*}
C_1  = (1,0)- e_1-e_2-e_3,\quad
C_2  = (0,1) - \ol e_3 - e_4.
\end{align*}
By using the former, we compute as
\begin{align*}
(F - \alpha_1, C_1)_Z & = 
\left( K_S^{-1} - (e_1 - \ol e_1),\,   (1,0)- e_1-e_2-e_3 \right)_S\notag\\
& = \left( (2,2) - 2e_1 - e_2 - e_3 - e_ 4 - \ol e_2 - \ol e_3 - \ol e_4 ,\,  (1,0)- e_1-e_2-e_3 \right)_S \notag \\
& = 2 + (-2 - 1 - 1) = -2.
\end{align*}
Hence the curve $C_1$ is a fixed component of
the system $|(F-\alpha_1)|_S|$.
In a similar way we further find
$((F-\alpha_1)|_S - C_1,\, C_2)_S = -1,$
meaning that $C_2$ is also a fixed component of the same system.
We then have
\begin{align}\label{subtract1}
(F-\alpha_1)|_S - C_1 - C_2 = (1,1) - e_1 - \ol e_2 - \ol e_4.
\end{align}
Now, since the points $\epsilon(e_1),
\epsilon(\ol e_2)$ and $\epsilon(\ol e_4)$
are smooth points on the cycle 
$C_1 + C_2 + \ol C_1 + \ol C_2$ on $\FF_0$,
there exists a unique member of the 
linear system  \eqref{subtract1},
and it is not contained in the cycle $C$.
Thus we get $h^0((F-\alpha_1)|_S) =1$, as claimed.
In the same manner, we obtain 
$((F - \alpha_2)|_S,\, C_1)_S<0$, $((F-\alpha_2)|_S - C_1,\,  C_2)_S <0$
and  $(F-\alpha_2)|_S - C_1 - C_2 = (1,1) - \ol e_1 - e_2 -\ol e_4$,
which again imply $h^0((F-\alpha_2)|_S) =1$.
We also obtain that the unique member of $|(F-\alpha_2)|_S|$ is not
contained in the cycle $C$.

Thus we have obtained $h^0((F-\aaa_i)|_S)=1$
for $i\in\{1,2\}$.
Then the same argument using the exact sequence
\eqref{ses:40} and Hitchin vanishing theorem
implies $h^0(F-\aaa_i)=1$.
The irreducibility of the unique member of
$|F-\aaa_i|$ $(i\in\{1,2\}$) can be
also obtained in a similar way to
the proof of Proposition \ref{prop:redivisorI},
by using the fact that 
the unique member of the system $|F-\aaa_i|_S|$ has a component
which is not contained in the cycle $C$.
\proofend

\vspace{2mm}
Finally we consider the case of type III.
For this case, we can choose the 
exceptional classes of 
$\epsilon:S\to \FF_0$ in such a
way that the irreducible components
of the anticanonical cycle are given by
\begin{align}\label{curvesIII}
C_1 = (1,0) - e_1 - e_2 - e_3, \,\,
C_2 = (0,1) - \ol e_2 - e_4, \,\,
C_3 = \ol e_2 - \ol e_3,\,\, C_4 = \ol e_3
\end{align}
in $H^2(S,\ZZ)$.
As before 
let $\{\aaa_1,\aaa_2,\aaa_3,\aaa_4\}$ be 
orthonormal basis of $H^2 ( 4\CP^2,\ZZ)$
determined from $\{e_1,e_2,e_3,e_4\}$ 
via the twistor fibration $t:Z\to4\CP^2$.

\begin{proposition}\label{prop:redivisorIII}
For the case of type III, under the
above choice of the basis of $H^2 ( 4\CP^2,\ZZ)$,
the linear systems $|F+t^*\alpha_1|$ and
$|F-t^*\alpha_1|$ consist of a single element
respectively.
Moreover, both of these divisors are irreducible.
\end{proposition}

\noindent
{\em Outline of a proof.}
By computing intersection numbers,
utilizing \eqref{curvesIII}, we can show  that $|( F-\alpha_1)|_S|$ has $C_1+C_2+C_3$
as fixed components. 
This system has a unique member, 
and the member has a component
which is not contained in the cycle $C$.
Then by the same argument for 
Proposition \ref{prop:redivisorI} 
using the exact sequence \eqref{ses:40} we deduce $h^0(F-\alpha_1)= 1$.
The irreducibility of the unique
member of the system $|F-\aaa_1|$ can 
be readily derived as in the proof
of Proposition \ref{prop:redivisorI}.
\proofend

\subsection{Finding double curves on the branch divisor}
\label{ss:double}
As before let $\Phi:Z\to Y\subset \CP^4$ be
the anticanonical map of a twistor space on $4\CP^2$
of double solid type,
and $B$  the branch divisor of $\Phi$.
By Theorem \ref{thm:DS},
$B$ is an intersection of the scroll $Y$
with a quartic hypersurface in $\CP^4$.
Let $D$ be an irreducible divisor on $B$.
We say that a curve $\ms C$ in $B$ is a 
{\em double curve with respect to the divisor $D$}
if the restriction $B|_D$ is a non-reduced curve
of multiplicity two.
Geometrically this means that
$D$ touches $B$ along the curve $\ms C$.
In this subsection, 
we will find {\em five} double curves on $B$,
regardless of the type of $Z$.
These double curves will bring a strong constraint
for the form of the defining equation of the quartic hypersurface.

The double curves we are going to consider consist of 
two groups. We begin with the first one.
We use the notations in 
Definition \ref{def:term2}.
By the commutative diagram \eqref{CD1},
for each reducible member 
$S_i^++S_i^-\in |F|$, we have 
$\Phi^{-1}(P_i) = S_i^+\cup S_i^-$.
If $Z$ is of type I, 
by Proposition \ref{prop:elimI3}, 
the anticanonical map $\Phi$ gives 
by restriction the birational maps
$S_i^+\to P_i$ and
$S_i^-\to P_i$ for $i\in \{1,2\}$.
Further, by the same proposition, 
the image $\ms C_i:=\Phi(L_i)$
of the twistor line $L_i$ is a conic in the plane.
Since $B|_{P_i}$ is a curve of degree four,
these imply that the conic $\ms C_i$ 
is a double curve with respect to the plane $P_i$
for $i\in \{1,2\}$.
Thus in the case of type I we have obtained two double curves $\ms C_1$ and $\ms C_2$.
For the case of type II, 
by Proposition \ref{prop:elimII3} this time,
the images $\ms C_1=\Phi(L_1)$ and 
$\ms C_2=\Phi(L_2)$ of the twistor lines
are conics in the planes $P_1$ and $P_2$ respectively,
and these conics have to be double curves with 
respect to the planes by the same reason to the
case of type I.
Similarly, for the case of type III 
(resp.\,type IV), 
by Proposition \ref{prop:elimIII3}
(resp.\,Proposition \ref{prop:elimIV2}),
$\ms C_i=\Phi(L_i)$ is a conic in $P_i$
for $i\in \{1,2,3\}$ (resp.\,$i\in\{1,2,3,4\}$) 
and these are double curves with respect to $P_i$.
Regardless of the types, we call these double curves 
 as {\em double conics}.
We have, for each double conic $\ms C_i$
\begin{align}\label{dcrest}
\ms C_i = B\cap P_i,
\end{align}
as sets.

On the other hand, for the case of type II, III
and IV, by Propositions \ref{prop:elimII3}, 
\ref{prop:elimIII3} and \ref{prop:elimIV2}
respectively,
the remaining reducible member
of $|F|$ in each case is mapped by $\Phi$ to two lines on the plane.
Namely if we put
\begin{align}
\ms C_3:= \Phi(S_3^+\cup S_3^-)\subset P_3
&\quad {\text{for the case of type II}},
\label{sdcII}\\
\ms C_4:= \Phi(S_4^+\cup S_4^-)\subset P_4
&\quad {\text{for the case of type III}},
\label{sdcIII}\\
\ms C_5:= \Phi(S_5^+\cup S_5^-)\subset P_5
&\quad {\text{for the case of type IV}},
\label{sdcIV}
\end{align}
then all these are a union of two lines.
From the computations in Section \ref{s:acs},
the normal bundle of $S_i^+$ and
$S_i^-$ in the manifold $Z_2$ is degree $-2$
along fibers of the restriction
$\Phi_2|_{S_i^+\cup S_i^-}:S_i^+
\cup S_i^-\to \ms C_i$.
This means that the above curves
$\ms C_3$ in the case of type II,
$\ms C_4$ in the case of type III,
and $\ms C_5$ in the case of type IV, are 
double curves on the branch divisor $B$ 
with respect to the planes $P_3$,
 $P_4$ and  $P_5$
respectively.
We call these double curves
\eqref{sdcII}--\eqref{sdcIV} as 
{\em splitting double conic}.
In the case of type I, there is no
splitting double conic.

Next we present double curves in the second group.
First for the case of type I,
by Proposition \ref{prop:redivisorI}, 
there exist exactly three reducible 
 anticanonical divisors on $Z$, 
each of which consists of two irreducible 
components.
Let 
$X_3+\ol X_3, X_4 + \ol X_4$ and $X_5 + \ol X_5$ be these three divisors.
Similarly, for the case of type II,
let $X_4+\ol X_4$ and $X_5 + \ol X_5$ be 
the two reducible anticanonical divisors 
found in Proposition \ref{prop:redivisorII},
and
for the case of type III,
let $X_5 + \ol X_5$ be 
the reducible anticanonical divisors 
found in Proposition \ref{prop:redivisorIII}.
Let  $X_i+\ol X_i$ be any one of these
reducible anticanonical divisors,
and $H_i\subset\CP^4$  the hyperplane
which satisfies
 $\Phi^{-1}(H_i) = X_i + \ol X_i$.
Then the curve
$X_i\cap \ol X_i$ is a ramification divisor
of the restriction $\Phi:\Phi^{-1}(H_i)
\to H_i\cap Y$.
Therefore the image
$$
\ms C_{i}:= \Phi(X_i\cap \ol X_i)
$$
is a double curve with respect
to the cone $H_i\cap Y$.
Then because $B$ is a cut of $Y$ by a quartic hypersurface, we have $\mathscr C_{i}\in |\mathscr O_{Y\cap H_{i}}(2)|$ , where $\mathscr O_{Y\cap H_i}(2)
:= \mathscr O_{\mathbb{CP}^4}(2)|_{Y\cap H_i}$.
Namely, these double curves 
are intersection of the quadratic cone $Y\cap H_i$ with a
quadric in $H_i=\mathbb{CP}^3$.
Therefore these double curves are of degree 4 in $\mathbb{CP}^4$.
So in the following we call these double curves as {\em double quartic curves}.
From our choice of the hyperplanes in $\CP^4$, the double quartic curves can be written as, as sets,
\begin{align}\label{dc3}
\mathscr C_i = B \cap H_i.
\end{align}
Since $B$ and $H_i$ are real,  these are real curves.

 We now display all the double curves
we found for every types in the following table.
The one associated by a dagger represents a splitting double conic:
 \begin{center}
\begin{tabular}{c||c|c|c|c}
                & type I & type II & type III & type IV\\ \hline
double conics  & $\ms C_1,  \ms C_2$ &    
$\ms C_1,  \ms C_2, \ms C_3^{\dag}$ &  
$\ms C_1,  \ms C_2, \ms C_3,\ms C_4^{\dag}$  & 
$\ms C_1,  \ms C_2, \ms C_3,\ms C_4, \ms C_5^{\dag}$  \\ \hline
double quartic curves & $\ms C_3,\ms C_4, \ms C_5$      &    
$\ms C_4, \ms C_5$    &  $\ms C_5$      & -\\ \hline
total number & 5     &    5   &   5      & 5\\ 
\end{tabular} 
\end{center}
Note that the number of the double conics is always equal to the 
number of  reducible members of the pencil $|F|$.

\subsection{Quadrics containing 
the double curves}\label{ss:Q}
Next  we investigate  quadrics in $\CP^4$ containing all
the double curves found in the last subsection. 
For this purpose we first make it clear how the five double curves on $B$ intersect each other.

\begin{figure}
\includegraphics{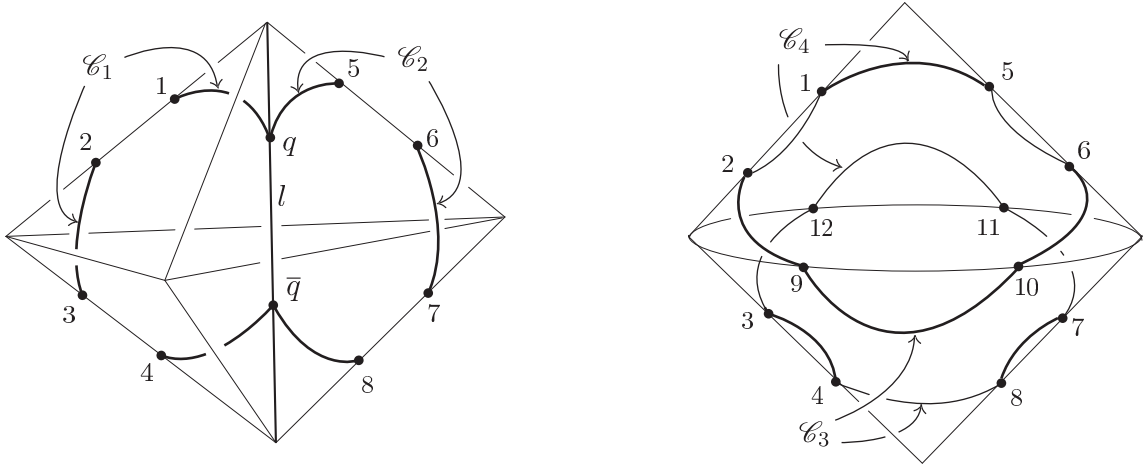}
\caption{The intersection of two double conics (left)
and two double quartic curves (right).
The left picture indicates all coordinate hyperplanes, 2-planes and lines, as well as double conics.
In the right picture the upper and lower halves are the cones $Y\cap H_4$ and $Y\cap H_3$ respectively, on which 
the double quartic curves
$\mathscr C_4$ and $\mathscr C_3$ lie.
The common numbers represent the same point, all of which
are ordinary double points of $B$.}
\label{fig:dc}
\end{figure}

If $\ms C_i$ is a double conic, we have
$\ms C_i = B\cap P_i$ as sets by \eqref{dcrest}.
Since $B\cap l=\{q,\ol q\}$ by \eqref{q} and $l\subset P_i$, it follows that $\ms C_i\cap l
= \{q,\ol q\}$.
Thus  $\ms C_i\cap \ms C_j=\{q,\ol q\}$
for any distinct double conics $\ms C_i$ and 
$\ms C_j$.
(See the left picture in Figure \ref{fig:dc}.)

Next we consider the intersection of two
double quartic curves.
Let $\ms C_i$ be a double quartic curve
and  $H_i$
 the hyperplane which satisfies $B\cap H_i= \ms C_i$ as above.
Then as seen above, we have $\ms C_i\in 
|\ms O_{Y\cap H_i}(2)|$.
If $\ms C_j$ is another double quartic curve
and $H_j$ is a hyperplane satisfying $B\cap H_j = \ms C_j$,
we have
\begin{align}\label{intcq}
\ms C_i\cap \ms C_j = \ms C_i\cap (B\cap H_j)
=(\ms C_i\cap B) \cap H_j
=\ms C_i\cap H_j.
\end{align}
As $\ms C_i$ is a quartic curve, this means that 
the intersection $\ms C_i\cap \ms C_j$ consists
of four points.
Further, since both $\ms C_i$ and $H_j$ are 
real, these four points consist of two real pairs.
(See the right picture in Figure \ref{fig:dc}.)

Finally we consider the intersection of
a double conic with a double quartic curve.
Let $\ms C_i\subset P_i$ be a double conic and $\ms C_j=B\cap H_j$
a double quartic curve.
Then in the same way to \eqref{intcq}, the intersection $\ms C_i\cap \ms C_j$ 
is a hyperplane section $\ms C_i\cap H_j$.
This means that $\ms C_i\cap \ms C_j$ consists of two points,
which are on the  line $ P_i\cap H_j$.
Again since $\ms C_i$ and $H_j$ are real,
$\ms C_i\cap \ms C_j$ is a real pair of points.
(See the two pictures in Figure \ref{fig:dc}.)

Thus for the case of type I,
the total intersection $\{\ms C_i\cap \ms C_j\set 1\le i<j\le 5\}$
consists of 
\begin{itemize}
\item
$\ms C_1\cap \ms C_2=\{q,\ol q\}$.
\item
$\ms C_i\cap \ms C_j,\,\,(i,j) = (3,4), (3,5),(4,5).$
These consist of $ 12$ points.
\item 
$\ms C_i\cap \ms C_j,\,\,(i,j) = (1,3), (1,4),(1,5),
(2,3),(2,4),(2,5).$
These consist of $12$ points.
\end{itemize}
Hence the total intersection consists of $26$ points.

For the case of type II,
the total intersection 
consists of 
\begin{itemize}
\item
$\ms C_i\cap \ms C_j,\,\, (i,j) = (1,2), (1,3), (2,3).$
These consist of just two points $\{q,\ol q\}$.
\item 
$\ms C_4\cap \ms C_5$.
This consists of $4$ points.
\item
$\ms C_i\cap \ms C_j,\,\,(i,j) = (1,4), (1,5),(2,4),
(2,5), (3,4), (3,5).$
These consist of $ 12$ points.
\end{itemize}
Hence the total intersection consists of $18$ points.

For the case of type III,
the total intersection 
consists of 
\begin{itemize}
\item
$\ms C_i\cap \ms C_j,\,\, 1\le i<j\le 4.$
These consist of just two points $\{q,\ol q\}$.
\item
$\ms C_i\cap \ms C_5,\,\,1\le i\le 4$.
These consist of $8$ points.
\end{itemize}
Hence the total intersection consists of $10$ points.

For the case of type IV, as there does not exist
double quartic curves, 
the two points $\{q,\ol q\}$ are exactly the total intersection.

With these information, we are now ready for proving
the existence of the quadric:

\begin{proposition}\label{prop:Q}
Let $Z$ be a twistor space of double solid type on $4\CP^2$.
Then there exists a quadratic hypersurface $Q$ which contains
all the double curves $\ms C_i$, $1\le i\le 5$, and 
which is different from the scroll $Y$.
 \end{proposition}

\noindent
{\em Proof of Proposition \ref{prop:Q} for the cases
of type I and type II.}
For these two types we can prove the assertion by a simple dimension counting.
For the case of type I,
we first consider all quadrics  going through
the following 8 points:
\begin{enumerate}
\item[
(a)] $\ms C_1\cap \ms C_2$
(consisting of 2 points),
\item[
(b)]  $\ms C_1\cap \ms C_3$
(consisting of 2 points),
\item[(c)]  $\ms C_2\cap \ms C_3$
(consisting of 2 points),
\item
[(d)] one of the two points $\ms C_1\cap \ms C_4$,
and \item[
(e)] one of the two points $\ms C_2\cap \ms C_4$.
\end{enumerate}
Since $h^0(\CP^4,\ms O(2))=15$,
these quadrics form at least $(14-8=)\,6$-dimensional linear system.
Let $Q$ be any quadric in this system.
As (a)-(e) include 5 points on the conic $\ms C_1$, the inclusion
$Q\supset \ms C_1$ follows.
By the same reason, the inclusion $Q\supset\ms C_2$ holds.
These trivially mean $Q\supset (\ms C_1\cap \ms C_3)
\cup (\ms C_2\cap \ms C_3)$,
$Q\supset (\ms C_1\cap \ms C_4)
\cup (\ms C_2\cap \ms C_4)$ and
$Q\supset (\ms C_1\cap \ms C_5)
\cup (\ms C_2\cap \ms C_5)$.
We next impose the constraint
$Q\supset \ms C_3\cap \ms C_4$.
Since $\ms C_3\cap \ms C_4$ consists of 4 points as above,
these quadrics still form a linear system which is at least
$(6-4=)\,2$-dimensional.
Then from the first one of the above three inclusions,
$Q$ contains 8 points on $\ms C_3$.
Hence as $(\ms C_3,\ms C_3)_{Y\cap H_3} = 8$,
the inclusion $Q\supset \ms C_3$ follows.
By the same reason 
the constraint $Q\supset \ms C_3\cap \ms C_4$ means
$Q\supset \ms C_4$ and $Q\supset \ms C_5$.
Thus $Q$  contains all the five double conics.
Further we can choose $Q$ which is different 
from $Y$ since these quadrics form
a 2-dimensional family at least.
Thus we obtain the assertion in the case of type I.

For the case of type II, 
we consider quadrics in $\CP^4$ which 
go through all the following 12 points:
\begin{enumerate}
\item[(a)] $\mathscr C_1\cap\mathscr C_2$ (consisting of 2 points),
\item[(b)] $\mathscr C_4\cap \mathscr C_5$ 
(consisting of 4 points),
\item[(c)]
$\mathscr C_1\cap \mathscr C_4$,
and one of the two points $\mathscr C_1\cap \mathscr C_5$
(consisting of 3 points),
\item[(d)] $\mathscr C_2\cap \mathscr C_4$
(consisting of 2 points),
\item[(e)] one of the 2 points $\mathscr C_2\cap \mathscr C_5$
(consisting of 1 point).
\end{enumerate}
Again by dimension counting
these quadrics form a $2$-dimensional linear system at least.
We can suppose $Q\neq Y$.
We show that any quadric $Q$ in this family  contains all the five double curves.
In fact, from (a)--(c), $Q$ contains 5 points on the conic $\mathscr C_1$,
and hence $ Q\supset \ms C_1$.
Further from (b)--(d),
$Q$ passes through 8 points on $\mathscr C_4$,
which means $ Q\supset\ms C_4$.
Furthermore from  (a), (d) and (e),
$Q$ goes through 5 points on $\mathscr C_2$, which means
$ Q\supset\ms C_2$.
This implies that $Q$ passes through 8 points on $\mathscr C_5$,
and therefore $ Q\supset \ms C_5$.
These mean that $Q$ contains 6 points on $\mathscr C_3$,
meaning $ Q\supset\ms C_3$.
Thus the quadric $Q$ contains all the five double curves.
Hence we have shown the assertion 
for the case of type II.
\proofend

\vsp
For the cases of type III and type IV,
a similar dimension counting does not work for proving Proposition \ref{prop:Q}, and we need to
derive a constraint for the double conics by using
the fact that all the double curves are lying on 
a common quartic hypersurface.
Recall that for any double conic $\ms C_i$ we have
$\ms C_i\cap l = \{q,\ol q\}$ as 
shown at the beginning of this subsection.
In particular, the tangent lines $T_q\ms C_i$ 
and $T_{\ol q}\ms C_i$ of the conic
are different from the line $l$ itself.

\begin{lemma}\label{lemma:common}
Suppose that $Z$ is of type III or type IV.
Then
there exists a 3-dimensional linear subspace $V\subset T_q\CP^4=\CC^4$ 
which satisfies 
$T_q\ms C_i\subset V$ for any double conic $\ms C_i$ on $B$.
\end{lemma}
 
\proof
Let $\nu_1:\tilde Y\to Y$ be the blowup of $Y$ along the line $l$, 
and $\Sigma_1$ the exceptional divisor.
This $\tilde Y$ is smooth and biholomorphic to the 
$\CP^2$-bundle
$\mathbb P(\ms O(2)^{\oplus 2}\oplus \ms O)$ over 
the conic $\Lmd$.
We have $\Sigma_1\simeq\mathbb F_0$.
As $q,\ol q\in l$, the pre-images $\nu_1^{-1}(q)$ and $\nu_1^{-1}(\ol q)$ are
isomorphic to $\CP^1$.
Let $\nu_2:\tilde Y'\to \tilde Y$ be the blowup of $\tilde Y$
along
$\nu_1^{-1}(q)\sqcup\nu_1^{-1}(\ol q)$,
and $\Sigma_2$ and $\ol\Sigma_2$ the exceptional divisors
over  $\nu_1^{-1}(q)$ and $\nu_1^{-1}(\ol q)$ respectively.
We readily have $\Sigma_2\simeq\ol\Sigma_2\simeq\mathbb F_{2}$,
and that the co-normal bundle $\ms O_{\Sigma_2}(-\Sigma_2)$ is
represented by a $(+2)$-section of the ruling
$\Sigma_{2}\to\CP^1$.

Let $T_qY$ be the tangent cone of $Y$ at $q$.
For any  $v\in T_qY$ 
satisfying $v\not\subset l$,
we can consider a natural lift $\tilde v$ to the blowup $\tilde Y$,
which is an element of the tangent space $T_{[v]}\tilde Y$,
where $[v]\in \Sigma_1$ is the point
represented by the tangent vector $v$.
The vector $\tilde v$ can further be lifted to 
a tangent vector $\tilde v'$ on $\tilde Y'$.
We denote $[\tilde v']\in \Sigma_2$ 
for the initial point of $\tilde Y'$.
Under these notations, we now show that, 
for any collection $v_1,v_2,\cdots, v_k$ of 
elements in $T_qY$ satisfying $v_i\not\in l$,
these are linearly independent 
(as elements in $T_q\CP^4$) if and only if
the points
$[\tilde v'_1],[\tilde v'_2],\cdots,
[\tilde v'_k]\in\Sigma_2$ belong to a common $(+2)$-section of the ruling $\Sigma_2\to\CP^1$.

For this let $X$ be an intersection of $Y$ with a
hyperplane which goes through the point $q$ but 
which does not contain the line $l$.
$X$ is a quadratic cone in the hyperplane whose vertex is $q$.
We write $\tilde X:=\nu_1^{-1}(X)$.
Clearly $\nu_1^{-1}(q)\subset \tilde X$,
and $\tilde X\in |\nu_1^*\ms O_Y(1)|$.
Let $\tilde X'$ be the strict transform of $\tilde X$ into $\tilde Y'$.
Then $\tilde X'\in |\nu_2^*\nu_1^*\ms O_Y(1) - \Sigma_2|$.
Since the line bundle $\nu_1^*\ms O_Y(1)$ is trivial
over $\nu_1^{-1}(q)$,
we have $\tilde X'|_{\Sigma_2}\in |-\Sigma_2|_{\Sigma_2}|$.
Thus we have obtained the assignment 
$X\mapsto \tilde X'|_{\Sigma_2}$.
As above the last system is linearly equivalent
to the class of  $(+2)$-sections of the ruling $\Sigma_2 \simeq\mathbb F_2
\to \CP^1$.
Also we can show that the above assignment gives a bijection from
 the space of the hyperplane sections satisfying the above properties
(i.e.\,$X\ni p$ and $X\not\supset l$)
to the complete linear system of $(+2)$-sections.
Therefore, the above criterion for the linearly dependence follows.

Thus in order to prove the assertion of the lemma,
it is enough to show that all the points on $\Sigma_2$ 
determined from the tangent line
$T_q\ms C_i$ of the double conics belong to a common $(+2)$-section on $\Sigma_2$.
To show this recall that all the double conics are lying on the branch divisor $B$,
and $B$ is an intersection of $Y$ with a quartic hypersurface in $\CP^4$.
 In particular we have $B\in |\ms O_Y(4)|$.
 Hence the inverse image $\tilde B:=\nu_1^{-1}(B)$ belongs to $ |\nu_1^*\ms O_Y(4)|$.
Clearly $\nu_1^{-1}(q)\sqcup\nu_1^{-1}(\ol q)\subset\tilde B$, 
 and further, since we are considering the case of type III or type IV,
the divisor  $\tilde B$ has $A_2$- or $A_3$-singularities along
 $\nu_1^{-1}(q)\sqcup\nu_1^{-1}(\ol q)$ according to the type respectively.
Hence the strict transform $\tilde B'$ of $\tilde B$ into $\tilde Y'$ belongs
to the system $|\nu_2^*\nu_1^*\ms O_Y(4) - 2\Sigma_2 - 2\ol\Sigma_2 |$,
and moreover $\tilde B'|_{\Sigma_2}$ is a non-reduced curve of
multiplicity two.
Since $\tilde B'|_{\Sigma_2}\in |-2\Sigma_2|_{\Sigma_2}|$,
this implies that the set theoretical intersection $\tilde B'\cap \Sigma_2$ is a $(+2)$-section of $\Sigma_2\to\CP^1$.
Therefore since the strict transforms of the double conics are clearly contained
in $\tilde B'$, we obtain the  claim of the lemma.
\proofend

\vsp
\noindent
{\em Proof of Proposition \ref{prop:Q} for the cases
of type III and type IV.}
Let $V$ be the hyperplane in Lemma \ref{lemma:common},
and $\ol V$ the image of $V$ under the real structure on
$\CP^4$.
Then since any double conic $\ms C_i$ is real,
the hyperplane $\ol V$ satisfies $T_{\ol q}\ms C_i\subset
\ol V$ for any double conic $\ms C_i$.

For the case of type III,
 we consider the set of all quadrics
which satisfy the following three conditions:
(a) they go through the two points $q$ and $\ol q$,
(b) they are non-singular at $q$ and $\ol q$,
and the tangent space at these points are
$V$ and $\ol V$ respectively,
(c) they go trough the 8 points $\{\ms C_i\cap\ms C_5
\set 1\le i\le 4\}$.
By dimension counting, these quadrics form a family 
whose dimension is at least $14- (2+2+8)=2$.
We show that any member $Q$ of this family contains 
the double curves $\ms C_i$ for any $1\le i\le 5$.
$\ms C_5\subset Q$ is obvious since $Q$ already goes through
the 8 points (c).
For any $i\in \{1,2,3,4\}$, 
 the intersection $ P_i\cap Q$ is tangent to
the double conic $\ms C_i$ at the points $q$
and $\ol q$  from the condition (b).
Further, since $Q\cap  P_i$ goes through 
the two points $\ms C_i\cap \ms C_5$ by (c),
the conic $Q\cap P_i$ goes through
in effect 6 points on $\ms C_i$, so
we obtain $\ms C_i=Q\cap P_i$.
Hence
the inclusion $Q\supset \ms C_i$ follows.

For the case of type IV,
we consider  the set of all quadrics 
which satisfy the above conditions (a), (b), and
the following (c') instead of (c);
(c') for each $1\le i\le 5$, they go through 
a point on $\ms C_i$ which is different from $q$ and $\ol q$.
These quadrics form a family whose dimension is at least
$14-(2+2+5)=5$.
Then by considering the intersection $P_i\cap
Q$, we again obtain that $\ms C_i\subset Q$ for any $1\le
i\le 5$.
\proofend

\subsection{Defining equation of the branch divisor}
\label{ss:de}
We are ready for determining defining equation
of the  branch divisor of the anticanonical map.
As before let $\pi:\CP^4 \to \CP^2$ be a linear projection,
$\Lambda$ a conic in the target plane, and 
$Y\subset\CP^4$ the scroll over $\Lambda$. 
We know that $Y$ is the image of the anticanonical map of the twistor spaces.
We fix any homogeneous coordinates 
$(z_0,z_1,z_2)$ on the above plane
such that each $z_i$ is real and that the conic $\Lambda$ is defined by 
the equation $z_0^2 = z_1z_2$.

\begin{theorem}\label{thm:main}
Let $Z$ be  a Moishezon twistor space on $4\CP^2$ which is of double solid type,
and $\Phi:Z\to Y\subset\CP^4$ the anticanonical map.
Then the branch divisor of \,$\Phi$ is an intersection of the scroll \,$Y$ with a 
quartic hypersurface
defined by the following equation:
\begin{align}
z_0z_3 z_4  f (z_0,z_1,z_2,z_3,z_4) = Q(z_0,z_1,z_2,z_3,z_4)^2&
\quad{\text{for the case of type I,}}
\label{eqn:BI}\\
z_0(z_0-z_1)z_3z_4  = Q(z_0,z_1,z_2,z_3,z_4)^2&
\quad{\text{for the case of type II,}}
\label{eqn:BII}\\
z_0(z_0-z_1) (z_0-a z_1) z_4 = 
Q(z_0,z_1,z_2,z_3,z_4)^2
& \quad{\text{for the case of type III,}}
\label{eqn:BIII}\\
z_0(z_0-z_1)(z_0-a_1 z_1)(z_0- a_2 z_1) = 
Q(z_0,z_1,z_2,z_3,z_4)^2
& \quad{\text{for the case of type IV,}}
\label{eqn:BIV}
\end{align}
where $ f (z_0,z_1,z_2,z_3,z_4)$ 
and $Q (z_0,z_1,z_2,z_3,z_4)$ are linear 
and quadratic polynomials respectively,
and $a$, $a_1$ and $a_2$ are
real numbers satisfying 
$\{a,a_1,a_2\}\cap\{0,1\}=\emptyset$ and $a_1\neq a_2$.
Further, in the cases of type II, III and IV,
the conic defined by the equation   $z_0=z_1=Q(z_0,z_1,z_2,z_3,z_4)=0$ is reducible.
\end{theorem}

\begin{remark}{\em
Before proceeding to a proof,
we give a comment about the equations
\eqref{eqn:BI}--\eqref{eqn:BIV}.
All the equations are of the form
\begin{align}\label{gform}
f_1f_2f_3f_4= Q^2,
\quad \deg f_i = 1,\,\deg Q = 2.
\end{align}
In particular, they are of the same form
as the case of $3\CP^2$ obtained 
by Poon \cite{P92} and Kreussler-Kurke
\cite{KK92},
except the number of variables.
According to each of the three degenerations
I $\to$ II $\to$ III $\to$  IV,
exactly one of the four linear polynomials degenerates from 
those with 5 variables $z_0,z_1,z_2,z_3,z_4$ to those with 3 variables 
$z_0,z_1,z_2$; 
geometrically this means that the hyperplane (in $\CP^4$) defined by the linear polynomial
degenerates from a general one to those which contains the
ridge $l$ of the scroll $Y$.
}
\end{remark} 

\vsp
\noindent
{\em Proof of Theorem \ref{thm:main}.}
First we choose appropriate homogeneous coordinates 
$(z_0,z_1,z_2,z_3,z_4)$ on $\CP^4$ for each of the four types.
By applying a projective transformation of 
the plane $\CP^2$ which preserves the conic $\Lmd$, we can suppose that 
the point $\lmd_i$ on the conic $\Lmd$ corresponding to reducible members
of $|F|$ is placed as in the following table:

 \begin{center}
\begin{tabular}{c||c|c|c|c|c}
& $\lmd_1$ & $\lmd_2$ 
& $\lmd_3$ & $\lmd_4$ & $\lmd_5$\\
\hline
type I & $(0,1,0)$ & (0,0,1) & - & - & -\\
\hline
type II & $(0,1,0)$ & $(1,1,1)$ & $(0,0,1)$ & - & -
  \\ \hline
type III & $(0,1,0)$ & $(1,1,1)$ 
& $(a,1,a^2)$ & $(0,0,1)$ & - 
\\ \hline
type IV  & $(0,1,0)$ & $(1,1,1)$ & $(a_1,1,a_1^2)$ & $(a_2,1,a_2^2)$ &  $(0,0,1)$ 
\end{tabular} 
\end{center}
Here, in the case of type II,
$a$ is a real number satisfying $a\not\in\{0,1\}$,
and in the case of type III,
$a_1$ and $a_2$ are distinct real numbers satisfying $\{a_1,a_2\}\cap\{0,1\}=\emptyset$.
Then we have
\begin{align}
\Lmd\cap \{z_0 = 0\} = \{\lmd_1,\lmd_2\}
&{\text{ for the case of type I}},\label{I1}\\
\Lmd\cap \{z_0 (z_0-z_1)= 0\}  = \{\lmd_1,\lmd_2,\lmd_3\}
&{\text{ for the case of type II}},\label{II1}\\
\Lmd\cap \{z_0 (z_0-z_1) 
(z_0- a z_1)= 0\} = \{\lmd_1,\lmd_2,\lmd_3,\lmd_4\}
&{\text{ for the case of type III}},\\
\Lmd\cap \{z_0 (z_0-z_1) (z_0- a_1 z_1)
(z_0- a_2 z_1) = 0\} = \{\lmd_1,\lmd_2,\lmd_3,\lmd_4,\lmd_5\}
&{\text{ for the case of type IV.}}
\end{align}
Here we note that for the case of type II, III and IV 
the final point $(0,0,1)$ is included as a multiple point,
where the multiplicity is 2, 3 and 4 
for the cases of type II, III and IV respectively.

Next we choose the remaining two coordinates $z_3$ 
and $z_4$, as well as 
a linear polynomial $f\in\CC[z_0,z_1,z_2,z_3,z_4]$ in such a way that
\begin{itemize}
\item
for the case of type I: $(z_3)= H_3,\,(z_4) = H_4$
and $(f) = H_5$,
\item
for the case of type II: $(z_3)= H_3,\,(z_4) = H_4$,
\item
for the case of type III: $(z_3)= H_3$
\end{itemize}
are satisfied, where $H_i$-s are the 
hyperplanes in $\CP^4$ found in Section \ref{ss:double}, which cut out the double quartic
curves $\ms C_i$-s.

Under these homogeneous coordinates, 	
for an algebraic subset $X\subset\CP^4$,
we denote by $I_X\subset\mathbb C[z_0,\cdots,z_4]$ for
the homogeneous ideal of $X$.
Let $F=F(z_0,\cdots,z_4)$ be a  quartic
polynomial which defines a branch divisor 
of the anticanonical map $\Phi$ as 
the intersection with the scroll $Y$.
Obviously $F$ is defined only up to polynomials in the ideal
$I_{Y}\subset\mathbb C[z_0,\cdots,z_4]$
of the scroll.
Let $Q$ be a defining polynomial of a quadric
in $\CP^4$ containing all the double conics,
whose existence was proved in 
Proposition \ref{prop:Q}.
Let $\ms C_i\subset P_i$ be a double curve
on $B$.
Then as $(F|_{ P_i}) = 2\mathscr C_i
= (Q^2|_{ P_i})$
as divisors on the plane $ P_i$,
there exists a constant $c_i$ such that $F-c_iQ^2 \in
I_{ P_i}\subset\CC[z_0,\cdots,z_4]$.
If $c_i\neq c_j$ for some $i\neq j$,
we obtain $Q^2\in I_{ P_i} + I_{ P_j}$.
Further the last ideal is readily seen to be equal to
$ I_{ P_i\cap P_j}$, and therefore
equals to $I_l=(z_0,z_1,z_2)\subset\CC[z_0,\cdots,z_4]$.
Hence $Q\in (z_0,z_1,z_2)$.
But this means that the divisor  $(Q|_{ P_i})$
contains $l$, which contradicts
$\ms C_i\not\subset l$ (and $\ms C_i\subset P_i)$.
Therefore $c_i=c_j$ for any double conics
$\mathscr C_i$ and $\mathscr C_j$.

Next let $\ms C_i$ be any double quartic curve.
Then since $(F|_{H_i\cap Y}) = 2\mathscr C_i
= (Q^2|_{H_i\cap Y})$,
there exists a constant $c_i$ such that 
$F - c_iQ^2\in I_{H_i\cap Y}=(z_i) + I_Y$.
So taking a difference with $F-c_1Q^2\in I_{ P_1}$,
we obtain that $(c_1-c_i)Q^2 \in (z_i) + I_Y + 
I_{ P_1}$.
But since $ P_1\subset Y$, we have
$I_{ P_1}\supset I_Y$,
and therefore 
$(c_1-c_i)Q^2 \in (z_i) + I_{ P_1}$.
Hence if $c_1\neq c_i$ we have $Q^2\in (z_i) + I_{ P_1}$, which means $Q^2|_{ P_1}\in (z_i|_{ P_1})$. 
Since $i>2$, this means that the divisor $(Q^2)|_{ P_1}$ contains a line $(z_i)$ on 
the plane $ P_1$ as an irreducible component,
which  again contradicts the irreducibility of $\mathscr C_1$.
Therefore we have $c_1=c_i$ for any double 
quartic curve $\mathscr C_i$.
By rescaling we can suppose $c_i=1$ 
for any $1\le i\le 5$.
Thus  we have
\begin{align}
F - Q^2 &\in I_{ P_1}\cap 
I_{ P_2}
\cap \big((z_3) + I_Y\big)
\cap  \big((z_4) + I_Y\big)
\cap \big((f) + I_Y\big)
&{\text{for the case of type I}},\label{I2}
\\
F - Q^2 &\in I_{ P_1}\cap 
I_{ P_2}\cap 
I_{ P_3}\cap \big((z_3) + I_Y\big)
\cap  ((z_4) + I_Y)
&{\text{for the case of type II}},
\label{gaeof9}
\\
F - Q^2 &\in I_{ P_1}\cap 
I_{ P_2}\cap 
I_{ P_3}\cap 
I_{ P_4}
\cap  \big((z_4) + I_Y\big)
&{\text{for the case of type III}},\\
F - Q^2 &\in I_{ P_1}\cap 
I_{ P_2}\cap 
I_{ P_3}\cap 
I_{ P_4}\cap I_{P_5}
&{\text{for the case of type IV}}.
\end{align}

Now for the case of type I, by \eqref{I1}, we have
$$
I_{P_1}\cap I_{P_2} =
I_{P_1\cup P_2}
= I_{\{z_0=0\}\cap Y}
= (z_0) + I_Y.
$$
So from \eqref{I2} we have
\begin{align}\label{f-q^2I}
F- Q^2 \in \bigcap_{i=0,3,4} \big((z_i)+ I_Y\big)
\cap \big( (f) + I_Y \big).
\end{align}
But it is elementary to see that 
the right hand side of \eqref{f-q^2I} is equal to 
the ideal $(z_0z_3z_4f) + I_Y$.
This directly implies that the equation $F$ of the 
quartic hypersurface in $\CP^4$ is exactly as in \eqref{eqn:BI}
in Theorem \ref{thm:main}.

For the case of type II, in a similar way, using \eqref{II1}
we obtain 
\begin{align}\label{f-q^2II}
F- Q^2 \in \bigcap_{i=0,3,4} \big((z_i)+ I_Y\big)
\cap I_{P_2}.
\end{align}
But clearly $I_{P_2}$ can be replaced by
$I_{P_1}\cap I_{P_2}$, which is equal to
$(z_0-z_1)+I_Y$.
Then we can see that the right hand side of \eqref{f-q^2II} is equal to
the ideal $(z_0(z_0-z_1)z_3z_4) + I_Y$.
This again implies \eqref{eqn:BII} in Theorem \ref{thm:main}.

The defining equations for the cases of type III and type IV
can be derived in a similar way, by disposing
$I_{P_1}$ and replacing the remaining  
$I_{P_i}$ by $I_{P_i}\cap I_{P_1}$.

It remains to show the final statement concerning reducibility of a conic.
For this suppose that $Z$ is not of type I and let $\ms C_i$ be the (unique) splitting 
double conic (see  \eqref{sdcII}--\eqref{sdcIV}).
Then from the table at the beginning of this proof,
the equation of the plane $P_i$ is given by $z_0=z_1=0$.
Therefore from \eqref{eqn:BII}--\eqref{eqn:BIV} we have
$\ms C_i=\{z_0=z_1=Q=0\}$.
Since this is reducible, we obtain the
assertion.
\proofend

\end{document}